\documentclass[a4paper,dvipsnames,fleqn]{article}

\usepackage{geometry}
\usepackage{amsmath}
\usepackage{amssymb}
\usepackage{amsthm}
\usepackage{bm}
\usepackage{xcolor}
\definecolor{darkred}{rgb}{0.5,0.,0.}
\usepackage[colorlinks=true,allcolors=darkred]{hyperref}
\usepackage[capitalise]{cleveref}
	\crefformat{equation}{\textup{#2(#1)#3}}
\usepackage[british]{babel}
\usepackage[UKenglish]{isodate}
\usepackage{dsfont}
\usepackage{enumitem}
\usepackage[OT2,OT1]{fontenc}
\usepackage[cal=cm,scr=boondoxo]{mathalfa}
\usepackage{pifont}
\usepackage{stmaryrd}
\usepackage{textcomp}
\usepackage{tikz}
	\usetikzlibrary{arrows}
	\usetikzlibrary{patterns}
\usepackage{tikz-cd}
\usepackage[textwidth=3cm,colorinlistoftodos]{todonotes}
\usepackage{xfrac}

\numberwithin{equation}{section}			
\swapnumbers						

\newcommand\cyr{
\renewcommand\rmdefault{wncyr}%
\renewcommand\sfdefault{wncyss}%
\renewcommand\encodingdefault{OT2}%
\normalfont
\selectfont}
\DeclareTextFontCommand{\textcyr}{\cyr}

\newcounter{Enum}				
\newenvironment{Enumerate}{\begin{enumerate}[label={\rm({\roman*})}]}{\end{enumerate}}

\newcommand{\descriptionlabelsave}{}		
\newenvironment{Itemize}{%
	\renewcommand{\descriptionlabelsave}{\descriptionlabel}\renewcommand{\descriptionlabel}{$\triangleright$}%
	\begin{description}[leftmargin=15pt,itemindent=-5.2pt]}{%
	\end{description}\renewcommand{\descriptionlabel}{\descriptionlabelsave}}

\newcounter{StepsCount}				
\newenvironment{Elist}{%
	\begin{list}{\ding{\value{StepsCount}}}{\usecounter{StepsCount} \leftmargin=0pt \labelwidth=12pt \itemindent=\labelwidth%
	\itemsep=5pt\listparindent=\parindent} \setcounter{StepsCount}{191}}{\end{list}}
\newcounter{StepsRefCount}

\newenvironment{Ilist}{
	\begin{list}{$\triangleright$}{\leftmargin=0pt \labelwidth=11pt \itemindent=\labelwidth%
	\itemsep=5pt\listparindent=\parindent}}{\end{list}}

\theoremstyle{plain}
	\newtheorem{lemma}{Lemma}[section]
	\newtheorem{proposition}[lemma]{Proposition}
	\newtheorem{theorem}[lemma]{Theorem}
	\newtheorem{corollary}[lemma]{Corollary}
	\newcommand{\GenericTheoremName}{}
\theoremstyle{definition}
	\newtheorem{definition}[lemma]{Definition}
	\newcommand{\GenericDefinitionName}{}\newtheorem{genericdefinition}[lemma]{\GenericDefinitionName}
	\newtheorem{notation}[lemma]{Notation}
\theoremstyle{remark}
	\newtheorem{remark}[lemma]{Remark}
	\newtheorem{example}[lemma]{Example}
	\newcommand{\GenericRemarkName}{}\newtheorem{genericremark}[lemma]{\GenericRemarkName}
\newenvironment{Definition}{\begin{definition}}{\hfill$\blacktriangleleft$\end{definition}}

\newenvironment{Notation}{\begin{notation}}{\par\noindent\hfill$\blacktriangleleft$\end{notation}}
\newenvironment{Remark}{\begin{remark}}{\hfill$\vartriangleleft$\end{remark}}
\newenvironment{Example}{\begin{example}}{\hfill$\vartriangleleft$\end{example}}

\newcommand{\mc}[1]{{\mathcal{#1}}}			
\newcommand{\ms}[1]{{\mathscr{#1}}}			
\newcommand{\mf}[1]{{\mathfrak{#1}}}			
\newcommand{\bb}[1]{{\mathbb{#1}}}			
\newcommand{\mr}{\mathring}				
\newcommand{\wt}{\widetilde}				

\DeclareMathOperator{\RE}{Re}				
\renewcommand{\Re}{\RE}
\DeclareMathOperator{\IM}{Im}				
\renewcommand{\Im}{\IM}

\DeclareMathOperator{\dom}{dom}
\DeclareMathOperator{\domr}{dom_r}
\DeclareMathOperator{\tr}{tr}
\DeclareMathOperator{\osc}{osc}

\DeclareMathOperator{\ran}{ran}
\DeclareMathOperator{\essinf}{ess\,inf}
\DeclareMathOperator{\esssup}{ess\,sup}

\newcommand{\Side}[1]{\hfill{#1}\kern10pt}		
\newcommand{\FD}[5]{
	\DF\left\{\begin{array}{rcl}{#1}&\to &{#2}\\[#3pt] {#4}&\mapsto &{#5}\end{array}\right.}

\newcommand{\smmatrix}[4]{\Bigl(			
\begin{smallmatrix}
\hspace*{-0.2ex} #1 \hspace*{0.2ex} & \hspace*{0.2ex} #2 \hspace*{-0.2ex}
\\[0.5ex]
\hspace*{-0.2ex} #3 \hspace*{0.2ex} & \hspace*{0.2ex} #4 \hspace*{-0.2ex}
\end{smallmatrix}
\Bigr)}

\newcommand{\Dummy}{\text{\textvisiblespace\kern1pt}}	
\newcommand{\Smallo}{{\rm o}}				
\newcommand{\BigO}{{\rm O}}				

\DeclareMathOperator{\Id}{id}				
\DeclareMathOperator{\Span}{span}			

\newcommand{\DS}{\mid\mkern3mu}				
\newcommand{\DQ}{\mkern6mu}				
\newcommand{\DP}{\kern2pt{\mathrel{\mathop:}\kern5pt}}	
\newcommand{\DF}{\colon}				
\newcommand{\DE}{\mathrel{\mathop:}=}			
\newcommand{\ED}{=\mathrel{\mathop:}}			
\newcommand{\DI}{\mathrel{\mathop:}\Leftrightarrow}	
\newcommand{\D}{\mathrm{d}}			
\newcommand{\DD}{\mkern4mu\mathrm{d}}			
\newcommand{\CAS}{&\text{if}\ }				

\DeclareMathOperator{\Ind}{ind}			

\begin{document}

\begin{flushleft}
	{\Large\bf Eigenvalue distribution of canonical systems:\\[2mm] trace class and sparse spectrum}
	\\[5mm]
	\textsc{
	Matthias Langer
	\,\ $\ast$\,\
	Jakob Reiffenstein
	\,\ $\ast$\,\
	Harald Woracek
		\hspace*{-14pt}
		\renewcommand{\thefootnote}{\fnsymbol{footnote}}
		\setcounter{footnote}{2}
		\footnote{This work was supported by the project P~30715-N35 of the
			Austrian Science Fund (FWF).
			The second author was supported by the Sverker Lerheden Foundation.
			The second and third authors were supported by the project I~4600 of the Austrian
			Science Fund (FWF).}
		\renewcommand{\thefootnote}{\arabic{footnote}}
		\setcounter{footnote}{0}
	}
	\\[6mm]
\end{flushleft}
	{\small
	\textbf{Abstract:}
		In this paper we consider two-dimensional canonical systems with discrete spectrum and study their eigenvalue
		densities.  We develop a formula that determines the Stieltjes transform of the eigenvalue counting function up
		to universal multiplicative constants.  An explicit criterion is given for the resolvents of the model operator
		to belong to a Schatten--von~Neumann class with index $p \in (0,2)$, thus giving an answer to the long-standing
		question which canonical systems have trace class resolvents.

		For canonical systems with two limit circle endpoints we develop an algorithm for determining the
		growth of the monodromy matrix up to a small error.  Moreover, we present examples to illustrate our results,
		show their sharpness and prove an inverse result giving explicit formulae.
	}
	\\
\begin{flushleft}
	{\small
	\textbf{AMS MSC 2020:} 34L15, 37J99, 30D15
	\\
	\textbf{Keywords:} canonical system, discrete spectrum, density of eigenvalues, inverse spectral theorem, 
	regularly varying function
	}
\end{flushleft}


\pagenumbering{arabic}
\setcounter{page}{1}
\section{Introduction}

A two-dimensional \emph{canonical system} is a differential equation of the form 
\begin{equation}
\label{X92}
	y'(t)=zJH(t)y(t)
\end{equation}
on some interval $[a,b)$, where $H$ is a $2\times 2$-matrix-valued function which is locally integrable and takes non-negative 
matrices with real entries as values, $z\in\bb C$, and $J=\smmatrix 0{-1}10$. 
Equations of this form originate from Hamiltonian mechanics with
one-dimensional phase space, and the function $H$ is called the \emph{Hamiltonian} of the system.
Canonical systems provide a unifying treatment for various equations of mathematical physics and pure mathematics, 
such as Dirac systems \cite{sakhnovich:2002}, Schr\"odinger equations \cite{remling:2002}, 
Krein strings \cite{kaltenbaeck.winkler.woracek:bimmel}, Jacobi matrices \cite{kac:1999}, 
and generalised indefinite strings \cite{eckhardt.kostenko:2016}. 
Every self-adjoint operator with simple spectrum can be realised in a very specific way as the differential operator 
induced by some canonical system (and appropriate boundary conditions).
Cornerstones of the theory of two-dimensional canonical systems were established by M.G.~Krein \cite{gohberg.krein:1970} and 
L.~de~Branges \cite{debranges:1968}.  Some recent references on the topic, including the construction of an operator model
for the equation \eqref{X92}, are \cite{romanov:1408.6022v1,remling:2018,behrndt.hassi.snoo:2020}.

Given a Hamiltonian $H$, each two self-adjoint realisations of \eqref{X92} are at most two-dimensional perturbations 
of each other. 
Denote by $\sigma_H$ the spectrum of one of them, and consider the situation when $\sigma_H$ is
discrete. Recall at this point that discreteness of $\sigma_H$ is characterised by an explicit and simple condition on 
$H$ in \cite{romanov.woracek:ideal}.  We address the question of how to determine the asymptotic density of $\sigma_H$ from $H$.
Here we think of the density of a sequence having no finite accumulation point in ways familiar from complex analysis. 

Let us start with explaining a simple instance of our results.  Consider the following question:
\begin{equation}
\label{X93}
	\text{\emph{Given $p>0$, does the series}}
	\sum_{\lambda\in\sigma_H\setminus\{0\}}\frac{1}{|\lambda|^p}
	\ \ \text{\emph{converge?}}
\end{equation}
For $p>1$ this question is answered in \cite{romanov.woracek:ideal}.  The formulae characterising convergence in this case 
are explicit and simple; they are of a form very much related to the discreteness criterion 
and depend only on the behaviour of the diagonal entries of $H$ at the right endpoint $b$ (for details see \cref{X137}).
When $p\le 1$, the situation is much more complicated. 
In the setting of Krein strings, which translates to the case when $H$ is a.e.\ a diagonal matrix, 
a criterion for convergence of the series in \eqref{X93} is given in \cite{kac:1986}.  
The formulae are explicit but cumbersome. 
For a Hamiltonian of arbitrary form and $p\le 1$, characterising convergence of the series in \eqref{X93} 
remained an open problem.  As a consequence of our present results, we can solve this problem. 
In particular, we obtain a characterisation of the case when the resolvents
of the model operator of \eqref{X92} are of trace class.
The criterion itself is quite involved, but explicit in terms of $H$. 
It should be noted that the characterisation also depends on the off-diagonal entries of $H$ and, moreover,
on properties of $H$ on the entire interval $[a,b)$ and not just the asymptotic behaviour at~$b$.

Our theorems go far beyond characterising convergence of \eqref{X93}. 
Denote by $n_H$ the counting function of the spectrum:
\[
	n_H(r) \DE \#\bigl\{\lambda\in\sigma_H\DS|\lambda|<r\bigr\}.
\]
In \cref{X79}, which is the first main result of the paper, we give a formula that determines 
the Stieltjes transform of $\frac{n_H(\sqrt t)}{t}$ up to universal multiplicative constants: 
there exist $C_-,C_+>0$ such that, for every Hamiltonian $H$
and $r>\bigl(\det\int_a^b H(t)\DD t\bigr)^{-\frac12}$, 
\begin{equation}
\label{X32}
	\frac{C_-}{r^2}\int_a^b K_H(t;r)\DD t \le \int_0^\infty\frac{1}{t+r^2}\cdot\frac{n_H(\sqrt t)}{t}\DD t
	\le \frac{C_+}{r^2}\int_a^b K_H(t;r)\DD t;
\end{equation}
the integrand $K_H(t;r)$ is a function defined explicitly in terms of $H$; see \cref{X72}. 

Applying appropriate Tauberian and Abelian theorems one can pass from the Stieltjes transform in 
\eqref{X32} to integrability or boundedness properties (with respect to regularly varying comparison functions) 
of the function $n_H$; see \cref{X6,X7}.
For example, returning to the question in \eqref{X93}, we obtain the equivalence
\begin{align}\label{X122}
	\sum_{\lambda\in\sigma_H\setminus\{0\}}\frac{1}{|\lambda|^p} \quad\Longleftrightarrow\quad
	\int_1^\infty\frac{1}{r^{p+1}}\int_a^b K_H(t;r)\DD t\DD r < \infty
\end{align}
for any $p\in(0,2)$.

For the proof of \cref{X79} we develop an analytic approach to asymptotic density.  It is worth comparing it to the
methods used in \cite{romanov.woracek:ideal} or \cite{kac:1986} since it is different from both. 
The approach in \cite{romanov.woracek:ideal} is operator-theoretic: growth properties of $n_H$ are encoded as 
membership of resolvents in certain operator ideals, and proofs exploit the theory of such ideals.
The argument breaks down slightly above the trace class ideal $\mf S_1$, and this is necessarily so:
the results of \cite{romanov.woracek:ideal} do not hold for $\mf S_1$.
The theorem from \cite{kac:1986} is a closer relative to our \cref{X79}.
It can be interpreted as being based on a differential inequality for an associated Riccati 
equation, and relies on the fact that the system \eqref{X93} can easily be rewritten as a scalar second-order equation 
when $H$ is diagonal.  It breaks down when off-diagonal entries are non-zero on some set of positive measure. 

The central piece of our present argument is a formula that determines the logarithm 
of the modulus of an entry of the fundamental solution up to universal multiplicative constants; 
this is given in \cref{X1}. 
In the first step of the proof of the latter theorem we use the canonical
differential equation \eqref{X92}.
Then it is crucial to realise that one can apply the estimates for Weyl coefficients of 
canonical systems from the recent work \cite{reiffenstein:imq,langer.pruckner.woracek:heniest} 
in a rather tricky way.  After that some computations lead to \cref{X1}. 
To obtain the final form \eqref{X32}, we rewrite the estimates in \cref{X1} using 
the connection between growth and zero distribution of entire functions
and carry out a limiting process.

One important issue in relation with \eqref{X32} which we need to address is applicability. 
The formula defining $K_H(t;r)$ is an explicit expression involving the entries of $H$, but also their integrals and some
inverse functions.  Hence, $K_H(t;r)$ and even more so its integral is hard to evaluate in general. 
We believe that this is in the nature of things: for slow growth, meaning $n_H(r)=\Smallo(r)$, there cannot exist a 
plain and simple formula relating $n_H$ to $H$. 
This belief is supported by several earlier works, for example \cite{kac:1986,kac:1990} about Krein strings, or 
\cite{berg.szwarc:2014} about Jacobi matrices.

Our second main result is \cref{X3}. 
There we consider the limit circle case, i.e.\ the case when $H$ is integrable over 
the whole interval $[a,b)$.  We introduce an algorithm that produces a function $\kappa_H$ satisfying
\begin{equation}\label{X34}
	\log 2\cdot\kappa_H(r)-\bigl(\log r+\BigO(1)\bigr) \le \int_a^b K_H(t;r)\DD t
	\le 2e\cdot\kappa_H(r)\bigl(\log r+\BigO(1)\bigr)
\end{equation}
for $r>\bigl(\det\int_a^b H(t)\DD t\bigr)^{-\frac12}$.
Here the $\BigO(1)$-terms involve constants that depend only on $\int_a^b\tr H(t)\DD t$.
This result is methodologically related to the covering theorem \cite[Theorem~2]{romanov:2017}, 
which deals with diagonal Hamiltonians in the limit circle case. 
We also point out that the lower bound in \eqref{X34} is the more significant part of the theorem. 
Upper bounds for $n_H$ can be obtained also from other sources, for example \cite{romanov:2017,pruckner.woracek:sinqA}, 
while lower bounds are usually hard to find.

\Cref{X3} is a trade-off: the function $\kappa_H$ is found constructively and can be handled more easily than 
$K_H$, but we pay for computability by accepting less precision due to the multiplicative logarithmic term on the
right-hand side. 
In most examples where $\int_a^b K_H(t;r)\DD t$ can be computed, $\kappa_H(r)$ actually gives the 
correct value in the sense that there are constants $C_\pm'>0$ 
with $C_-'\kappa_H(r)\le\int_a^b K_H(t;r)\DD t\le C_+'\kappa_H(r)$. 
However, there do exist examples where the upper bound is attained, meaning that
$C_-'\kappa_H(r)\log r\le\int_a^b K_H(t;r)\DD t\le C_+'\kappa_H(r)\log r$. 
Thus both estimates given in \cref{X3} are sharp. 

An almost immediate consequence of \cref{X3} is \cref{X19}, where we establish a variant of
the following intuitive statement: if one cuts out pieces of the interval $[a,b)$, the spectrum will not become denser. 
In \cite{pruckner.woracek:sinqA} a theorem is shown which establishes this statement under a certain (in some cases quite
restrictive) assumption on the piece of $H$ that is deleted.  Using \cref{X3} it is straightforward to show that the
statement holds up to a possible logarithmic error without any additional assumptions.

As another application of the methods that are used in the proof of \cref{X3}
we consider a particular situation in more detail, namely Hamiltonians $H$ for which $\det H(t)=0$ a.e.\ 
and $H$ oscillates at the left endpoint $a$.
Using an infinite partition of the interval $(a,b)$ we obtain tight estimates for $\int_a^b K_H(t;r)\DD t$
from above and below, which, in many cases, yield the precise growth up to multiplicative constants;
see \cref{X214}.
As a consequence we obtain in \cref{X213} an inverse theorem:
for any given regularly varying comparison function $\ms f$ with index in $(\frac12,1)$ (and some $\ms f$ with index $\frac12$),
we construct an explicit Hamiltonian for which $\int_a^b K_H(t;r)\DD t$ grows like $\ms f(r)$ (up to multiplicative constants).

Let us also mention that the method leading to \cref{X3} has striking consequences in the theory of Jacobi matrices. This,
however, is beyond the scope of this paper and is presented in the forthcoming work 
\cite{reiffenstein:kachamA,reiffenstein:kachamB}.

Finally, a word about what we cannot do: we cannot characterise actual eigenvalue asymptotics. 
This is because of an intrinsic limitation of our method, which originates in the Weyl coefficient estimates 
\cite{reiffenstein:imq,langer.pruckner.woracek:heniest} and leads to constants $C_\pm$ in \eqref{X32}
that cannot be made arbitrarily close to each other.
Therefore we cannot detect small oscillations of the Stieltjes transform.
On the other hand, our results are still applicable in situations where the eigenvalues
do not have a simple asymptotic behaviour.

\subsubsection*{Structure of the manuscript}

After this introduction the paper is divided into six sections, of which we give a brief outline. 

In \cref{X131} we collect various preliminaries.  This includes basic facts from the theory of canonical systems, 
which we present to make the paper as self-contained as possible.  Furthermore, we study the function $\det\Omega_H(s,t)$
where $\Omega_H(s,t)\DE\int_s^t H(x)\DD x$.  This function plays a crucial role in the whole paper; we think of it as a
measure for the `speed of rotation' of $H$ when $H$ has rank one and $\tr H(t)=1$ a.e.;
see, in particular, \cref{X54}, which backs this intuition and expresses $\det\Omega_H(s,t)$ in terms
of the `angle'  $\varphi(t)$ that is such that $\binom{\cos\varphi(t)}{\sin\varphi(t)}\perp\ker H(t)$.

\Cref{X132} is devoted exclusively to the statement and proof of our first main result, \cref{X79}.
In \cref{X133} we use two variants of Karamata's Tauberian and Abelian theorems to deduce knowledge about $n_H$ itself, 
where the growth of $n_H$ is measured relative to comparison functions $\ms f$ that are regularly varying 
in Karamata's sense (we recall this notion and some frequently used properties in \cref{X144}). 
The major results in this section are \cref{X6,X7}.  The first one characterises `convergence class'
conditions, the second one deals with `finite type' properties (terminology borrowed from complex analysis). 
We also provide uniform estimates for the monodromy matrix on circles with radius $r$ in the limit circle case;
see \cref{X45}.
Particular attention is paid to the case of summability with respect to the comparison function $\ms f(r)\DE r$, 
which corresponds to trace class resolvents.  Notably, in \cref{X53} we give a plain and simple formula for the trace 
of the inverse of the model operator provided it is of trace class. 

In \cref{X134} we turn to our second main result, \cref{X3}, which gives a pointwise estimate for
the Stieltjes transform of $\frac{n_H(\sqrt{t})}{t}$ in terms of the more computable function $\kappa_H(r)$.  
After having proved this theorem, we formulate and prove \cref{X19} about cutting out pieces of the domain $(a,b$).
Further, we discuss a few statements that emphasise the property that 
our theorems hold with universal constants.  
In particular, we prove pointwise estimates of the fundamental solution, instead of estimates for the 
limit superior as $|z|\to\infty$, which are common in the literature. 
This includes estimates when $\det H(t)=0$, $\tr H(t)=1$ a.e.\ and the angle $\varphi$ satisfies a H\"older condition
or is of bounded variation; see \cref{X36}.

The situation where $H$ satisfies $\det H(t)=0$, $\tr H(t)=1$ a.e.\ and $H$ oscillates at the left endpoint
is treated in \cref{X210}.  We start with the main estimates in \cref{X214} and continue with 
a couple of examples in \cref{X151}.  In particular, in \cref{X42} we study $H$ with $\varphi$ being
a chirp signal, where we determine the exact growth of $\int_a^b K_H(t;r)\DD t$ up to multiplicative constants.
In \cref{X212} we state and prove the above mentioned inverse theorem.

The final \cref{X135} contains further applications and illustrations of our main results. 
In \cref{X8} we present a comparison result where we compare two different Hamiltonians.
\Cref{X152,X153} contain some examples:
in \cref{X127} we consider a Hamiltonian with $\varphi$ being a Weierstra{\ss} function,
which proves sharpness of the bound in \cref{X36}\,(ii) for $\varphi$ satisfying a H\"older condition;
in \cref{X5} we study a non-trivial Hamiltonian that shows sharpness of the upper bound in \cref{X3}.

\section{Preliminaries}
\label{X131}
\subsection{Nevanlinna functions}
\label{X136}

In the spectral theory of self-adjoint operators a particular class of analytic functions plays a prominent role. 

\begin{Definition}\label{X60}
	We call $q$ a \emph{Nevanlinna function} if $q$ is an analytic function on the open 
	upper half-plane $\bb C^+=\{z\in\bb C\mid\Im z>0\}$ and maps this half-plane into its closure $\bb C^+\cup\bb R$.
	The set of all Nevanlinna functions is denoted by $\mc N$. 
\end{Definition}

\noindent
Sometimes, also the function constant equal to $\infty$ is included in $\mc N$. 
We will not do this, and write $\mc N\cup\{\infty\}$ whenever necessary. 

In the literature such functions are also often called \emph{Herglotz functions}. 
Indeed, G.~Herglotz proved that they admit an integral representation, which is recalled
in the following theorem.

\begin{theorem}\label{X71}
	A function $q\DF\bb C^+\to\bb C$ is a Nevanlinna function if and only if there exist $\alpha\in\bb R$, $\beta\ge 0$, 
	and a positive Borel measure $\mu$ on $\bb R$ with $\int_{\bb R}\frac{\D\mu(t)}{1+t^2}<\infty$ such that 
	\begin{equation}\label{X18}
		q(z) = \alpha+\beta z+\int_{\bb R}\Bigl(\frac 1{t-z}-\frac t{1+t^2}\Bigr)\DD\mu(t), \qquad z\in\bb C^+.
	\end{equation}
\end{theorem}

\noindent
The data in this integral representation can be obtained from $q$ by explicit formulae, namely, $\alpha=\Re q(i)$, 
$\beta=\lim_{y\to\infty}\frac 1{iy}q(iy)$, and $\mu$ by the Stieltjes inversion formula; see, e.g.\ 
\cite[\S5.4]{rosenblum.rovnyak:1994}. 

\subsection{Canonical systems}
\label{X137}

As already said in the introduction, we study systems \eqref{X92} whose Hamiltonian $H$ 
has some integrability and positivity properties. 

\begin{Definition}\label{X22}
	Let $-\infty<a<b\le\infty$.  We denote by $\bb H_{a,b}$ the set of all measurable functions 
	$H\DF(a,b)\to\bb R^{2\times 2}$ such that 
	\begin{Itemize}
	\item $H\in L^1\bigl((a,c),\bb R^{2\times 2}\bigr)$ for all $c\in(a,b)$,
	\item $H(t)\ge 0$ for $t\in(a,b)$ a.e.,
	\item $\{t\in(a,b)\DS H(t)=0\}$ has measure zero.
	\end{Itemize}
	Functions that coincide a.e.\ will tacitly be identified.
	We say that $H$ is in the \emph{limit circle case} if $H\in L^1\bigl((a,b),\bb R^{2\times 2}\bigr)$, and that $H$ is in 
	the \emph{limit point case} otherwise.
	
	With a slight abuse of notation we set $\dom H\DE[a,b]$ if $H$ is in the limit circle case,
	and $\dom H\DE[a,b)$ if $H$ is in the limit point case.
\end{Definition}

\noindent
If there is need to refer to the entries of $H$, we always use the generic notation 
\begin{equation}\label{X83}
	H(t) = \begin{pmatrix} h_1(t) & h_3(t) \\[0.5ex] h_3(t) & h_2(t) \end{pmatrix}.
\end{equation}
Due to positivity, $H$ is in the limit circle case if and only if 
\[
	\int_a^b\tr H(t)\DD t<\infty.
\]
Intervals where $H$ is a scalar multiple of a constant rank $1$ matrix play an intrinsically exceptional role. 
Let us set
\begin{equation}\label{X104}
	\xi_\phi \DE \binom{\cos\phi}{\sin\phi}, \qquad\phi\in\bb R.
\end{equation}

\begin{Definition}\label{X43}
	Let $H\in\bb H_{a,b}$.  A non-empty interval $(c,d)\subseteq(a,b)$ is called \emph{indivisible} for $H$
	if there exists $\phi\in\bb R$ such that 
	\[
		H(t) = \tr H(t)\cdot\xi_\phi\xi_\phi^T \qquad\text{for } t\in(c,d)\text{ a.e.}
	\]
	The number $\phi$, which is determined up to additive integer multiples of $\pi$, 
	is called the \emph{type} of the indivisible interval $(c,d)$. 

	We denote by $\domr H$ the set of all $t\in(a,b)$ that are not inner point of an indivisible interval. 
	
	We call $H$ \emph{definite} if $(a,b)$ is not indivisible.
\end{Definition}

\noindent
Next, let us very briefly recall the operator model for $H\in\bb H_{a,b}$.  
A Hilbert space $L^2(H)$ is defined as a certain closed subspace of the usual $L^2$-space of the 
matrix measure $H(t)\DD t$, namely, as the subspace containing all those (equivalence classes $H(t)\DD t$--a.e.\ of) 
functions $f$ that have the following property: for every indivisible interval
$(c,d)\subseteq(a,b)$ the function $\xi_\phi^Tf(t)$ is constant a.e.\ on $(c,d)$ where $\phi$ is the type of $(c,d)$. 
We denote the equivalence class of a function $f$ by $[f]_H$ if we have to explicitly refer to it.

In the definition below we introduce the model operator $A_H$, distinguishing between limit circle and limit point cases. 

\begin{Definition}
Let $H \in\bb H_{a,b}$ be definite. We define (the graph of) the model operator (or linear relation) $A_H$ as follows.
\begin{Ilist}
\item 
	When $H$ is in the limit circle case, set 
	\[
		A_H \DE \biggl\{
		(f,g)\in L^2(H)\times L^2(H)\;\Big|\; 
		\raisebox{1.5ex}{\parbox[t]{44ex}{\small $f$ has an absolutely continuous representative \\[0.5ex] 
		with $f'=JHg$ a.e., $(1,0)f(a)=(0,1)f(b)=0$}}
		\biggr\}.
	\]
\item 
	When $H$ is in the limit point case, set 
	\[
		A_H\DE\biggl\{
		(f,g)\in L^2(H)\times L^2(H)\;\Big|\; 
		\raisebox{1.5ex}{\parbox[t]{44ex}{\small $f$ has an absolutely continuous representative \\[0.5ex]
		with $f'=JHg$ a.e., $(1,0)f(a)=0$}}
		\biggr\}.
	\]
\end{Ilist}
In both cases $A_H$ is a self-adjoint linear relation.  We write $\sigma_H$ for the spectrum of (the operator part of) $A_H$. 
\end{Definition}

\noindent
Definiteness of $H$ is needed to ensure that absolutely continuous representatives are unique and 
hence boundary values are well defined.  The spectrum of $A_H$ is always simple. 
If $H$ is in the limit circle case, then $A_H$ is invertible and $A_H^{-1}$ is compact. 

\medskip

\begin{remark}\label{X113}
	Let us mention here that $A_H$ is multi-valued if and only if there exists $c\in(a,b)$ 
	such that $h_2(t)=0$ for $t\in(a,c)$ a.e., that is, $(a,c)$ is an indivisible interval 
	of type $0$.
\end{remark}

\begin{Remark}\label{X69}
	Let us consider the situation when $H$ ends with an indivisible interval, $(c,b)$,
	of type $\frac{\pi}{2}$ with $c\in(a,b)$, and let the interval $(c,b)$ 
	be maximal with this property.  We distinguish two cases.
	\begin{Enumerate}
	\item
		Assume that $H$ is in the limit point case.
		One can then identify $L^2(H)$ with $L^2(H|_{(a,c)})$
		since each element $f\in L^2(H)$ contains a representative that vanishes on $(c,b)$.
		Also the relations $A_H$ and $A_{H|_{(a,c)}}$ are equivalent under this identification
		as each absolutely continuous representative in $L^2(H)$ satisfies $(0,1)f(c)=0$.
	\item
		If $H$ is in the limit circle case, then $L^2(H)$ is a one-dimensional extension
		of $L^2(H|_{(a,c)})$, but the operator parts of $A_H$ and $A_{H|_{(a,c)}}$
		are unitarily equivalent.
	\end{Enumerate}
\end{Remark}

\noindent
Discreteness of the spectrum, and membership in Schatten--von Neumann ideals (for $p>1$) can be characterised. 
We recall the following two facts from \cite[Theorems~1.1 and 1.3]{romanov.woracek:ideal} 
(the first equivalence is also shown in \cite[Theorem~1.3]{remling.scarbrough:2020a} with a different method): 
\begin{Itemize}
\item 
	$A_H$ has compact resolvents, i.e.\ $\sigma_H$ is discrete if and only if 
	there exists $\phi\in\bb R$ such that 
	\begin{align}
		\label{X123}
		& \int_a^b\xi_\phi^TH(t)\xi_\phi\DD t<\infty,
		\\[1ex]
		\nonumber
		& \lim_{t\to b} \biggl(\int\limits_a^t\xi_{\phi+\frac\pi 2}^TH(s)\xi_{\phi+\frac\pi 2}\DD s
		\cdot\int\limits_t^b\xi_\phi^TH(s)\xi_\phi\DD s\biggr) 
		= 0.
	\end{align}
\item
	For $p>1$, $A_H$ has resolvents belonging to the Schatten--von Neumann class $\mf S_p$, i.e.\
	$\sum_{\lambda\in\sigma_H\setminus\{0\}}\frac{1}{\lambda^p}<\infty$, if and only if
	there exists $\phi\in \bb R$ such that
	\[
		\hspace*{-3ex}\eqref{X123} \;\;\text{holds and}\quad \int\limits_a^b 
		\bigg(\int\limits_t^b\xi_\phi^TH(s)\xi_\phi\DD s\bigg)^{\frac{p}{2}-1}
		\bigg(\int\limits_a^t\xi_{\phi+\frac\pi 2}^TH(s)\xi_{\phi+\frac\pi 2}\DD s\bigg)^{\frac{p}{2}}
		\cdot\xi_\phi^TH(s)\xi_\phi\DD t 
		< \infty.
	\]
	Note that the expression simplifies when $p=2$; this particular case (i.e.\ the characterisation of Hilbert--Schmidt class)
	was already proved in \cite{kaltenbaeck.woracek:hskansys}.
\end{Itemize}

\begin{Remark}\label{X16}
	Applying rotation isomorphisms always allows one to reduce to the case $\phi=0$.
	To be more precise, consider the Hamiltonian
	\begin{equation}\label{X203}
		H_\alpha \DE N_\alpha HN_\alpha^{-1} \qquad\text{with}\quad
		N_\alpha \DE \begin{pmatrix} \cos\alpha & \sin\alpha \\[0.5ex] -\sin\alpha & \cos\alpha \end{pmatrix}.
	\end{equation}
	The operator (or relation) $A_{H_\alpha}$ is unitarily equivalent to $A_H^{(\alpha)}$,
	where the latter is defined as $A_H$ but with different boundary conditions, 
	namely $\xi_\alpha^Tf(a)=0$ at $a$, and $\xi_{\alpha+\frac\pi2}f(b)=0$ in the limit circle case.
	Since $(A_H-z)^{-1}-(A_H^{(\alpha)}-z)^{-1}$ is a finite-rank operator for $z\in\bb C\setminus\bb R$,
	the asymptotic properties of $\sigma_H$ and $\sigma_{H_\alpha}$ are the same.
	In the following we therefore often assume that \eqref{X123} holds for $\phi=0$, i.e.\
	\begin{equation}\label{X154}
		\int_a^b h_1(t)\DD t < \infty.
	\end{equation}
	Note that, under this assumption, $0\notin\sigma_H$; this follows since the only constant
	satisfying the boundary condition at $a$ is $\binom01$; however, this does not satisfy
	the boundary condition at $b$ when $H$ is in the limit circle case, and it is not
	in the space $L^2(H)$ when $H$ is in the limit point case.
	If we assume that \eqref{X154} holds, then the characterisations for discrete spectrum
	and Schatten--von Neumann resolvents simplify:
	\begin{align*}
		\sigma_H\text{ discrete}
		\quad&\Longleftrightarrow\quad
		\lim_{t\to b} \bigg(\int\limits_t^b h_1(s)\DD s\cdot\int\limits_a^t h_2(s)\DD s \bigg) = 0,
		\\
		\sum_{\lambda\in\sigma_H}\frac{1}{\lambda^p} < \infty	
		\quad&\Longleftrightarrow\quad
		\int\limits_a^b 
		\bigg(\int\limits_t^b h_1(s)\DD s\bigg)^{\frac{p}{2}-1}
		\bigg(\int\limits_a^t h_2(s)\DD s\bigg)^{\frac{p}{2}}\cdot h_1(s)\DD t < \infty
	\end{align*}
	for $p>1$.
\end{Remark}

\subsection{Construction of the Weyl coefficient}
\label{X138}

For the proofs of our main theorems we need also the concept of the Weyl coefficient.
Given $H\in\bb H_{a,b}$, the initial value problem 
\begin{equation}\label{X51}
	\left\{
	\begin{array}{l}
		\dfrac{\partial}{\partial t}W_H(t;z)J=zW_H(t;z)H(t), \qquad t\in(a,b),
		\\[3ex]
		W_H(a;z)=I.
	\end{array}
	\right.
\end{equation}
has a unique solution $W_H(t;z)$ on $[a,b)$ where $W_H(t;z)\in\bb C^{2\times2}$. 
We refer to $W_H(t;z)$ as the \emph{fundamental solution} of $H$.
Note that the transposes of the rows of $W_H(\,\cdot\,;z)$ are linearly independent solutions of \eqref{X92}.
If $H$ is in the limit circle case, the solution $W_H$ can be extended continuously to $b$, 
and we call $W_H(b;\cdot)$ the \emph{monodromy matrix} of $H$.
Note that in both cases (limit point and limit circle) the fundamental solution is defined on $\dom H$.
For each fixed $t\in\dom H$ the function $W_H(t;z)$ is an entire function in $z$. 
Its entries $w_{H,ij}(t;z)$, $i,j=1,2$, are real along the real axis, have only real and simple zeros, 
and are of bounded type in $\bb C^+$; the latter means that their restrictions to $\bb C^+$ are
quotients of bounded analytic functions.
Note also that the orders (as entire functions) of all entries coincide;
we denote by $\rho_H$ the common order of the monodromy matrix when $H$ is in the limit circle case.
\begin{Ilist}
\item
	Assume that $H$ is in the limit circle case.  Then the function 
	\[
		q_H(z) \DE \frac{w_{H,12}(b;z)}{w_{H,22}(b;z)}
	\]
	is a Nevanlinna function.  The measure in its Herglotz integral representation \eqref{X18}
	is a scalar spectral measure for the operator part of $A_H$.
	In particular, $\sigma_H$ is discrete and coincides with the set of zeros of the function $w_{H,22}(b;z)$.
\item 
	Assume that $H$ is in the limit point case.  Then the limit
	\[
		q_H(z) \DE \lim_{t\to b}\frac{w_{H,12}(t;z)}{w_{H,22}(t;z)}
	\]
	exists locally uniformly as a function from $\bb C\setminus\bb R$ to $\bb C\cup\{\infty\}$ 
	and belongs to $\mc N\cup\{\infty\}$.
	The measure in its Herglotz integral representation (in the case when $q_H\ne\infty$) 
	is a scalar spectral measure for the operator part of $A_H$. 
	We shall refer to the function $q_H$ as the \emph{Weyl coefficient} of $H$. 
\end{Ilist}

\noindent
Note that $A_H$ is multi-valued if and only if the linear term $bz$ in the integral representation \eqref{X18} 
of $q_H$ is present, which, in turn, is equivalent to the fact that $H$ starts with 
an indivisible interval of type $0$ at the left endpoint $a$.

Let us also recall reparameterisations of Hamiltonians.
Let $H_1$ and $H_2$ be Hamiltonians on the intervals $(a_1,b_1)$ and $(a_2,b_2)$, respectively.
We say that $H_1$ and $H_2$ are \emph{reparameterisations} of each other
if there exists a strictly increasing bijection  $\gamma:(a_1,b_1)\to (a_2,b_2)$
such that both $\gamma$ and $\gamma^{-1}$ are locally absolutely continuous and
\begin{equation}\label{X170}
	H_1(t) = H_2(\gamma(t))\gamma'(t),  \qquad t\in (a_1,b_1)\ a.e.
\end{equation}
If \eqref{X170} holds and $y$ is a solution of \eqref{X92} with $H$ replaced by $H_2$,
then $y\circ\gamma$ is a solution of \eqref{X92} with $H$ replaced by $H_1$.
For each Hamiltonian $H$ on an interval $(a,b)$ there is a unique
\emph{trace-normalised} reparameterisation $H_2$ of $H$ defined on $(0,b_2)$,
i.e.\ $\tr H_2(t)=1$ for a.e.\ $t\in(0,b_2)$.
One obtains this reparameterisation by choosing $\gamma$ in \eqref{X170} to be
the inverse function of $x\mapsto\int_a^x \tr H(t)\DD t$.
It is easy to see that the Weyl coefficients of Hamiltonians that are reparameterisations
of each other coincide.

\subsection[The function $\det\Omega_H$]{The function \boldmath{$\det\Omega_H$}}
\label{X139}

Given a Hamiltonian $H\in\bb H_{a,b}$, we set
\[
	\Omega_H(s,t) \DE \int_s^t H(u)\DD u,
\]
where $s,t\in\dom H$ with $s\le t$.
We shall use the notation (recall \eqref{X83})
\begin{equation}
\label{X48}
	\omega_{H,j}(s,t)\DE\int_s^t h_j(u)\DD u, \quad j=1,2,3,
\end{equation}
so that 
\[
	\Omega_H(s,t)=
	\begin{pmatrix}
		\omega_{H,1}(s,t) & \omega_{H,3}(s,t)
		\\[0.5ex]
		\omega_{H,3}(s,t) & \omega_{H,2}(s,t)
	\end{pmatrix}.
\]
Since $H(t)\ge 0$ and $\tr H(t)>0$ a.e., we have $\Omega_H(s,t)\ge 0$ for all $s\le t$ and $\Omega_H(s,t)=0$ 
if and only if $s=t$. 

We make extensive use of the function 
\[
	\det\Omega_H\FD{\{(s,t)\in\dom H\times\dom H\DS s\le t\}}{[0,\infty)}{6}{(s,t)}{\det\Omega_H(s,t).}
\]
If there is no risk of ambiguity, we drop explicit notation of $H$ and write 
$\Omega(s,t)$ and $\omega_j(s,t)$ instead of $\Omega_H(s,t)$ and $\omega_{H,j}(s,t)$, respectively.

\begin{lemma}\label{X9}
	Let $H\in\bb H_{a,b}$. 
	\begin{Enumerate}
	\item 
		The function $\det\Omega$ is continuous. 
	\item 
		We have $\det\Omega(s,t)=0$ if and only if $(s,t)$ is indivisible.
	\item 
		Let $s,t\in\dom H$ with $s<t$ and assume that $\det\Omega(s,t)>0$.  
		Then, for all $s'<t'$ with $[s,t] \subseteq [s',t'] \subseteq \dom H$ we have 
		\[
			\frac{\det\Omega(s,t)}{\omega_{i}(s,t)} \le \frac{\det\Omega(s',t')}{\omega_{i}(s',t')},
			\qquad i\in\{1,2\}.
		\]
	\item 
		Fix $s\in[a,b)$ and set
		\[
			c(s) \DE
			\begin{cases}
				\sup\{t\in(s,b]\DS(s,t)\text{ indivisible}\} \CAS \text{$(s,x)$ is indivisible for some $x\in(s,b]$},
				\\
				s & \text{otherwise}.
			\end{cases}
		\]
		Then $\det\Omega(s,t)=0$ for $t\in[s,c(s))$, and the function $t\mapsto\det\Omega(s,t)$ 
		is strictly increasing on $[c(s),b)$. 
	\item 
		Fix $t\in\dom H\setminus\{a\}$ and set 
		\[
			d(t) \DE
			\begin{cases}
				\inf\{s\in[a,t)\DS(s,t)\text{ indivisible}\} \CAS \text{$(x,t)$ is indivisible for some $x\in[a,t)$},
				\\
				t & \text{otherwise}.
			\end{cases}
		\]
		Then $\det\Omega(s,t)=0$ for $s\in[d(t),t]$, and the function $s\mapsto\det\Omega(s,t)$ is strictly decreasing on $[a,d(t)]$. 
	\end{Enumerate}
\end{lemma}

\begin{proof}
	Item \textup{(i)} is clear since $\Omega(s,t)$ is an integral over an $L^1_{\rm loc}$-function. 
	For the proof of \textup{(ii)} note that, for each $x\in\bb C^2$, 
	\[
		\big(\Omega(s,t)x,x\big) = \int_s^t\big(H(u)x,x\big)\DD u.
	\]
	Using that $H(u)$ and $\Omega(s,t)$ are positive semi-definite we see that $\Omega(s,t)x=0$ if and only if 
	$H(u)x=0$ for a.a.\ $u\in(s,t)$.  Existence of a non-zero vector $x$ with this property is equivalent
	to $(a,b)$ being indivisible.

	Item (iii) follows by differentiation; see \cite[(2.17)]{langer.pruckner.woracek:heniest}.
	
	The first statement in (iv) and the fact that $\det\Omega(s,t)>0$ for $t\in(c(s),b)$ follow from (ii).
	For the second statement in (iv) let $t_1,t_2\in(c(s),b)$ with $t_1<t_2$.
	Since $t\mapsto\tr\Omega(s,t)$ is strictly increasing, we have $\omega_i(s,t_1)<\omega_i(s,t_2)$
	for at least one $i\in\{1,2\}$.  Now the claim follows from (iii).
	
	Item (v) is analogous to (iv).
\end{proof}

\begin{lemma}\label{X10}
	Let $H\in\bb H_{a,b}$ be definite.  Then $H$ is in the limit point case 
	if and only if $\lim_{t\to b}\det\Omega(a,t)=\infty$.
\end{lemma}

\begin{proof}
	For a matrix $M \in \bb C^{2 \times 2}$ with $M \ge 0$, we denote by $\lambda_1(M)$ and $\lambda_2(M)$ 
	its eigenvalues enumerated such that $\lambda_1(M) \ge \lambda_2(M)$.  Since 
	\[
		\int_a^t \tr H(u)\DD u = \lambda_1\big(\Omega(a,t)\big)+\lambda_2\big(\Omega(a,t)\big),
	\]
	limit point case takes place if and only if $\lim_{t\to b}\lambda_1(\Omega(a,t))=\infty$.
	The definiteness of $H$ implies that we can find $c\in(a,b)$ such that $\lambda_2(\Omega(a,c))>0$, 
	and hence 
	\[
		\lim_{t\to b}\lambda_1(\Omega(a,t))=\infty 
		\quad\Longleftrightarrow\quad 
		\lim_{t\to b}\Big(\lambda_1\big(\Omega(a,t)\big)\cdot\lambda_2\big(\Omega(a,t)\big)\Big) = \infty
	\]
	because $\lambda_2(\Omega(a,t))\ge\lambda_2(\Omega(a,c))$ for $t\ge c$.
\end{proof}

\noindent
If $H$ is in the limit point case, we set
\[
	\det\Omega(s,b) \DE \lim_{t\to b}\det\Omega(s,t), \qquad s\in[a,b);
\]
note that the improper limit exists because of the monotonicity of $\det\Omega(s,\cdot)$. 
The function $\det\Omega$, now defined on $\{(s,t)\in[a,b)\times[a,b]\DS s\le t\}$, 
again has all properties listed in \cref{X9}.

The following notion plays a central role in our present investigations.

\begin{Definition}
\label{X11}
	Let $H\in\bb H_{a,b}$, $r_0\ge 0$, and $c_-,c_+>0$. 
	\begin{Enumerate}
	\item 
		We call $\hat t$ a \emph{compatible function} for $H,r_0$ with constants $c_-,c_+$ if 
		\begin{align*}
			& \hat t\DF(r_0,\infty)\to(a,b),
			\\
			& \forall r\in(r_0,\infty)\DP
			\frac{c_-}{r^2} \le \det\Omega\bigl(a,\hat t(r)\bigr) \le \frac{c_+}{r^2}.
		\end{align*}
	\item 
		If $\hat t$ is a compatible function for $H,r_0$, we set
		\[
			\Gamma(\hat t) \DE \bigl\{(t,r)\in(a,b)\times(r_0,\infty)\DS t \ge \hat t(r)\bigr\}.
		\]
	\item 
		We call $(\hat t,\hat s)$ a \emph{compatible pair} for $H,r_0$ with constants $c_-,c_+$ if 
		$\hat t$ is a compatible function for $H,r_0$ with constants $c_-,c_+$ and 
		\begin{align*}
			& \hat s\DF\Gamma(\hat t)\to[a,b),
			\\
			& \forall (t,r)\in\Gamma(\hat t)\DP \hat s(t;r)\le t, \quad 
			\frac{c_-}{r^2} \le \det\Omega\bigl(\hat s(t;r),t\bigr) \le \frac{c_+}{r^2}.
		\end{align*}
	\end{Enumerate}
\end{Definition}

\noindent
The properties of $\det\Omega$ established in \cref{X9} imply that, for definite Hamiltonians, 
compatible pairs exist as the following proposition shows.

\begin{proposition}\label{X39}
	Let $H\in\bb H_{a,b}$ be definite, let $c>0$ and set $r_0=\bigl(\frac{c}{\det\Omega(a,b)}\bigr)^{\frac 12}$ 
	if $H$ is in the limit circle case and $r_0=0$ if $H$ is in the limit point case.  
	Then there exists a unique compatible pair $(\hat t,\hat s)$ for $H,r_0$ with constants $c,c$.
\end{proposition}

\begin{proof}
	Let $r>r_0$.  Then $0<\frac{c}{r^2}<\det\Omega(a,b)$ since $\Omega(a,b)=\infty$ in the limit point case
	by \cref{X10}.  Further, let $c(a)$ be as in \cref{X9}\,(iv). 
	Since the function $t\mapsto\det\Omega(a,t)$ is a continuous and increasing bijection 
	from $[c(a),b)$ onto $[0,\det\Omega(a,b))$, we find a unique point $\hat t(r)\in(a,b)$ such that 
	\[
		\det\Omega\bigl(a,\hat t(r)\bigr) = \frac{c}{r^2}.
	\]
	Now let $(t,r)\in\Gamma(\hat t)$.  Then 
	\[
		0 < \frac{c}{r^2} = \det\Omega\bigl(a,\hat t(r)\bigr) \le \det\Omega(a,t),
	\]
	and since $s\mapsto\det\Omega(s,t)$ is a continuous and decreasing bijection from $[a,d(t)]$ onto
	$[0,\det\Omega(a,t)]$, we find a unique point $\hat s(t;r)\in[a,t)$ 
	such that 
	\[
		\det\Omega\bigl(\hat s(t;r),t\bigr) = \frac c{r^2}.
	\]
	Uniqueness of the compatible pair in this situation is clear.
\end{proof}

\begin{remark}\label{X46}
\rule{0ex}{1ex}
\begin{Enumerate}
\item
	For many purposes it would be enough to work with the compatible pairs constructed in \Cref{X39}, 
	even using only $c=1$.  Our reasons to introduce the terminology of a `compatible pair' at all, 
	and to use the additional parameters $r_0,c_-,c_+$, are that this makes it easier to trace constants 
	in estimates and to actually apply the abstract results in concrete situations. 
\item
	Let $r_0,c_-,c_+$ be as in \cref{X11} and let $(\hat t,\hat s)$ be a compatible pair for $H,r_0$
	with constants $c_-,c_+$.  Further, let $r_0'\ge r_0$, $c_-'\le c_-$, $c_+'\ge c_+$.
	Then $(\hat t',\hat s')$ with $\hat t'=\hat t|_{(r_0',\infty)}$ and $\hat s'=\hat s|_{\Gamma(\hat t)\cap((a,b)\times(r_0',\infty))}$ 
	is a compatible pair for $H,r_0'$ with constants $c_-',c_+'$.
\end{Enumerate}
\end{remark}

\subsection[Behaviour of $q_H$ towards $i\infty$]{%
Behaviour of \boldmath{$q_H$} towards \boldmath{$i\infty$}}
\label{X140}

We use the following notation to compare functions up to multiplicative constants.

\begin{Notation}\label{X47}
\rule{1ex}{0ex}
\begin{Enumerate}
\item
	Let $X$ be a set and let $f,g\DF X\to[0,\infty)$ be functions. 
	Then we write $f\lesssim g$ (or $f(x)\lesssim g(x)$) to say that 
	\[
		\exists c>0\DQ\forall x\in X\DP f(x) \le cg(x).
	\]
	Moreover, we write
	\[
		f\gtrsim g \;\;\DI\;\; g\lesssim f, \qquad f\asymp g \;\;\DI\;\; f\lesssim g \;\wedge\; f\gtrsim g.
	\]
	When $f,g:[x_0,\infty)\to(0,\infty)$ with some $x_0>0$, we write 
	\begin{alignat*}{3}
		f(x)&\sim g(x) \quad&&\text{as}\;\; x\to\infty \qquad&&\text{if}\;\; \lim_{x\to\infty}\tfrac{f(x)}{g(x)}=1,
		\\[1ex] 
		f(x)&\ll g(x) \quad&&\text{as}\;\; x\to\infty \qquad&&\text{if}\;\; \lim_{x\to\infty}\tfrac{f(x)}{g(x)}=0.
	\end{alignat*}
\item
	For $A,B\in[0,\infty]$ we write $A\diamond B$ to say that $A$ and $B$ are either both infinite,
	both finite and positive, or both zero.
\end{Enumerate}
\end{Notation}

\noindent
The following result is shown in \cite[Theorem~1.1]{reiffenstein:imq}. It is a crucial tool for the present paper. 

Recall the notation $\omega_{H,j}$ (abbreviated to $\omega_j$) from \eqref{X48}.

\begin{theorem}\label{X49}
	Let $H\in\bb H_{a,b}$ be definite and in the limit point case. 
	Further, let $r_0\ge 0$, $c_-,c_+>0$, and let $\hat t$ be a compatible function for $H,r_0$ 
	with constants $c_-,c_+$.  Then 
	\begin{alignat*}{2}
		|q_H(ir)| &\asymp \sqrt{\frac{\omega_{1}(a,\hat t(r))}{\omega_{2}(a,\hat t(r))}\,},
		\qquad && r\in(r_0,\infty),
		\\[1ex]
		\Im q_H(ir) &\asymp \frac{1}{r\omega_{2}(a,\hat t(r))},
		\qquad && r\in(r_0,\infty).
	\end{alignat*}
	The constants implicit in the two relations `\/$\asymp$' depend on $c_-,c_+$ 
	but not on $H,r_0,\hat t$. 
\end{theorem}

\section{Growth of the eigenvalue counting function}
\label{X132}

In the theorem below we give a formula that determines the growth of the Stieltjes transform of the counting function 
of the spectrum of $H$ provided that $\sigma_H$ is discrete 
and $\sum_{\lambda\in\sigma_H\setminus\{0\}}\frac{1}{\lambda^2}<\infty$.

Recall the notion of a compatible pair from \cref{X11}, and that, for $c_-\le c_+$, such pairs always exist by \cref{X39}.
Moreover, recall \cref{X16}, and denote by $\mathds{1}_M$ the characteristic function of a set $M$.
We start with the definition of an integral kernel, which is a central piece in our estimates.

\begin{Definition}\label{X72}
	Let $H\in\bb H_{a,b}$ be a definite Hamiltonian, let $r_0\ge0$, $c_-,c_+>0$, 
	assume that $(\hat t,\hat s)$ is a compatible pair
	for $H,r_0$ with constants $c_-,c_+$
	\begin{align}
		K_H(t;r) \DE \mathds{1}_{[a,\hat t(r))}(t)
		\frac{\omega_{H,2}(a,t)h_1(t)}{\frac{c_+}{r^{2}}+(\omega_{H,3}(a,t))^2}
		+ \mathds{1}_{[\hat t(r),b)}(t) \frac{h_1(t)}{\omega_{H,1}\bigl(\hat s(t;r),t\bigr)},
		\hspace*{10ex} &
		\label{X15}
		\\[1ex]
		t\in\dom H,\,r\in(r_0,\infty). &
		\nonumber
	\end{align}
\end{Definition}

\noindent
The following theorem is one of our main results.

\begin{theorem}\label{X79}
	Let $H\in\bb H_{a,b}$ be a definite Hamiltonian with discrete spectrum 
	that satisfies $\int_a^b h_1(t)\DD t<\infty$.
	Let $r_0\ge 0$, $c_-,c_+>0$, assume that $(\hat t,\hat s)$ is a compatible pair 
	for $H,r_0$ with constants $c_-,c_+$, and let $K_H$ be as in \cref{X72}.
	Further, set 
	\begin{equation}\label{X81}
		n_H(r)\DE\#\big\{\lambda\in\sigma_H\DS |\lambda|<r\big\}.
	\end{equation}
	Then
	\begin{equation}\label{X90}
		\int_0^\infty\frac{1}{t+r^2}\cdot\frac{n_H(\sqrt t)}t\DD t
		\asymp \frac{1}{r^2}\int_a^b K_H(t;r)\DD t,
		\qquad r>r_0,
	\end{equation}
	where this relation includes that one side is finite if and only if the other side is. 
	The constants implicit in `\/$\asymp$' depend on $c_-,c_+$ but not on $H,r_0,\hat t,\hat s$.
\end{theorem}

\noindent
The formula \eqref{X90} gives a meaningful result about $n_H$ only if the integrals on either side are finite. 
The integral on the left-hand side of \eqref{X90} is finite for some (equivalently, for every) $r>0$ if and only if 
\[
	\int_0^\infty \frac{n_H(\sqrt t)}{t^2}\DD t<\infty,
\]
which, in turn, is equivalent to $\sum\limits_{\lambda\in\sigma_H}\frac{1}{\lambda^2}<\infty$. 

\begin{Remark}\label{X17}
	\phantom{}
	\begin{Enumerate}
	\item 
		If we are willing to accept that the constants in \eqref{X90} depend also on $H,r_0,\hat t,\hat s$, then the first
		summand of $K_H(t;r)$ can be dropped.  Details are given in \cref{X62}. 
	\item 
		Knowing the behaviour of the Stieltjes transform of the measure $\frac{n_H(\sqrt t)}t\DD t$ allows to employ
		classical Tauberian arguments to extract knowledge about $n_H$ itself.  Details are given in \cref{X6,X7}.
	\item 
		Applying \cref{X79} in practice is often challenging.  We return to this topic with an algorithmic approach in \cref{X3}.  
		We consider non-trivial examples where the integral over $K_H(t;r)$ can be evaluated explicitly in \cref{X42,X127,X5}.  
	\end{Enumerate}
\end{Remark}

\noindent
The proof of \cref{X79} proceeds in three stages, which are given in Sections~\ref{X141}--\ref{X143} below.
First, we use the recent result \cref{X49} about the high-energy behaviour of Weyl coefficients to determine the growth 
of the monodromy matrix $W_H(b;z)$ for a Hamiltonian in the limit circle case.  This is the essence of the argument.  
Second, we rewrite the growth of $W_H(b;z)$ in terms of the counting function $n_H(r)$.  
Finally, we make a limit argument to include the limit point case. 

\subsection{Growth of the fundamental solution}
\label{X141}

In this subsection we prove a formula for the growth of (one entry of) the fundamental solution, 
which is valid for any definite Hamiltonian $H\in\bb H_{a,b}$. 
If $H$ is in the limit circle case, the formula describes the growth of the monodromy matrix, 
and as such it is the centrepiece of the proof of \cref{X79}.

Recall the notation $W_H(t;z)=(w_{H,ij}(t;z))_{i,j=1}^2$.

\begin{theorem}
\label{X1}
	Let $H\in\bb H_{a,b}$ be definite.  Further, let $r_0\ge 0$, $c_-,c_+>0$, assume that $(\hat t,\hat s)$
	is a compatible pair for $H,r_0$ with constants $c_-,c_+$, and let $K_H(t;r)$ be as in \eqref{X15}. 
	Then, for $x\in\dom H$,
	\begin{equation}\label{X65}
		\log|w_{H,22}(x;ir)| \asymp \int_a^x K_H(t;r)\DD t,
		\qquad r>r_0,
	\end{equation}
	The constants implicit in the relation `\/$\asymp$' depend on $c_-,c_+$ but not on $H,x,r_0,\hat t,\hat s$.
\end{theorem}

\noindent 
The first step in the proof of \cref{X1} (which is contained in \cref{X2,X58}) 
is to observe that $\log|w_{H,22}(x;ir)|$ can be written as an integral over $(a,x)$,
where the integrand involves the imaginary parts of Weyl coefficients of Hamiltonians depending 
on the integration variable. 

\begin{lemma}\label{X2}
	Let $H\in\bb H_{a,b}$, $t_1,t_2\in\dom H$ with $t_1\le t_2$, and let $r>0$.  Then
	\[
		\log|w_{H,22}(t_2;ir)|-\log|w_{H,22}(t_1;ir)|
		= r\int_{t_1}^{t_2}\Im\biggl(-\frac{w_{H,21}(t;ir)}{w_{H,22}(t;ir)}\biggr)h_1(t)\DD t.
	\]
\end{lemma}

\begin{proof}
	The right lower entry of the matrix equation \eqref{X51} reads as
	\begin{equation}\label{X198}
		\frac{\partial}{\partial t}w_{H,22}(t;z) = z\Bigl(w_{H,22}(t;z)h_3(t)+w_{H,21}(t;z)h_1(t)\Bigr).
	\end{equation}
	The function $w_{H,22}(t;\cdot)$ is zero-free on $\bb C^+$, and $w_{H,22}(t;0)=1$.
	Let $\log w_{H,22}(t;\cdot)$ be the branch of the logarithm analytic in a domain
	containing $\bb C^+\cup\{0\}$ with $\log w_{H,22}(t;0)=0$. Then the function $t\mapsto\log w_{H,22}(t;z)$ is 
	differentiable a.e., and 
	\[
		\frac{\partial}{\partial t}\log w_{H,22}(t;z)
		= \frac{\frac{\partial}{\partial t}w_{H,22}(t;z)}{w_{H,22}(t;z)},
		\qquad z\in\bb C^+\cup\{0\},\, t\in(t_1,t_2).
	\]
	Hence, with \eqref{X198} we arrive at
	\begin{align*}
		\frac{\partial}{\partial t}\log|w_{H,22}(t;z)|
		&= \frac{\partial}{\partial t}\Bigl(\Re\log w_{H,22}(t;z)\Bigr)
		= \Re\frac{\frac{\partial}{\partial t}w_{H,22}(t;z)}{w_{H,22}(t;z)}
		\\[1ex]
		&= (\Re z)h_3(t) + \Re\biggl(z\frac{w_{H,21}(t;z)}{w_{H,22}(t;z)}\biggr)h_1(t).
	\end{align*}
	For $z=ir$ this yields
	\[
		\frac{\partial}{\partial t}\log|w_{H,22}(t;ir)|
		= r\Im \biggl( -\frac{w_{H,21}(t;ir)}{w_{H,22}(t;ir)}\biggr)h_1(t).
	\]
	Integrating both sides of this relation we obtain the stated formula.
\end{proof}

\begin{lemma}\label{X58}
	Let $H\in\bb H_{a,b}$ and $t\in\dom H\setminus\{a\}$.  Define the Hamiltonian
	\begin{equation}\label{X27}
		H_{(t)}(s) \DE
		\begin{cases}
			RH(a+t-s)R \CAS s\in[a,t),
			\\[1ex]
			\smmatrix 0001 \CAS s\in[t,\infty),
		\end{cases}
	\end{equation}
	with $R\DE\smmatrix0110$.  Then $H_{(t)}\in\bb H_{a,\infty}$, $H_{(t)}$ is in the limit point case, and 
	\begin{equation}\label{X63}
		q_{H_{(t)}}(z) = -\frac{w_{H,21}(t;z)}{w_{H,22}(t;z)},
		\qquad z\in\bb C\setminus\bb R.
	\end{equation}
\end{lemma}

\begin{proof}
	An elementary calculation shows that, on the interval $[a,t]$,
	the fundamental solution of $H_{(t)}$ is given by
	\[
		W_{H_{(t)}}(s;z) = RW_H(t;z)^{-1}W_H(a+t-s;z)R,
		\qquad s\in[a,t],
	\]
	and hence
	\[
		W_{H_{(t)}}(t;z) = RW_H(t;z)^{-1}R =
		\begin{pmatrix}
			w_{H,11}(t;z) & -w_{H,21}(t;z)
			\\[1ex]
			-w_{H,12}(t;z) & w_{H,22}(t;z)
		\end{pmatrix}
		.
	\]
	Since the interval $[t,\infty)$ is $H_{(t)}$-indivisible of type $\frac{\pi}{2}$,
	relation \eqref{X63} follows.
\end{proof}

\noindent
The above lemmata allow to invoke \cref{X49}.

\begin{lemma}\label{X52}
	Let $H\in\bb H_{a,b}$ be definite, $r_0\ge 0$, $c_-,c_+>0$, and let $(\hat t,\hat s)$ be a compatible pair 
	for $H,r_0$ with constants $c_-,c_+$.  Then 
	\begin{multline}
	\label{X57}
		\Im\biggl(-\frac{w_{H,21}(t;ir)}{w_{H,22}(t;ir)}\biggr)
		\asymp \mathds{1}_{[a,\hat t(r))}(t)\frac{\omega_{H,2}(a,t)}{r\bigl[\frac{c_+}{r^2}+(\omega_{H,3}(a,t))^2\bigr]}
		+ \frac{\mathds{1}_{[\hat t(r),b)}(t)}{r\omega_{H,1}\bigl(\hat s(t;r),t\bigr)},
		\\[1ex]
		r>r_0,\ t\in\dom H.
	\end{multline}
	The constants implicit in the relation `\/$\asymp$' depend on $c_-,c_+$ but not on $H,r_0,\hat t,\hat s$. 
\end{lemma}

\begin{proof}
	If $t=a$ or $(a,t)$ is indivisible of type $0$, then both sides of \eqref{X57} are equal to zero, 
	and hence \eqref{X57} holds with any constants in `\/$\asymp$' for such $t$. 
	For the rest of the proof we may thus assume that $t>a$ and $(a,t)$ is not indivisible of type $0$. 

	For $r>r_0$ and $t\in\dom H\setminus\{a\}$ we define 
	\[
		\hat t_{(t)}(r)\DE
		\begin{cases}
			t+\dfrac{1}{\omega_{H,2}(a,t)}\Big(\dfrac{c_+}{r^2}-\det\Omega_H(a,t)\Big) \CAS t\in(a,\hat t(r)),
			\\[2ex]
			a+t-\hat s(t;r) \CAS t\in[\hat t(r),b).
		\end{cases}
	\]
	Note here that $\omega_{H,2}(a,t)>0$ since $t>a$ and $(a,t)$ is not indivisible of type $0$.
	Further, observe that $\hat t_{(t)}(r)\in[a,t]$ when $\hat t(r)\le t$ and that $\hat t_{(t)}(r)\in[t,\infty)$
	when $\hat t(r)>t$;
	the latter follows from the fact that $\det\Omega_H(a,t)\le\det\Omega_H(a,\hat t(r))\le \frac{c_+}{r^2}$.

	The definition of $H_{(t)}$ in \eqref{X27} yields 
	\begin{align}\label{X111}
		\Omega_{H_{(t)}}(a,s)
		=
		\begin{cases}
			R\Omega_H(a+t-s,t)R \CAS s\in[a,t],
			\\[1ex]
			R\Omega_H(a,t)R+(s-t)\smmatrix 0001 \CAS s\in(t,\infty),
		\end{cases}
	\end{align}
	and, consequently,
	\begin{align}\label{X112}
		\det\Omega_{H_{(t)}}(a,s)
		=
		\begin{cases}
			\det\Omega_H(a+t-s,t) \CAS s\in[a,t],
			\\[1ex]
			\det\Omega_H(a,t)+(s-t)\omega_{H,2}(a,t) \CAS s\in(t,\infty).
		\end{cases}
	\end{align}
	Plugging in the definition of $\hat t_{(t)}(r)$ we obtain
	\begin{equation}\label{X64}
		\det\Omega_{H_{(t)}}\big(a,\hat t_{(t)}(r)\big)
		=
		\begin{cases}
			\det\Omega_H\bigl(\hat s(t;r),t\bigr) \CAS \hat t(r)\in[a,t],
			\\[1ex]
			\dfrac{c_+}{r^2} \CAS \hat t(r)\in(t,\infty).
		\end{cases}
	\end{equation}
	It follows from \eqref{X64} that $\hat t_{(t)}$ is a compatible function for $H_{(t)},r_0$ with constants $c_-,c_+$. 

	Now we apply \cref{X49}, which, together with \cref{X58}, yields that 
	\begin{equation}\label{X40}
		\Im\biggl(-\frac{w_{H,21}(t;ir)}{w_{H,22}(t;ir)}\biggr)
		= \Im q_{H_{(t)}}(ir)
		\asymp \frac{1}{r\omega_{H_{(t)},2}(a,\hat t_{(t)}(r))},
		\qquad r>r_0,
	\end{equation}
	where the constants in `\/$\asymp$' depend on $c_-,c_+$ but not on $H,t,r_0,\hat t,\hat s$. 

	The remaining task is to express $\omega_{H_{(t)},2}\bigl(a,\hat t_{(t)}(r)\bigr)$ in terms of $H$. 
	Relation \eqref{X111} gives 
	\begin{equation}\label{X182}
		\omega_{H_{(t)},2}(a,s)
		=
		\begin{cases}
			\omega_{H,1}(a+t-s,t) \CAS s\in[a,t],
			\\[1ex]
			\omega_{H,1}(a,t)+(s-t) \CAS s\in(t,\infty),
		\end{cases}		
	\end{equation}
	and plugging in the definition of $\hat t_{(t)}(r)$ we obtain
	\[
		\omega_{H_{(t)},2}\bigl(a,\hat t_{(t)}(r)\bigr)= 
		\begin{cases}
			\omega_{H,1}\bigl(\hat s(t;r),t\bigr) \CAS \hat t(r)\in[a,t],
			\\[2ex]
			\dfrac{1}{\omega_{H,2}(a,t)}\Bigl(\dfrac{c_+}{r^2}+\bigl(\omega_{H,3}(a,t)\bigr)^2\Bigr)
			\CAS \hat t(r)\in(t,\infty).
		\end{cases}
	\]
	Combining this with \eqref{X40} we arrive at \eqref{X57}. 
\end{proof}

\noindent
We are now ready to prove \cref{X1}.

\begin{proof}[Proof of \cref{X1}]
	Let $x\in\dom H$.  It follows from \cref{X2,X52} that
	\begin{align*}
		\log|&w_{H,22}(x;ir)|
		= r\int_a^x\Im\biggl(-\frac{w_{H,21}(t;ir)}{w_{H,22}(t;ir)}\biggr)h_1(t)\DD t
		\\[1ex]
		&\asymp r\int_a^x\bigg[
		\mathds{1}_{[a,\hat t(r))}(t)\frac{\omega_{H,2}(a,t)}{r\bigl[\frac{c_+}{r^2}+(\omega_{H,3}(a,t))^2\bigr]}
		+\frac{\mathds{1}_{[\hat t(r),b)}(t)}{r\omega_{H,1}\bigl(\hat s(t;r),t\bigr)}
		\bigg]h_1(t)\DD t
		\\[1ex]
		&= \int_a^x K_H(t;r)\DD t
	\end{align*}
	for $r>r_0$.
\end{proof}

\subsection[Relating the growth of the fundamental solution to $n_H$]{%
Relating the growth of the fundamental solution to \boldmath{$n_H$}}
\label{X142}

We settle the limit circle case in \cref{X79} by reformulating the statement of \Cref{X1} in terms of $n_H(r)$. 
To achieve this, we use symmetrisation and common tools from complex analysis.

\begin{lemma}\label{X26}
	Let $f$ be an entire function with $f(0)=1$ that is real along the real axis, has only real zeros, and is 
	of bounded type in $\bb C^+$.  Further, denote by $n_f(r)$ the number of zeros of $f$ in the interval $(-r,r)$ 
	\textup{(}counted according to their multiplicities\textup{)}.  Then 
	\[
		\log|f(ir)| = \frac{r^2}{2}\int_0^\infty\frac{1}{t+r^2}\cdot\frac{n_f(\sqrt t)}{t}\DD t.
	\]
\end{lemma}

\begin{proof}
	Let $\lambda_1,\lambda_2,\ldots$ be the sequence of zeros of $f$ arranged according to 
	non-decreasing modulus.  Since $f$ is of bounded type, the genus of $f$ is $0$ or $1$ and 
	\begin{equation}\label{X89}
		f(z) = \lim_{R\to\infty}\prod_{|\lambda_n|<R}\Bigl(1-\frac z{\lambda_n}\Bigr);
	\end{equation}
	see, e.g.\ \cite[V.Theorem~11]{levin:1980}.

	Set $\kappa_n\DE\lambda_n^2$. Then $\sum\limits_n\frac{1}{\kappa_n}<\infty$, 
	and we may consider the canonical product 
	\[
		g(z) \DE \prod_n\Big(1+\frac z{\kappa_n}\Big).
	\]
	Clearly, the counting functions $n_g$ of the zeros of $g$ is related to $n_f$ by 
	\[
		n_g(r) = \#\big\{n\in\bb N\DS \kappa_n<r\big\} = n_f(\sqrt r).
	\]
	Using \eqref{X89} and the fact that $f$ is symmetric w.r.t.\ the real axis we see that 
	\[
		g(r) = f(i\sqrt r)f(-i\sqrt r) = |f(i\sqrt r)|^2.
	\]
	On the other hand, we have 
	\[
		\log g(r)=r\int_0^\infty\frac 1{t+r}\cdot\frac{n_g(t)}t\DD t
	\]
	by \cite[(4.1.4)]{boas:1954}.
\end{proof}

\begin{proof}[Proof of \Cref{X79} \textup{(}limit circle case\textup{)}]
	The spectrum $\sigma_H$ is nothing but the zero set of $w_{H,22}(b;\cdot)$, 
	and $w_{H,22}(b;\cdot)$ is real along the real axis, has only real and simple zeros, 
	and is of bounded type in $\bb C^+$. 
	The above lemma thus shows that the statement \eqref{X65} of \cref{X1} is equivalent to 
	\[
		r^2\int_0^\infty\frac 1{t+r^2}\cdot\frac{n_H(\sqrt t)}t\DD t\asymp\int_a^b K_H(t;r)\DD t,
		\qquad r>r_0,
	\]
	which is \eqref{X90}.
\end{proof}

\subsection{Making a limit argument}
\label{X143}

It remains to make the passage to the limit point case.
The basis of this step is to thoroughly understand the 
relation between the operator models of a Hamiltonian and its truncations. 

Let $H\in\bb H_{a,b}$ and $c\in\domr H$; recall \cref{X43} for the definition of $\domr H$.  Then we have the maps
\begin{align*}
	& \iota_c\FD{L^2(H|_{(a,c)})}{L^2(H)}{5}{[f]_{H|_{(a,c)}}}{\bigg[t\mapsto
	\begin{cases} f(t) \CAS t\in(a,c) \\ 0 \CAS t\in[c,b)\end{cases}\bigg]_H}
	\\[1ex]
	& \rho_c\FD{L^2(H)}{L^2(H|_{(a,c)})}{3}{[f]_H}{[f|_{(a,c)}]_{H|_{(a,c)}.}}
\end{align*}
Clearly, $\iota_c$ is an isometry, $\rho_c$ is a partial isometry with $\ker\rho_c=(\ran\iota_c)^\perp$ and 
$\rho_c\circ\iota_c=\Id$, and $P_c\DE\iota_c\circ\rho_c$ is the orthogonal projection acting as $P_cf=\mathds{1}_{(a,c)}f$. 

The following result that relates the operators $A_H$ and $A_{H|_{(a,c)}}$ is known from the literature. 
In fact, a more general version is shown in \cite[Theorem~3.4]{pruckner.woracek:sinqA}. 
The proof in the presently considered situation is much simpler than in the general case; 
for the convenience of the reader we provide the argument.

\begin{lemma}\label{X21}
	Let $H\in\bb H_{a,b}$ be definite and assume that $\int_a^b h_1(t)\DD t<\infty$. 
	Moreover, let $c\in\domr H$.  Then $A_H$ is invertible and
	\[
		P_cA_H^{-1}P_c\circ\iota_c = \iota_c\circ A_{H|_{(a,c)}}^{-1}.
	\]
\end{lemma}

\begin{proof}
	First note that $A_H$ is invertible by \cref{X16}.
	Let $g\in L^2(H|_{(a,c)})$ and set $\hat g\DE\iota_cg$.  Let $f$ be the unique 
	absolutely continuous representative of $A_{H|_{(a,c)}}^{-1}g$, so that 
	\[
		f' = JHg \text{ a.e.}, \quad (1,0)f(a)=(0,1)f(c)=0.
	\]
	Set 
	\[
		\hat f(t)
		\DE
		\begin{cases} 
			f(t)\CAS t\in[a,c], 
			\\[0.5ex]
			f(c)\CAS t\in(c,b].
		\end{cases}
	\]
	Then $\hat f$ is absolutely continuous and satisfies 
	\begin{align*}
		& \hat f'(t)=JH(t)\hat g(t), \qquad t\in(a,b)\text{ a.e.},
		\\
		& (1,0)\hat f(a)=0, \qquad (0,1)\hat f(t)=0,\ t\in[c,b].
	\end{align*}
	We see that $\hat f\in L^2(H)$ and $A_H^{-1}\hat g=\hat f$. 
	It remains to note that $P_c\circ\iota_c=\iota_c$ and $P_c\hat f=\iota_c f$. 
\end{proof}

\noindent
We need a simple (and well-known) fact about compact operators, which we state explicitly for completeness but otherwise 
skip details. For a compact operator $T$ denote by $s_n(T)$ its $n^{\textup{th}}$ $s$-number, and set
\[
	n(T,\delta) \DE \#\bigl\{n\ge 1\DS s_n(T)>\delta\bigr\}.
\]

\begin{lemma}\label{X23}
	Let $T$ be a compact operator, and let $B_l,B,C_l,C$ be bounded operators 
	such that $\lim_{l\to\infty}B_l=B$ and $\lim_{l\to\infty}C_l^*=C^*$ strongly.  Then 
	\[
		\forall \delta\in(0,\infty)\setminus\big\{s_n(BTC)\DS n\ge1\big\}\DP
		\lim_{l\to\infty} n(B_lTC_l,\delta) = n(BTC,\delta).
	\]
\end{lemma}

\begin{lemma}\label{X24}
	Let $H\in\bb H_{a,b}$ be definite and assume that $\sigma_H$ is discrete 
	and $\int_a^b h_1(t)\DD t<\infty$. 
	\begin{Enumerate}
	\item 
		For each $r>0$ and $c,c'\in\domr H$ with $a<c\le c'<b$ we have 
		\[
			n_{H|_{(a,c)}}(r) \le n_{H|_{(a,c')}}(r).
		\]
	\item 
		If $\sup\domr H=b$, then, for each $r>0$ that is a point of continuity of $n_H$, the relation
		\[
			\lim_{\substack{c\to b \\[0.2ex] c\in\domr H}}n_{H|_{(a,c)}}(r) = n_H(r)
		\]
		holds.
	\item
		If $H$ is in the limit point case and $\sup\domr H<b$, then, for $r>0$ and $c\in[\sup\domr H,b)$,
		\[
			n_{H|_{(a,c)}}(r) = n_H(r).
		\]
	\end{Enumerate}
\end{lemma}

\begin{proof}
	(i)
	Observe that $n_H(r)=n(A_H^{-1},\frac 1r)$, and a similar relation holds for $H|_{(a,c)}$. 
	Now we can apply the above lemmata.  First, by \cref{X21}, 
	\[
		n\big(A_{H|_{(a,c)}}^{-1},\delta\big) = n\big(P_cA_H^{-1}P_c,\delta\big),
	\]
	and the assertion in \textup{(i)} follows since $P_c=P_cP_{c'}$.
	
	(ii)
	The assumption $\sup\domr H=b$ ensures that the union of the ranges of $P_c$ is dense and hence 
	\begin{align}\label{X121}
		\lim_{\substack{c\to b \\[0.2ex] c\in\domr H}}P_c = I
	\end{align}
	strongly.  Now the statement follows from \cref{X23}.
	
	(iii)
	Set $x_0\DE\sup\domr H$ and assume that $x_0<b$.  The interval $(x_0,b)$ must then 
	be an indivisible interval of type $\frac{\pi}{2}$ since $\int_a^b h_1(t)\DD t<\infty$ 
	by assumption.  Let $c\in(x_0,b)$ be arbitrary.
	It follows from \cref{X69} that $A_H$ is unitarily equivalent to $A_{H|_{(a,x_0)}}$,
	and that the operator parts of $A_{H|_{(a,x_0)}}$ and $A_{H|_{(a,c)}}$ are unitarily equivalent.
	From this it follows that $n_H=n_{H|_{(a,x_0)}}=n_{H|_{(a,c)}}$.
\end{proof}

\noindent
We are now ready to finish the proof of \cref{X79}.

\begin{proof}[Proof of \cref{X79} \textup{(}limit point case\textup{)}]
	For $c\in(a,b)$ set
	\begin{align*}
		& r_{0,c}\DE\inf\{r>r_0\DS\hat t(r)<c\},
		\\[0.5ex]
		& \hat t_c\FD{(r_{0,c},\infty)}{(a,c)}{2}{r}{\hat t(r).}
	\end{align*}
	Then, clearly, $\hat t_c$ is a compatible function for $H|_{(a,c)},r_{0,c}$ with constants $c_-,c_+$. 
	Further, set 
	\[
		\hat s_c\FD{\Gamma(\hat t_c)}{[a,c)}{2}{(t,r)}{\hat s(t;r);}
	\]
	then $(\hat t_c,\hat s_c)$ is a compatible pair for $H|_{(a,c)},r_{0,c}$ with constants $c_-,c_+$. 
	Let $K_{H|_{(a,c)}}$ be the function \eqref{X15} constructed with $(\hat t_c,\hat s_c)$
	instead of $(\hat t,\hat s)$.  Then $K_{H|_{(a,c)}}=(K_H)|_{[a,c)\times(r_{0,c},\infty)}$.

	Let $r>r_0$ be fixed.  If $c\in(a,b)$ is sufficiently close to $b$, 
	then $r>r_{0,c}$ and hence $K_{H|_{(a,c)}}(t;r)$ is well defined. 
	Thus we can apply the Monotone Convergence Theorem to conclude that 
	\begin{equation}\label{X76}
		\lim_{c\to b}\int_a^c K_{H|_{(a,c)}}(t;r)\DD t
		= \int_a^b K_H(t;r)\DD t,
	\end{equation}
	where, in the case when $\sup\domr H=b$, one can replace the limit on the left-hand side
	by $\lim_{c\to b,\,c\in\domr H}$.

	It follows from \cref{X24} and the Monotone Convergence Theorem that 
	\[
		\lim_{\substack{c\to b \\[0.2ex] c\in\domr H}}
		\int_0^\infty\frac{1}{t+r^2}\cdot\frac{n_{H|_{(a,c)}}(\sqrt t)}t\DD t
		= \int_0^\infty\frac{1}{t+r^2}\cdot\frac{n_H(\sqrt t)}t\DD t
	\]
	when $\sup\domr H=b$; in the case when $\sup\domr H<b$, the same relation holds with
	the limit on the left-hand side replaced by $\lim_{c\to b}$.
	This, together with \eqref{X76} and the already established limit circle version
	of the theorem, implies \eqref{X90}.
\end{proof}

\section{Growth relative to a regularly varying comparison \\ function}
\label{X133}

In this section we make the transition from the description of the Stieltjes transform 
obtained in \Cref{X79} to that of the function $n_H$ itself.  
There are two types of results; see \cref{X6,X7}.
In the language of complex analysis these are convergence class properties and finite or 
minimal type properties.

\subsection{Functions of regular variation}
\label{X144}

In order to measure the growth of $n_H$, we use regularly varying functions in Karamata's sense as comparison functions. 
This is a very fine scale, and includes all proximate orders in the sense of Valiron. 
In particular, of course, the classical case of comparing growth with powers is covered.

\begin{Definition}\label{X96}
	A function $\ms f\DF[1,\infty)\to(0,\infty)$ is called \emph{regularly varying} (at $\infty$) with 
	\emph{index} $\alpha\in\bb R$ if it is measurable and 
	\begin{equation}\label{X108}
		\forall \lambda\in(0,\infty)\DP \lim_{r\to\infty}\frac{\ms f(\lambda r)}{\ms f(r)} = \lambda^\alpha.
	\end{equation}
	We write $\Ind\ms f$ for the index of a regularly varying function $\ms f$. 
	A regularly varying function with index $0$ is also called \emph{slowly varying}.
\end{Definition}

\noindent
Our standard reference for the theory of such functions is \cite{bingham.goldie.teugels:1989}.

Regularly varying functions $\ms f$ are used to quantify growth for $r\to\infty$, 
and hence the behaviour of $\ms f(r)$ for small $r$ is not controlled by \eqref{X108}. 
Sometimes, we make the additional mild assumption that $\ms f|_Y\asymp1$ for every compact $Y\subset[1,\infty)$.
Note that this is essentially not a restriction since modifying $\ms f$ on a finite interval 
does not affect the validity of condition \eqref{X108} and regularly varying functions 
and their reciprocals are locally bounded for large enough arguments 
by \cite[Corollary~1.4.2]{bingham.goldie.teugels:1989}.
 
Proofs for the following facts can be found in \cite{bingham.goldie.teugels:1989}.

\begin{theorem}[Karamata's Theorem]\label{X95}
	Let $\ms f$ be regularly varying with index $\alpha\in\bb R$ such that $\ms f|_Y \asymp 1$ 
	for every compact set $Y\subseteq[1,\infty)$. 
	\begin{Enumerate}
	\item 
		Assume that $\alpha\ge 0$.  Then the function $t\mapsto\int_1^t \ms f(s)\frac{\D s}{s}$ 
		is regularly varying with index $\alpha$, and 
		\[
			\lim_{t\to\infty}\bigg(
			\raisebox{3pt}{$\ms f(t)$}\Big/\,\raisebox{-2pt}{$\int\limits_1^t \ms f(s)\frac{\D s}{s}$}
			\bigg) 
			= \alpha.
		\]
	\item 
		Assume that $\alpha\le 0$ and $\int_1^\infty \ms f(s)\frac{\D s}{s}<\infty$. 
		Then the function $t\mapsto\int_t^\infty \ms f(s)\frac{\D s}{s}$ is regularly varying 
		with index $\alpha$, and 
		\[
			\lim_{t\to\infty}\bigg(
			\raisebox{3pt}{$\ms f(t)$}\Big/\,\raisebox{-2pt}{$\int\limits_t^\infty \ms f(s)\frac{\D s}{s}$}
			\bigg)
			= -\alpha.
		\]
	\end{Enumerate}
\end{theorem}

\begin{Remark}\label{X206}
	If $\ms f$ is as in \cref{X95} and $\alpha\le0$, then
	\[
		\int_1^t \ms f(s)\frac{\D s}{s} \gg \ms f(t)
	\]
	as $t\to\infty$.  For $\alpha=0$ this follows from \cref{X95}\,(i).
	When $\alpha<0$, then $\ms f(t)\to0$ by \cite[Theorem~1.5.6\,(iii)]{bingham.goldie.teugels:1989}.
\end{Remark}

\noindent
We also need Karamata's theorem for Stieltjes transforms in the form proved 
in \cite[Theorem~A.7]{langer.woracek:kara}.

\begin{theorem}[Karamata's theorem for the Stieltjes transform, direct half]
\label{X91}
	Let $\mu$ be a positive Borel measure on $[0,\infty)$ that
	satisfies $\int_{[0,\infty)}\frac{\D\mu(t)}{1+t}<\infty$.
	If $t\mapsto\mu([0,t))$ is regularly varying with index $\alpha$, then $\alpha\in[0,1]$ and
	\[
		\int_{[0,\infty)}\frac{\D\mu(t)}{t+x}
		\sim
		\frac{\pi\alpha(1-\alpha)}{\sin(\pi\alpha)}\int_x^\infty\frac{\mu([0,t))}{t^2}\DD t
		\sim
		\begin{cases}
			\displaystyle\frac{\pi\alpha}{\sin(\pi\alpha)}\cdot\frac{\mu([0,x))}{x} \CAS \alpha\in[0,1),
			\\[3ex]
			\displaystyle\int_x^\infty\frac{\mu([0,t))}{t^2}\DD t \CAS \alpha=1,
		\end{cases}
	\]
	as $x\to\infty$ where $\frac{\pi\alpha(1-\alpha)}{\sin(\pi\alpha)}$ is understood as $1$ when $\alpha\in\{0,1\}$,
	and $\frac{\pi\alpha}{\sin(\pi\alpha)}$ is understood as $1$ when $\alpha=0$.
\end{theorem}

\noindent
A regularly varying function $\ms f$ with positive index is---at least asymptotically---invertible. 
The following result can be found in \cite[Theorem~1.5.12]{bingham.goldie.teugels:1989}.

\begin{theorem}\label{X100}
	Let $\ms f$ be regularly varying with index $\alpha>0$.  
	Then there exists a regularly varying function $\ms f^-$ with index $\frac 1\alpha$ such that
	\begin{equation}\label{X101}
		(\ms f\circ\ms f^-)(x) \sim (\ms f^-\circ\ms f)(x)\sim x.
	\end{equation}
\end{theorem}

\noindent
Any regularly varying function $\ms f^-$ with the property \eqref{X101} is called an \emph{asymptotic inverse} 
of $\ms f$, and asymptotic inverses are determined uniquely up to $\sim$.

In \cref{X212} we also need the notion of smooth variation.  For the definition
see, e.g.\ \cite[\S1.8]{bingham.goldie.teugels:1989}.

\begin{Definition}\label{X220}
	A positive function $f$ is called \emph{smoothly varying} with index $\alpha$ if 
	it is in $C^\infty$ and $h(x)\DE\log f(e^x)$ satisfies
	\begin{equation}\label{X221}
		h'(x) \to \alpha, \qquad h^{(n)}(x) \to 0, \quad n\in\{2,3,\ldots\},
	\end{equation}
	as $x\to\infty$.
\end{Definition}

\noindent
One can show that \eqref{X221} is equivalent to 
\begin{equation}\label{X222}
	\lim_{x\to\infty}\frac{x^nf^{(n)}(x)}{f(x)} = \alpha(\alpha-1)\cdots(\alpha-n+1),
	\qquad n\in\bb N;
\end{equation}
see \cite[(1.8.1')]{bingham.goldie.teugels:1989}.

The next theorem, \cite[Theorem~1.8.2]{bingham.goldie.teugels:1989}, shows that we can often 
assume, without loss of generality, that a regularly varying function is smoothly varying.

\begin{theorem}[Smooth Variation Theorem]
\label{X74}
	Let $\ms f$ be regularly varying with index $\alpha$.  Then there exists a smoothly varying
	function $\ms g$ with index $\alpha$ such that $f(r)\sim g(r)$ as $r\to\infty$.
	If $\alpha>0$, one can choose $\ms f$ to be strictly increasing.
\end{theorem}

\noindent
A function $\ms f:(0,1]\to(0,\infty)$ is said to be regularly varying (or smoothly varying respectively) at $0$ 
with index $\alpha$ if $\ms g(x)\DE \ms f\bigl(\frac{1}{x}\bigr)$ is regularly varying (or smoothly varying respectively) 
with index $-\alpha$.
For a smoothly varying function $\ms f$ at $0$ with index $\alpha$ the relations in \eqref{X222}
remain valid with $\lim_{x\to\infty}$ replaced by $\lim_{x\to0}$.

\subsection{Conditions for convergence class}
\label{X145}

In the first main theorem of this section we give an explicit criterion for the set $\sigma_H$ to be of 
convergence class with respect to a regularly varying function $\ms g$, i.e.\ we characterise convergence of the series
\[
	\sum_{\lambda\in\sigma_H\setminus(-1,1)}\frac{1}{\ms g(|\lambda|)}.
\]
Here we assume that $\Ind \ms g \in [0,2]$ since our characterisation rests on the a priori assumption 
that $\sum_{\lambda\in\sigma_H\setminus\{0\}}\frac 1{\lambda^2}<\infty$. 
This, however, is no essential restriction since the case of dense spectrum (even with any index strictly larger than $1$) 
was settled in \cite{romanov.woracek:ideal}. 

Note that setting $\ms g(r) \DE r$ yields a criterion for $A_H$ to have trace-class resolvents, 
which is stated explicitly in \cref{X146}.  Moreover, it should be emphasised that the theorem applies also when $\Ind \ms g=0$. 
In particular, we can characterise occurrence of very sparse spectrum on the level of 
(again terminology from complex analysis) logarithmic order or other concepts of growth of order zero. 

Recall the discussion concerning discrete spectrum in \cref{X16} and the notion of compatible pairs from \cref{X11}.

\begin{theorem}\label{X6}
	Let $H\in\bb H_{a,b}$ be definite and assume that $\sigma_H$ is discrete and that $\int_a^b h_1(t)\DD t<\infty$.
	Further, let $r_0\ge 0$, $c_-,c_+>0$, suppose that $(\hat t,\hat s)$ is a compatible pair for $H,r_0$ 
	with constants $c_-,c_+$, and let $K_H(t;r)$ be as in \eqref{X15}. 
	Let $\ms g$ be a regularly varying function with $\Ind\ms g\in[0,2]$ such that $\ms g(t)\ll t^2$ 
	as $t\to\infty$ and $\ms g|_Y \asymp 1$ on every compact set $Y \subseteq [1,\infty)$. 
	\begin{Itemize}
	\item 
		If $\Ind\ms g\in(0,2)$, set $\ms g^*\DE\ms g$.
	\item 
		If $\Ind\ms g=0$, assume that we have a regularly varying function $\ms g^*$ 
		with $\frac{1}{\ms g^*}$ locally integrable and
		\begin{equation}\label{X12}
			\int_t^\infty\int_u^\infty\frac{1}{\ms g^*(s)}\frac{\D s}{s}\frac{\D u}{u} 
			\asymp \frac{1}{\ms g(t)}.
		\end{equation}
	\item 
		If $\Ind\ms g=2$, assume that we have a regularly varying function $\ms g^*$ with 
		\begin{equation}\label{X13}
			\int_1^t\frac{s^2}{\ms g^*(s)}\frac{\D s}{s} \asymp \frac{t^2}{\ms g(t)}.
		\end{equation}
	\end{Itemize}
	Then, with $r_0'\DE\max\{r_0,1\}$, the following equivalence holds:
	\[
		\sum_{\lambda\in\sigma_H\setminus(-1,1)}\frac{1}{\ms g(|\lambda|)} < \infty
		\quad\Longleftrightarrow\quad
		\int_{r_0'}^\infty\frac 1{r\ms g^*(r)}\int_a^b K_H(t;r)\DD t\DD r<\infty.
	\]
\end{theorem}

\bigskip

\noindent
Before we prove the theorem, let us briefly discuss the case when $\Ind\ms g\in\{0,2\}$.
It seems to be unknown whether \eqref{X12} or \eqref{X13}, respectively, 
can always be satisfied, cf.\ \cite[Remark~3.7]{pesa:2024}. 
However, note that if $\ms g^*$ exists, then $\Ind\ms g^*=\Ind\ms g$ but necessarily $\ms g(t)\ll\ms g^*(t)$. 
We give examples where one can find $\ms g^*$ in \cref{X14,X56} below.

To deduce \cref{X6} from \cref{X79}, it is enough to apply the following simple Abelian--Tauberian lemma. 
It is shown by standard methods; for completeness we provide details.

\begin{lemma}\label{X20}
	Let $\ms f$ be regularly varying and locally integrable with 
	\[
		\int_1^\infty\ms f(r)\DD r=\infty, 
		\qquad 
		\int_1^\infty\frac{\ms f(r)}{1+r}\DD r<\infty.
	\]
	Let $\nu$ be a positive Borel measure on $[1,\infty)$, and set $n(t)\DE\nu([1,t))$. 
	Then 
	\[
		\int_1^\infty\ms f(r)\biggl(\int_1^\infty\frac{1}{t+r}\cdot\frac{n(\sqrt t)}t\DD t\biggr)\DD r < \infty
		\quad\Longleftrightarrow\quad
		\int_1^\infty\hat{\ms f}(t)\DD\nu(t) < \infty
	\]
	where 
	\[
		\hat{\ms f}(t)\DE\int_t^\infty\frac 1u\int_u^\infty\frac{1}{s^3}
		\int_1^s \ms f(r^2)r\DD r\DD s\DD u.
	\]
	If $-1<\Ind\ms f<0$, then $\hat{\ms f}(t)\asymp\ms f(t^2)$. 
	If $\Ind\ms f\in\{-1,0\}$, then $\hat{\ms f}(t)\gg\ms f(t^2)$. 
\end{lemma}

\begin{proof}
	Consider the measure $\DD\mu(t)\DE\ms f(t)\DD t$.
	Its distribution function
	\[
		\mu([1,r))=\int_1^r\ms f(t)\DD t
	\]
	is regularly varying by \cref{X95}, and hence its Stieltjes transform satisfies 
	\[
		\int_1^\infty\frac 1{t+r}\cdot\ms f(r)\DD r
		\asymp \int_t^\infty\frac 1{s^2}\int_1^s\ms f(r)\DD r\DD s,
		\qquad t>1,
	\]
	by \cref{X91}.
	We use Fubini's theorem, a change of variables, and integration by parts
	(in the form of \cite[Lemma~A.3]{langer.woracek:kara}) to compute 
	(recall the notation `$\diamond$' from \cref{X47})
	\begin{align*}
		& \int\limits_1^\infty\ms f(r)\bigg( 
		\int\limits_1^\infty\frac 1{t+r}\cdot\frac{n(\sqrt t)}t\DD t\bigg)\DD r
		= \int\limits_1^\infty\bigg(\int\limits_1^\infty\frac 1{t+r}\cdot\ms f(r)\DD r\bigg)\cdot\frac{n(\sqrt t)}t\DD t
		\\[1ex]
		&\diamond
		\int\limits_1^\infty\bigg(\int\limits_t^\infty\frac 1{s^2}\int\limits_1^s\ms f(r)\DD r\DD s\bigg)
		\cdot\frac{n(\sqrt t)}t\DD t
		\\[1ex]
		&\diamond
		\int\limits_1^\infty\frac 1u\int\limits_{u^2}^\infty\frac 1{s^2}\int\limits_1^s\ms f(r)\DD r\DD s\cdot n(u)\DD u
		\\[1ex]
		&= 
		\int\limits_1^\infty\bigg(\int\limits_t^\infty\frac 1u\int\limits_{u^2}^\infty\frac 1{s^2}
		\int\limits_1^s\ms f(r)\DD r\DD s\DD u\bigg)\DD\nu(t).
	\end{align*}
	It remains to evaluate
	\[
		\int\limits_{u^2}^\infty\frac{1}{s^2}\int\limits_1^s\ms f(r)\DD r\DD s
		= 2\int\limits_{u^2}^\infty\frac 1{s^2}\int\limits_1^{\sqrt s}\ms f(p^2)p\DD p\DD s
		= 4\int\limits_u^\infty\frac 1{v^3}\int\limits_1^v\ms f(p^2)p\DD p\DD v.
	\]
	To see the additional statements, note that 
	\begin{alignat*}{3}
		\Ind\ms f &> -1 \quad&&\Rightarrow\quad & \frac{1}{s^2}\int\limits_1^s\ms f(r^2)r\DD r &\asymp \ms f(s^2),
		\\[1ex]
		\Ind\ms f &< 0 \quad&&\Rightarrow\quad & \int\limits_u^\infty\frac{\ms f(s^2)}{s}\DD s &\asymp \ms f(s^2),
	\end{alignat*}
	while $\ms f(s^2)$ is a $\Smallo(\cdot)$ of the respective integral if $\Ind\ms f=-1$ or $\Ind\ms f=0$. 
\end{proof}

\begin{proof}[Proof of \cref{X6}]
	We check that \cref{X20} is applicable with the function 
	\[
		\ms f(r) \DE \frac 1{\ms g^*(\sqrt r)}.
	\]
	Note that $\Ind\ms f=-\frac 12\Ind\ms g$.  If $\Ind\ms g\in[0,2)$, it is clear 
	that $\int_1^\infty\ms f(t)\DD t=\infty$. 
	If $\Ind\ms g=2$, we have 
	\[
		\int_1^{t^2}\ms f(r)\DD r
		= 2\int_1^{t}\frac{s^2}{\ms g^*(s)}\frac{\D s}{s}
		\asymp \frac{t^2}{\ms g(t)} \to \infty
	\]
	by assumption.
	If $\Ind\ms g\in(0,2]$, we clearly have $\int_1^\infty\frac{\ms f(r)}{r}\DD r<\infty$.  If $\Ind\ms g=0$, then
	\[
		\int_1^\infty\int_{u^2}^\infty\frac{\ms f(r)}{r}\DD r\frac{\D u}{u}
		= 2\int_1^\infty\int_u^\infty\frac{1}{\ms g^*(s)}\frac{\D s}{s}\frac{\D u}{u}
		< \infty,
	\]
	and thus $\int_1^\infty\frac{\ms f(r)}{r}\DD r<\infty$. 
	
	Next, we identify $\hat{\ms f}$ from \cref{X20}.  If $\Ind\ms g\in(0,2)$, it is already stated in the lemma that 
	$\hat{\ms f}(t)\asymp\ms f(t^2)\asymp \frac{1}{\ms g(t)}$.  
	If $\Ind\ms g=0$ or $\Ind \ms g=2$, we get rid of the integrals in the definition of $\hat{\ms f}$ 
	by using Karamata's theorem (\cref{X95}) and one of \eqref{X12}, \eqref{X13}.
	Also in this case we arrive at $\hat{\ms f}(t)\asymp\frac{1}{\ms g(t)}$. 
	
	Hence we can use \cref{X20,X79} to obtain
	\begin{align*}
		\sum_{\lambda\in\sigma_H\setminus(-1,1)}\frac{1}{\ms g(|\lambda|)}
		&\diamond \int\limits_1^\infty\hat{\ms f}(t)\DD n_H(t)
		\diamond \int\limits_1^\infty\frac{1}{\ms g^*(\sqrt r)}
		\int\limits_1^\infty\frac{1}{t+r}\cdot\frac{n_H(\sqrt t)}t\DD t\DD r
		\\
		&\diamond \int\limits_1^\infty\frac{r}{\ms g^*(r)}
		\int\limits_1^\infty\frac{1}{t+r^2}\cdot\frac{n_H(\sqrt t)}t\DD t\DD r
		\\
		&\diamond \int\limits_{r_0'}^\infty\frac{r}{\ms g^*(r)}
		\int\limits_0^\infty\frac{1}{t+r^2}\cdot\frac{n_H(\sqrt t)}t\DD t\DD r
		\\
		&\diamond \int\limits_{r_0'}^\infty\frac{1}{r\ms g^*(r)}\int\limits_a^b K_H(t;r)\DD t\DD r,
	\end{align*}
	which finishes the proof.
\end{proof}

\noindent
Let us now discuss two classes of examples where one can find $\ms g^*$ when $\Ind\ms g\in\{0,2\}$.
In \cref{X14} below we show that \eqref{X12} and \eqref{X13} can be satisfied at least 
for all $\ms g$ that behave for $r\to\infty$ like
\begin{equation}\label{X94}
	r^\alpha\cdot\bigl(\log r\bigr)^{\beta_1}\cdot\bigl(\log\log r\bigr)^{\beta_2}
	\cdot\ldots\cdot
	\bigl(\underbrace{\log\cdots\log}_{\text{\footnotesize$N$\textsuperscript{th} iterate}}r\bigr)^{\beta_N},
\end{equation}
where $\alpha \in\bb R$ and $\beta_1,\ldots,\beta_N\in\bb R$.
In fact, we explicitly specify a suitable function $\ms g^*$; see \eqref{X28} and \eqref{X29}.
It is noteworthy, and somewhat unexpected, that the cases $\Ind\ms g=0$ and 
$\Ind\ms g=2$ behave differently when it comes to the gap between $\ms g$ and $\ms g^*$. 

\begin{Example}\label{X14}
	Let us start with some general considerations.
	To this end, let $N\in\bb N$ and $\gamma_1,\ldots,\gamma_N\in\bb R$, 
	and consider the slowly varying function (for $r$ sufficiently large)
	\begin{equation}\label{X77}
		\ms l(r) \DE \prod_{k=1}^N\bigl(\log^{[k]}r\bigr)^{\gamma_k},
	\end{equation}
	where $\log^{[k]}$ is the $k^{\textup{th}}$ iterate of the function $\log$.
	We assume in the following that 
	\begin{equation}\label{X109}
		\gamma_k=0, \quad k<n; \qquad \gamma_n\ne0
	\end{equation}
	with some $n\in\{1,\ldots,N\}$.
	Further, denote by $D$ the operator 
	\begin{equation}\label{X85}
		(Df)(r)\DE rf'(r).
	\end{equation}
	A computation using induction with respect to $m$ shows that, for each $m\in\bb N$,
	\begin{multline}\label{X84}
		(D^m\ms l)(r) \in \frac{\ell(r)}{(\log r)^m\prod_{j=2}^n\log^{[j]}r}
		\biggl[\gamma_n\prod_{j=1}^{m-1}(\gamma_1-j)
		\\[1ex]
		+ \Span\biggl\{\prod_{l=2}^N\bigl(\log^{[l]}r)^{\nu_l}\DS
		\nu_l\le0,\,\nu_l\in\mathbb{Z},\,(\nu_2,\ldots,\nu_N)\ne(0,\ldots,0)\biggr\}\biggr].
	\end{multline}
	Let us now consider the two cases $\Ind\ms g=0$ and $\Ind\ms g=2$ separately.
	\begin{Enumerate}
	\item 
		Assume that
		\[
			\ms g(r) = \prod_{k=n}^N\bigl(\log^{[k]}r\bigr)^{\beta_k}
		\]
		for large enough $r$ with $n,N\in\bb N$, $n\le N$ and $\beta_n>0$.
		Set $\ms l\DE\frac{1}{\ms g}$, which is of the form \eqref{X77}
		with $\gamma_n=-\beta_n<0$, and hence $\ms l(r)\to0$ as $r\to\infty$.
		It follows from \eqref{X84} that
		\[
			(D^2\ms l)(r) \sim \gamma_n(\gamma_1-1)\frac{\ms l(r)}{(\log r)^2\prod_{j=2}^n\log^{[j]}r}
		\]
		with $\gamma_n(\gamma_1-1)>0$.  Therefore the function
		\begin{equation}\label{X28}
			\ms g^*(r) \DE \ms g(r)\cdot(\log r)^2\prod_{j=2}^n\log^{[j]}r
		\end{equation}
		satisfies
		\begin{align*}
			\int_t^\infty\int_u^\infty\frac{1}{\ms g^*(s)}\frac{\D s}{s}\frac{\D u}{u}
			&= \int_t^\infty\int_u^\infty\frac{\ms l(s)}{(\log s)^2\prod_{j=2}^n\log^{[j]}s}\frac{\D s}{s}\frac{\D u}{u}
			\\[1ex]
			&\asymp\int_t^\infty\int_u^\infty(D^2\ms l)(s)\frac{\D s}{s}\frac{\D u}{u}
			= \ms l(t),
		\end{align*}
		which shows that \eqref{X12} holds;
		note that in the last equality we have also used that $(D\ms l)(r)\to0$ as $r\to\infty$,
		which also follows from \eqref{X84}.
	\item
		Let us now assume that
		\[
			\ms g(r) = r^2\prod_{k=n}^N\bigl(\log^{[k]}r\bigr)^{\beta_k}
		\]
		for large enough $r$ with $n,N\in\bb N$, $n\le N$ and $\beta_n<0$.
		Set $\ms l(r)\DE\frac{r^2}{\ms g(r)}$, which is of the form \eqref{X77} with $\gamma_n=-\beta_n>0$.
		It follows from \eqref{X84} that
		\[
			(D\ms l)(r) \sim \gamma_n\frac{\ms l(r)}{\prod_{j=1}^n\log^{[j]}r}.
		\]
		Hence the function
		\begin{equation}\label{X29}
			\ms g^*(r) \DE \ms g(r)\cdot\prod_{j=1}^n\log^{[j]}r
		\end{equation}
		satisfies (with $t_0$ sufficiently large)
		\[
			\int_{t_0}^t\frac{s^2}{\ms g^*(s)}\frac{\D s}{s}
			= \int_{t_0}^t \frac{\ms l(s)}{\prod_{j=1}^n\log^{[j]}s}\cdot\frac{\D s}{s}
			\asymp \int_{t_0}^t(D\ms l)(s)\frac{\D s}{s}
			\asymp \ms l(t) = \frac{t^2}{\ms g(t)},
		\]
		which shows that \eqref{X13} holds.
	\end{Enumerate}
\end{Example}

\noindent
As a second class of examples we consider slowly varying functions that are between logarithms
and power functions with positive exponent.  For these functions the gap that
appears between $\ms g$ and $\ms g^*$ gradually closes.

\begin{Example}\label{X56}
	Let $\gamma\in(0,1)$ and consider the function $\ms g(r)=e^{(\log r)^\gamma}$, 
	which is slowly varying but grows faster than any power of $\log r$.
	Set $\ms l(r)\DE \frac{1}{\ms g(r)}=e^{-(\log r)^\gamma}$.
	It is easy to see that, with the operator $D$ from \eqref{X85},
	we have
	\[
		(D^2\ms l)(r) \sim \gamma^2e^{-(\log r)^\gamma}(\log r)^{2(\gamma-1)}.
	\]
	Hence the function
	\begin{equation}\label{X86}
		\ms g^*(r) \DE \ms g(r)(\log r)^{2(1-\gamma)}
	\end{equation}
	satisfies
	\[
		\int_t^\infty\int_u^\infty\frac{1}{\ms g^*(s)}\frac{\D s}{s}\frac{\D u}{u}
		\asymp \int_t^\infty\int_u^\infty(D^2\ms l)(s)\frac{\D s}{s}\frac{\D u}{u}
		= \ms l(t) = \frac{1}{\ms g(t)},
	\]
	and therefore \eqref{X12} holds.
	Note that the exponent in the last factor on the right-hand side of \eqref{X86}
	tends to $0$ as $\gamma\to1$.
\end{Example}

\subsection{Trace class resolvents}
\label{X146}

\Cref{X6} gives, in particular, a criterion for resolvents of $A_H$ to belong to the trace class $\mathfrak{S}_1$, 
namely by setting $\ms g(r) \DE r$.
We also obtain an explicit expression for the trace of the inverse if the latter
is a trace class operator.
Note that under our standard assumption \eqref{X154} we have $0\notin\sigma_H$;
see \cref{X16}.

\begin{theorem}\label{X53}
	Let $H\in\bb H_{a,b}$ be definite and assume that $\sigma_H$ is discrete and that $\int_a^b h_1(t)\DD t<\infty$.
	Further, let $r_0\ge 0$, $c_-,c_+>0$, suppose that $(\hat t,\hat s)$ is a compatible pair for $H,r_0$ 
	with constants $c_-,c_+$, and let $K_H(t;r)$ be as in \eqref{X15}. 
	Then
	\begin{equation}\label{X199}
		A_H^{-1} \in \mathfrak{S}_1
		\quad\Longleftrightarrow\quad
		\int_1^\infty\frac 1{r^2}\int_a^b K_H(t;r)\DD t\DD r<\infty.
	\end{equation}
	Moreover, if $A_H^{-1}\in\mf S_1$, then
	\begin{equation}\label{X66}
		\tr\bigl(A_H^{-1}\bigr) = -\lim_{t\to b}\int_a^t h_3(s)\DD s.
	\end{equation}
\end{theorem}

\begin{proof}
	First note that \eqref{X199} follows directly from \cref{X6} with $\ms g(r)=r$.
	
	Let us now come to the proof of \eqref{X66}.
	Let $c\in\domr H$.
	From \cref{X21} we obtain
	\begin{equation}\label{X119}
		\sum_{\lambda\in\sigma_{H|_{(a,c)}}}\frac{1}{|\lambda|} 
		\le \sum_{\lambda\in\sigma_H}\frac{1}{|\lambda|} < \infty.
	\end{equation}
	Thus Hadamard's theorem yields the representation
	\[
		w_{H,22}(c;z) = \prod_{\lambda\in\sigma_{H|_{(a,c)}}}\Bigl(1-\frac{z}{\lambda}\Bigr),
	\]
	which leads to
	\begin{equation}\label{X30}
		\frac{\partial}{\partial z}\Bigl(w_{H,22}(c;z)\Bigr)\Big|_{z=0}
		= \frac{\partial}{\partial z}\Bigl[\log\Bigl(w_{H,22}(c;z)\Bigr)\Bigr]\Big|_{z=0}
		= -\sum_{\lambda\in\sigma_{H|_{(a,c)}} }\frac{1}{\lambda}.
	\end{equation}
	On the other hand, the fundamental solution $W_H(t;z)$ can be expanded into a power series
	(see, e.g.\ \cite[Section~2.2]{langer.pruckner.woracek:heniest}):
	\[
		W_H(t;z) = \sum_{n=0}^\infty W_n(t)z^n
	\]
	with
	\[
		W_0(t) = I, \qquad W_{n+1}(t) = \int_a^t W_n(s)H(s)(-J)\DD s,
	\]
	which gives
	\begin{equation}
	\label{X31}
		\frac{\partial}{\partial z}W_H(c;z)\big|_{z=0} = W_1(c) = \Omega_H(a,c)(-J).
	\end{equation}
	Combining \eqref{X30} and \eqref{X31} we arrive at
	\begin{equation}
	\label{X120}
		\tr (A_{H|_{(a,c)}}^{-1}) = \sum_{\lambda\in\sigma_{H|_{(a,c)}}}\frac{1}{\lambda}
		= -\omega_{H,3}(a,c).
	\end{equation}
	When $H$ is in the limit circle case, we choose $c=b$ to finish the proof.
	When $H$ is in the limit point case and $\sup\domr H<b$, then choose $c=\sup\domr H$
	and use \cref{X69}\,(i).
	
	It remains to consider the case when $H$ is in the limit point case and $\sup\domr H=b$.
	\Cref{X21} implies that
	\[
		\tr (A_{H|_{(a,c)}}^{-1}) = \tr(P_cA_H^{-1}P_c) = \tr(P_cA_H^{-1})
	\]
	and hence
	\[
		\tr(A_H^{-1})-\tr(A_{H|_{(a,c)}}^{-1}) = \tr\bigl((I-P_c)A_H^{-1}\bigr).
	\]
	By \eqref{X121}, we can apply \cref{X23} with $B_l=I-P_l$, $B=0$, $T=A_H^{-1}$, and $C_l=C=I$,
	which yields that, for every $n$, the numbers $s_n ((I-P_c)A_H^{-1})$ converge to zero 
	as $c\to b$, $c\in\domr H$. 
	Since $s_n\big((I-P_c)A_H^{-1}\big) \le s_n (A_H^{-1})$ for every $n$, 
	the dominated convergence theorem implies that
	\begin{align*}
		\limsup_{\substack{c\to b \\[0.2ex] c\in\domr H}}\big|\tr (A_H^{-1})-\tr (A_{H|_{(a,c)}}^{-1})\big|
		\le \lim_{\substack{c\to b \\[0.2ex] c\in\domr H}}\sum_{n\in\bb N}s_n\bigl((I-P_c)A_H^{-1}\bigr) 
		= 0.
	\end{align*}
	Taking the limit $c \to b$, $c\in\domr H$, in \eqref{X120} we obtain \eqref{X66}.
\end{proof}

\begin{remark}\label{X116}
	The method in the proof of \cref{X53} can be used to compute
	$\sum_{\lambda\in\sigma_H}\frac{1}{\lambda^p}$ for any integer $p>1$.
	Assume that $A_H^{-1}\in\mathfrak{S}_p$ (a criterion for the latter is stated in \cref{X16}).
	Then we can write
	\[
		w_{H,22}(c;z) = \prod_{\lambda\in\sigma_{H|_{(a,c)}}}\Bigl(1-\frac{z}{\lambda}\Bigr)
		e^{\frac{z}{\lambda}+\frac12(\frac{z}{\lambda})^2+\cdots+\frac{1}{p-1}(\frac{z}{\lambda})^{p-1}}, 
		\qquad z\in\bb C,
	\]
	by Hadamard's theorem.  Since
	\[
		\log\biggl[\Bigl(1-\frac{z}{\lambda}\Bigr)
		e^{\frac{z}{\lambda}+\frac12(\frac{z}{\lambda})^2+\cdots+\frac{1}{p-1}(\frac{z}{\lambda})^{p-1}}\biggr]
		= -\sum_{k=p}^{\infty}\frac{1}{k}\Bigl(\frac{z}{\lambda}\Bigr)^k
	\]
	for $|z|<|\lambda|$, we obtain
	\[
		\frac{\partial^p}{\partial z^p}\log w_{H,22}(c;z)\Big|_{z=0} 
		= -(p-1)!\sum_{\lambda\in\sigma_{H|_{(a,c)}}}\frac{1}{\lambda^{p}}.
	\]
	This expression can be rewritten in terms of iterated integrals over entries of $H$
	using the series representation of $W_H(t;z)$.  
	For $p=2$ we re-obtain the formula for the Hilbert--Schmidt norm given in \cite{kaltenbaeck.woracek:hskansys}:
	\[
		\sum_{\lambda\in\sigma_H}\frac{1}{\lambda^2} = 2\int_a^b \omega_{H,2}(a,t)h_1(t) \DD t;
	\]
	for $p=3$ we get
	\[
		\sum_{\lambda\in\sigma_H}\frac{1}{\lambda^3}
		= 12\lim_{t\to b}\int_a^t\int_a^s\omega_{H,2}(a,x)h_3(x)\DD x\,h_1(s)\DD s.
	\]
\end{remark}

\subsection{Conditions for finite and minimal type}
\label{X147}

Our second way of describing the density of eigenvalues is by comparing $n_H(r)$ with a regularly varying function $\ms g$. 
Contrasting \cref{X6}, it amounts to a characterisation only if $\Ind\ms g\in(0,2)$; 
note that, in the formulation of the theorem, $\ms g(r)\ll\ms g_*(r)$ if $\Ind\ms g\in\{0,2\}$.

\begin{theorem}\label{X7}
	Let $H\in\bb H_{a,b}$ be definite, assume that $\int_a^b h_1(t)\DD t<\infty$, and let $n_H$ be as in \eqref{X81}.
	Further, let $r_0\ge 0$, $c_-,c_+>0$, assume that $(\hat t,\hat s)$ is a compatible pair for $H,r_0$ 
	with constants $c_-,c_+$, and let $K_H(t;r)$ be as in \eqref{X15}. 
	\begin{Enumerate}
	\item
		We have
		\begin{equation}\label{X155}
			n_H(r) \lesssim \int_a^b K_H(t;r)\DD t, \qquad r>r_0,
		\end{equation}
		where the constant implicit in $\lesssim$ depends only on $c_-,c_+$
		but not on $H,r_0,\hat t,\hat s$.
	\item
		Let $\ms g$ be a regularly varying function with $\alpha\DE\Ind\ms g\ge 0$ 
		such that $\int_1^\infty\frac{\ms g(t)}{t^3}\DD t<\infty$, and set 
		\begin{equation}\label{X185}
			\ms g_*(r)\DE
			\begin{cases}
				\int_1^r\frac{\ms g(t)}t\DD t \CAS \alpha=0,
				\\[1ex]
				\frac{\frac{\pi}{2}}{\sin\frac{\pi\alpha}{2}}\ms g(r) \CAS \alpha\in(0,2),
				\\[2ex]
				r^2\int_r^\infty\frac{\ms g(t)}{t^3}\DD t \CAS \alpha=2.
			\end{cases}
		\end{equation}
		Then there exists $\Delta_->0$, which depends on $c_-,c_+$ 
		but not on $H,r_0,\hat t,\hat s,\ms g,\alpha$, such that
		\begin{equation}\label{X87}
			\limsup_{r\to\infty}\frac{\Delta_-}{\ms g_*(r)}\int_a^b K_H(t;r)\DD t
			\le \limsup_{r\to\infty}\frac{n_H(r)}{\ms g(r)}.
		\end{equation}
	\end{Enumerate}
\end{theorem}

\noindent
As an example let us again consider a comparison function $\ms g$ of the form
\[
    \ms g(r) = r^\alpha\prod_{k=1}^N(\log^{[k]}r)^{\beta_k}
\]
with $\alpha=\Ind\ms g\in\{0,2\}$ and assume that $\beta_k=-1$ for $k<n$
and $\beta_n\ne-1$.
If $\Ind\ms g=0$ and $\ms g(r)\to\infty$ as $r\to\infty$, then $\beta_1\ge0$ and hence $n=1$.
We can use the computation from \cref{X14} to obtain
\[
    \ms g_*(r) \asymp \ms g(r)\cdot\prod_{j=1}^n\log^{[j]}r.
\]
Contrasting the situation of \cref{X6}, here the boundary cases produce the same gap
between $\ms g$ and $\ms g_*$.

Before we prove \cref{X7}, we need a lemma that provides estimates for Stieltjes transforms
as in \cref{X79}.

\begin{lemma}\label{X233}
	Let $\eta:[0,\infty)\to[0,\infty)$ be a non-decreasing function that vanishes on a neighbourhood of $0$,
	and set
	\[
		v(r) \DE \frac12\int_0^\infty\frac{r^2}{t+r^2}\cdot\frac{\eta(\sqrt t)}t\DD t, \qquad r>0.
	\]
	\begin{Enumerate}
	\item
		We have
		\[
			\eta(r) \le \frac{2}{\log 2}v(r), \qquad r>0.
		\]
	\item
		Let $\ms g$ be a regularly varying function with $\alpha\DE\Ind\ms g\ge 0$ 
		such that $\int_1^\infty\frac{\ms g(t)}{t^3}\DD t<\infty$, and 
		let $\ms g_*$ be as in \eqref{X185}.
		Then
		\begin{equation}\label{X234}
			\limsup_{r\to\infty}\frac{v(r)}{\ms g_*(r)}
			\le \limsup_{r\to\infty}\frac{\eta(r)}{\ms g(r)}.
		\end{equation}
	\end{Enumerate}
\end{lemma}

\begin{proof}
	(i)
	Let $r>0$.  We use the monotonicity of $\eta$ to obtain
	\begin{align*}
		v(r) &\ge \frac12\int_{r^2}^\infty\frac{r^2}{t+r^2}\cdot\frac{\eta(\sqrt t)}t\DD t
		\\[0.5ex]
		&\ge \frac{\eta(r)}{2}\int_{r^2}^\infty\Bigl(\frac{1}{t}-\frac{1}{t+r^2}\Bigr)\DD t
		= \frac{\eta(r)}{2}\cdot\log\frac{t}{t+r^2}\Big|_{r^2}^\infty
		= \frac{\log 2}{2}\cdot\eta(r).
	\end{align*}

	(ii)
	The argument for the proof of \eqref{X87} is more involved. Set 
	$x_0\DE\sup\{t\ge 0\DS\eta(t)=0\}$, which satisfies $x_0>0$ by assumption.
	If the right-hand side of \eqref{X234} is infinite, there is nothing to prove. 
	Hence assume that $\limsup_{r\to\infty}\frac{\eta(r)}{\ms g(r)}<\infty$, choose
	\begin{equation}\label{X230}
		\gamma > \limsup_{r\to\infty}\frac{\eta(r)}{\ms g(r)}, 
	\end{equation}
	and fix $x>\max\{1,\sqrt{x_0}\}$ such that 
	\begin{equation}\label{X232}
		\forall t\ge x\DP \eta(t) \le \gamma\ms g(t).
	\end{equation}
	For $r>0$, we combine the two inequalities
	\begin{align*}
		\int_{x_0}^{x^2}\frac{r^2}{t+r^2}\cdot\frac{\eta(\sqrt t)}{t}\DD t
		&\le \frac{x^2\eta(x)}{x_0},
		\\[1ex]
		\int_{x^2}^\infty\frac{r^2}{t+r^2}\cdot\frac{\gamma\ms g(\sqrt t)}{t}\DD t
		&\ge \eta(x)\int_{x^2}^\infty\frac{r^2}{t+r^2}\cdot\frac{1}{t}\DD t
		= \eta(x)\log\Bigl(1+\frac{r^2}{x^2}\Bigr)
	\end{align*}
	and use \eqref{X232} again to obtain
	\begin{align}
		v(r)
		&= \frac12\int_{x_0}^{x^2}\frac{r^2}{t+r^2}\cdot\frac{\eta(\sqrt t)}t\DD t
		+ \frac12\int_{x^2}^\infty\frac{r^2}{t+r^2}\cdot\frac{\eta(\sqrt t)}t\DD t
		\nonumber
		\\[1ex]
		&\le \frac{x^2}{2x_0\log\bigl(1+\frac{r^2}{x^2}\bigr)}
		\int_{x^2}^\infty\frac{r^2}{t+r^2}\cdot\frac{\gamma\ms g(\sqrt t)}{t}\DD t
		+ \frac12\int_{x^2}^\infty\frac{r^2}{t+r^2}\cdot\frac{\gamma\ms g(\sqrt t)}t\DD t
		\nonumber
		\\[1ex]
		&\le \gamma\cdot\frac12\int_1^\infty\frac{r^2}{t+r^2}\cdot\frac{\ms g(\sqrt t)}{t}\DD t
		\cdot\biggl(1+\frac{x^2}{x_0\log\bigl(1+\frac{r^2}{x^2}\bigr)}\biggr).
		\label{X231}
	\end{align}
	Since the fraction with the logarithm on the right-hand side tends to $0$ as $r\to\infty$, 
	we only have to study the asymptotic behaviour of the integral on the right-hand side.
	
	To this end define the measure $\mu$ on $[0,\infty)$ by $\mu([0,1])=0$ 
	and $\mu([1,t))=\int_1^t\frac{\ms g(\sqrt s)}{s}\DD s$, $t>1$.
	Since $t\mapsto\mu([0,t))$ is regularly varying with index $\frac{\alpha}{2}$ by \cref{X95}\,(i),
	we can apply \cref{X91}, which, for $\alpha=2$, yields
	\[
		\frac12\int_1^\infty\frac{r^2}{t+r^2}\cdot\frac{\ms g(\sqrt t)}{t}\DD t
		\sim \frac12 r^2\int_{r^2}^\infty\frac{\ms g(\sqrt{t})}{t^2}\DD t
		= r^2\int_r^\infty\frac{\ms g(s)}{s^3}\DD s
		= \ms g_*(r).
	\]
	When $\alpha\in[0,2)$, we obtain from \cref{X91} that
	\begin{align*}
		\frac12\int_1^\infty\frac{r^2}{t+r^2}\cdot\frac{\ms g(\sqrt t)}{t}\DD t
		&= \frac{r^2}{2}\int_1^\infty\frac{1}{t+r^2}\DD\mu(t)
		\sim \frac{r^2}{2}\cdot\frac{\frac{\pi\alpha}{2}}{\sin(\frac{\pi\alpha}{2})}\frac{\mu([0,r^2))}{r^2}
		\\[1ex]
		&= \frac12\cdot\frac{\frac{\pi\alpha}{2}}{\sin(\frac{\pi\alpha}{2})}\int_1^{r^2}\frac{\ms g(\sqrt{t})}{t}\DD t
		= \frac{\frac{\pi\alpha}{2}}{\sin(\frac{\pi\alpha}{2})}\int_1^r\frac{\ms g(s)}{s}\DD s
		\sim \ms g_*(r)
	\end{align*}
	in the last step we used \cref{X95} when $\alpha\in(0,2)$ whereas there is equality when $\alpha=0$.
	It now follows from \eqref{X231} that $\limsup_{r\to\infty}\frac{v(r)}{\ms g_*(r)}\le\gamma$.
	Since $\gamma$ was arbitrary satisfying \eqref{X230}, the inequality \eqref{X234} follows.
\end{proof}	

\begin{Remark}\label{X235}
	Let $f$ be an entire function as in \cref{X26} and denote by $n_f(r)$ the number of
	zeros of $f$ in $(-r,r)$.
	Using \cref{X26} and \cref{X233} with $\eta(r)=n_f(r)$ one obtains estimates
	between $\log|f(ir)|$ and $n_f(r)$, e.g.\
	\begin{equation}\label{X228}
		n_f(r) \le \frac{2}{\log 2}\log|f(ir)|, \qquad r>0.
	\end{equation}
\end{Remark}

\begin{proof}[Proof of \cref{X7}]
	The statements of \cref{X7} follow directly from \cref{X79} and \cref{X233} applied with $\eta(r)=n_H(r)$.
\end{proof}

\noindent
Using \cref{X7} we obtain a corollary for the limit circle case. 
In \cref{X1} we saw that the integral over $K_H(t;r)$ governs the growth of the right lower entry $w_{H,22}(b;z)$ 
of the monodromy matrix along the imaginary axis. 
Using machinery from complex analysis we can pass to the growth of $\max_{|z|=r}\|W_H(b;z)\|$. 

\begin{proposition}\label{X45}
	Let $H\in\bb H_{a,b}$ be definite and in the limit circle case, and let $\ms g$ be a regularly varying function 
	with index $\alpha$ such that $\int_1^\infty\frac{\ms g(t)}{t^2}\DD t<\infty$ and $1\ll\ms g(r)\ll r$ as $r\to\infty$.
	Set 
	\begin{equation}\label{X237}
		\ms g^{**}(r)
		\DE
		\begin{cases}
			\int_1^r\frac{\ms g(t)}{t}\DD t \CAS \alpha=0,
			\\[1ex]
			\frac{1}{\alpha(1-\alpha)}\ms g(r) \CAS \alpha\in(0,1),
			\\[1.5ex]
			r\int_r^\infty\frac{\ms g(t)}{t^2}\DD t \CAS \alpha=1.
		\end{cases}
	\end{equation}
	If $\alpha\in(0,1]$, then 
	\begin{equation}\label{X88}
		\limsup_{r\to\infty}\frac{1}{\ms g^{**}(r)}\log\Bigl(\max_{|z|=r}\|W_H(b;z)\|\Bigr)
		\le \Delta_+\limsup_{r\to\infty}\frac{1}{\ms g(r)}\int_a^b K_H(t;r)\DD t
	\end{equation}
	where $\Delta_+$ depends only on $c_-,c_+$ but not on $H,r_0,\hat t,\hat s,\ms g,\alpha$.
	\\
	If $\alpha=0$ and $\int_a^b K_H(t;r)\DD t\lesssim\ms g(r)$, then
	\[
		\limsup_{r\to\infty}\frac{1}{\ms g^{**}(r)}\log\Bigl(\max_{|z|=r}\|W_H(b;z)\|\Bigr) = 0.
	\]
\end{proposition}

\medskip

\noindent
Before we prove \cref{X45}, we need a lemma about entire functions as in \cref{X26}.

\begin{lemma}\label{X236}
	Let $f$ be an entire function with $f(0)=1$ that is real along the real axis, has only real zeros, and is 
	of bounded type in $\bb C^+$.  Denote by $n_f(r)$ the number of zeros of $f$ in the interval $(-r,r)$ 
	\textup{(}counted according to their multiplicities\textup{)},
	and let $\ms g$ be a regularly varying function with index $\alpha$ such that $\int_1^\infty\frac{\ms g(t)}{t^2}\DD t<\infty$ 
	and $1\ll\ms g(r)\ll r$ as $r\to\infty$.  Moreover, assume that 
	\begin{equation}\label{X238}
		\log|f(ir)| \lesssim \ms g(r).
	\end{equation}
	Then
	\begin{equation}\label{X239}
		\limsup_{r\to\infty}\frac{\log\Bigl(\max\limits_{|z|=r}|f(z)|\Bigr)}{\ms g^{**}(r)}
		\le
		\begin{cases}
			\displaystyle \limsup_{r\to\infty}\frac{n_f(r)}{\ms g(r)} \CAS \alpha\in(0,1],
			\\[3ex]
			0 \CAS \alpha=0,
		\end{cases}
	\end{equation}
	where $\ms g^{**}$ is as in \eqref{X237}.
\end{lemma}

\begin{proof}
	The relations \eqref{X228} and \eqref{X238} imply that 
	\begin{equation}\label{X124}
		n_f(r)\lesssim\ms g(r),
	\end{equation}
	and hence
	\[
		\sum_{x\in\bb R: f(x)=0}\frac{1}{|x|} = \int_0^\infty\frac{n_f(t)}{t^2}\DD t < \infty.
	\]
	Using \cite[Theorem~11 in \S V.4]{levin:1980} we obtain 
	\[
		f(z) = \lim_{R\to\infty}\prod_{|x|<R,f(x)=0}\Bigl(1-\frac zx\Bigr)
		= \prod_{x\in\bb R:f(x)=0}\Bigl(1-\frac zx\Bigr).
	\]
	Let $\gamma$ be arbitrary such that
	\begin{equation}\label{X229}
		\gamma > \limsup_{r\to\infty}\frac{n_f(r)}{\ms g(r)},
	\end{equation}
	and choose $R>0$ such that $n_f(r)\le\gamma\ms g(r)$ for all $r\ge R$.
	We obtain from \cite[Lemma~3 in \S I.4]{levin:1980} that, for $z\in\bb C$ with $|z|=r\ge R$,
	\begin{align}
		\log|f(z)| &\le \int_0^r\frac{n_f(t)}{t}\DD t + r\int_r^\infty\frac{n_f(t)}{t^2}\DD t
		\nonumber \\[1ex]
		&\le \int_0^R\frac{n_f(t)}{t}\DD t 
		+ \gamma\biggl[\int_R^r\frac{\ms g(t)}{t}\DD t+r\int_r^\infty\frac{\ms g(t)}{t^2}\DD t\biggr].
		\label{X75}
	\end{align}
	\Cref{X95} implies that
	\[
		\int_R^r\frac{\ms g(t)}{t}\DD t \;
		\begin{cases}
			\sim\frac{1}{\alpha}\ms g(r) \CAS \alpha>0,
			\\[1ex]
			\gg \ms g(r) \CAS \alpha=0,
		\end{cases}
		\hspace*{10ex}
		r\int_r^\infty\frac{\ms g(t)}{t^2}\DD t \;
		\begin{cases}
			\sim\frac{1}{1-\alpha}\ms g(r) \CAS \alpha<1,
			\\[1ex]
			\gg \ms g(r) \CAS \alpha=1,
		\end{cases}
	\]
	and hence the expression within the square brackets on the right-hand side of \eqref{X75}
	is asymptotically equal to $\ms g^{**}(r)$.
	Since $\gamma$ was arbitrary satisfying \eqref{X229}, we obtain
	\begin{equation}\label{X105}
		\limsup_{r\to\infty}\frac{\log\Bigl(\max\limits_{|z|=r}|f(z)|\Bigr)}{\ms g^{**}(r)}
		\le \limsup_{r\to\infty}\frac{n_f(r)}{\ms g(r)},
	\end{equation}
	which proves \eqref{X239} when $\alpha\in(0,1]$.
	
	Assume now that $\alpha=0$ and suppose that the left-hand side of \eqref{X105}
	is strictly positive.  It follows from \cite[Corollary~6.1]{berg.pedersen:2007} that then
	\begin{equation}\label{X118}
		\limsup_{r\to\infty}\frac{\log\Bigl(\min\limits_{|z|=r}|f(z)|\Bigr)}{\ms g^{**}(r)} > 0.
	\end{equation}
	On the other hand we obtain from \eqref{X238} that
	\[
		\log\Bigl(\min_{|z|=r}|f(z)|\Bigr) \le \log|f(ir)|
		\lesssim \ms g(r) \ll \ms g^{**}(r),
	\]
	which contradicts \eqref{X118}.  Hence the left-hand side of \eqref{X105} vanishes.
\end{proof}

\begin{proof}[Proof of \cref{X45}]
	The function $w_{H,22}(b;\cdot)$ is an entire function with $w_{H,22}(b;0)=1$, 
	is real along the real axis, has only real zeros, and is of bounded type in $\bb C^+$.
	If the right-hand side of \eqref{X88} is infinite, there is nothing to prove.
	So assume that the latter is finite.  Then \cref{X1} implies that
	\[
		\log|w_{H,22}(b;ir)| \asymp \int_a^b K_H(t;r)\DD t \lesssim \ms g(r).
	\]
	When $\alpha\in(0,1]$, we can use \cref{X236} and \cref{X7}\,(i) to obtain
	\begin{equation}\label{X186}
		\limsup_{r\to\infty}\frac{\log\Bigl(\max\limits_{|z|=r}|w_{H,22}(b;z)|\Bigr)}{\ms g^{**}(r)}
		\le \limsup_{r\to\infty}\frac{n_H(r)}{\ms g(r)}
		\le \Delta_+\limsup_{r\to\infty}\frac{1}{\ms g(r)}\int_a^b K_H(t;r)\DD t
	\end{equation}
	with $\Delta_+$ depending only on $c_-,c_+$;
	when $\alpha=0$, the left-hand side of \eqref{X186} vanishes by \cref{X236}.

	To extend the estimates from $|w_{H,22}(b;z)|$ to $\|W_H(b;z)\|$, 
	we use \cite[Proposition~2.3]{baranov.woracek:smsub};
	note that $\ms g^{**}(r)\gg\log r$ by assumption.
\end{proof}

\noindent
We would like to mention that one can always choose a regularly varying function $\ms g$ such that 
\[
	\limsup_{r\to\infty}\frac{1}{\ms g(r)}\int_a^b K_H(t;r)\DD t=1;
\]
see \cite[Theorem~2.3.11]{bingham.goldie.teugels:1989}.

\section{Algorithmic approach}
\label{X134}

We saw in \cref{X132,X133} that the growth of the eigenvalue counting function in general, and of the 
monodromy matrix in the limit circle case, is governed by the integral over the function $K_H(t;r)$.  
The challenge thus is to handle that integral, and this is often difficult.

In this section we consider the limit circle case, and give an algorithm that leads to an evaluation of 
$\int_a^b K_H(t;r)\DD t$ up to a possible logarithmic error.  This algorithm enhances applicability tremendously, 
and has interesting consequences; see, e.g.\ \cref{X19}.  
It is methodologically related to the covering theorem \cite[Theorem~2]{romanov:2017}. 

\subsection[Evaluating $\int_a^b K_H(t;r)\DD t$ by partitioning]{%
Evaluating {\boldmath$\int_a^b K_H(t;r)\DD t$} by partitioning}
\label{X148}

The size of $\int_a^b K_H(t;r)\DD t$ can be determined (up to a logarithmic factor) by using a 
clever partitioning of the interval $[a,b]$, which is constructed algorithmically.

\begin{Definition}\label{X25}
	Let $H\in\bb H_{a,b}$ be in the limit circle case and let $c>0$.  For each $r>0$ we define 
	points $\sigma_j^{(r)}$ and a number $\kappa_H(r)$ by the following procedure:
	\begin{Itemize}
	\item 
		set $\sigma_0^{(r)}\DE a$;
	\item 
		if $\det\Omega\bigl(\sigma_{j-1}^{(r)},b\bigr)>\frac{c}{r^2}$, let $\sigma_j^{(r)}\in(\sigma_{j-1}^{(r)},b)$ 
		be the unique point such that 
		\[
			\det\Omega\bigl(\sigma_{j-1}^{(r)},\sigma_j^{(r)}\bigr) = \frac c{r^2};
		\]
		otherwise, set $\sigma_j^{(r)}\DE b$ and $\kappa_H(r)\DE j$, and terminate.
	\end{Itemize}
\end{Definition}

\noindent
The algorithm terminates after a finite number of steps with 
\begin{equation}\label{X241}
	\kappa_H(r) \le \biggl\lfloor r\cdot\frac{\sqrt{\det\Omega(a,b)}\,}{\sqrt c}\biggr\rfloor+1.
\end{equation}
This follows from Minkowski's determinant inequality,
\begin{equation}\label{X161}
	(\det(A+B))^{\frac12} \ge (\det A)^{\frac12}+(\det B)^{\frac12},
\end{equation}
for positive semi-definite $2\times2$-matrices $A$ and $B$.

\begin{remark}\label{X219}
	Like for the notion of compatible pairs from \cref{X11}, the precise value of $c$ is not important for our purposes. 
	In fact, if $\kappa_{H,c}$ and $\kappa_{H,c'}$ are the functions obtained from the algorithm in \cref{X25} 
	for the constants $c$ and $c'$, respectively, then $\kappa_{H,c'}(r)=\kappa_{H,c} \big(\sqrt{\frac{c}{c'}}r \big)$. 
	As a consequence of \cref{X4}\,(iii), (iv), we then get
	\[
		\frac{1}{\Bigl\lceil\sqrt{\frac{c'}{c}\,}\,\Bigr\rceil}\kappa_{H,c}(r) \le \kappa_{H,c'}(r) 
		\le \biggl\lceil\sqrt{\frac{c}{c'}\,}\,\biggr\rceil\kappa_{H,c}(r).
	\]
	The freedom of choosing $c$ is useful in \cref{X107} and for calculations in concrete examples.
\end{remark}

The following theorem, the main result of this section, gives estimates for the 
integral $\int_a^b K_H(t;r)\DD t$.

\begin{theorem}\label{X3}
	Let $H\in\bb H_{a,b}$ be definite and in the limit circle case.
	Further, let $c>0$, set $r_0 \DE \bigl(\frac{c}{\det\Omega(a,b)}\bigr)^{\frac 12}$ 
	and let $(\hat t,\hat s)$ be the unique compatible pair for $H,r_0$ with constants $c,c$
	\textup{(}see \cref{X39}\textup{)}.
	Then 
	\begin{equation}\label{X157}
		\log 2\cdot\kappa_H(r)-\bigl(\log r+\BigO(1)\bigr) \le \int_a^b K_H(t;r)\DD t
		\le 2e\cdot\kappa_H(r)\bigl(\log r+\BigO(1)\bigr)
	\end{equation}
	for $r>r_0$.
	The expressions $\BigO(1)$ on both sides depend only on $c$ and $\tr\Omega(a,b)$; 
	explicit formulae are \eqref{X99} and \eqref{X159}.
\end{theorem}

\noindent
Before we prove \cref{X3}, let us formulate a corollary.
In order to evaluate the growth compared against a regularly varying comparison function, 
it is sufficient to compute $\kappa_H(r_n)$ along a sufficiently dense sequence $r_n\to\infty$. 
The following corollary illustrates this principle by giving a formula for the order $\rho_H$;
for the definition of $\rho_H$ see the first paragraph of \cref{X138}.

\begin{corollary}\label{X110}
	Let $H\in\bb H_{a,b}$ be in the limit circle case, and let $\rho_H$ be the order of the monodromy matrix $W_H(b;z)$. 
	Further, let $(r_n)_{n=1}^\infty$ be a strictly increasing sequence of positive numbers with $\lim_{n \to \infty} r_n=\infty$ 
	and $\limsup_{n \to \infty}\frac{r_{n+1}}{r_n}<\infty$.
	For fixed $c>0$, let $\kappa_H(r_n)$ be the numbers produced by the algorithm in \cref{X25}.  Then 
	\begin{align*}
		\rho_H = \limsup_{n\to\infty}\frac{\log\kappa_H(r_n)}{\log r_n}.
	\end{align*}
\end{corollary}

\begin{proof}[Proof of \cref{X110}]
	If $(a,b)$ is indivisible, the statement is trivial.  We hence assume that $H$ is definite.
	Let $d \DE \limsup_{n\to\infty}\frac{\log\kappa_H(r_n)}{\log r_n}$ and
	let $r_{n_k}$ be a subsequence of $r_n$ such that $\lim_{k\to\infty}\frac{\log\kappa_H(r_{n_k})}{\log r_{n_k}}=d$.  
	From \cref{X1,X3} we obtain, with some $c_1>0$,
	\begin{align*}
		\rho_H &= \limsup_{r\to\infty}\frac{\log\log\max_{|z|=r}\|W_H(b;z)\|}{\log r}
		\ge \limsup_{r\to\infty}\frac{\log\log|w_{H,22}(b;ir)|}{\log r}
		\\[1ex]
		&= \limsup_{r\to\infty}\frac{\log\bigl[\int_a^b K_H(t;r)\DD t\bigl]}{\log r}
		\ge \limsup_{r\to\infty}\biggl[\frac{\log\kappa_H(r)}{\log r}
		+\frac{\log\bigl(\log 2-c_1\frac{\log r}{\kappa_H(r)}\bigr)}{\log r}\biggr]
		\\[1ex]
		&\ge \limsup_{k\to\infty}\biggl[\frac{\log\kappa_H(r_{n_k})}{\log r_{n_k}}
		+\frac{\log\bigl(\log 2-c_1\frac{\log r_{n_k}}{\kappa_H(r_{n_k})}\bigr)}{\log r_{n_k}}\biggr]
		= d.
	\end{align*}
	The reverse inequality $\rho_H \le d$ is trivially true if $d=1$ since the entries of $W_H(b;\cdot)$
	are of finite exponential type.  We therefore assume that $d<1$. 
	Let $\varepsilon>0$ be such that $d+\varepsilon<1$.  For $r>0$ sufficiently large, let $n(r)$ be the unique integer 
	so that $r \in [r_{n(r)-1},r_{n(r)})$.  Then, by \cref{X45,X3}, we have (with $c_2>0$)
	\begin{align*}
		&\limsup_{r\to\infty}\frac{1}{r^{d+\varepsilon}}\log\Bigl(\max_{|z|=r}\|W_H(b;z)\|\Bigr)
		\le \Delta_{d+\varepsilon}\limsup_{r\to\infty}\frac{1}{r^{d+\varepsilon}}\int_a^b K_H(t;r)\DD t 
		\\[1ex]
		&\le c_2\cdot\limsup_{r\to\infty}\frac{\kappa_H(r)\log r}{r^{d+\varepsilon}} 
		\le c_2\cdot\limsup_{r\to\infty}\frac{\kappa_H(r_{n(r)})\log r_{n(r)}}{r_{n(r)-1}^{d+\varepsilon}} 
		\\[1ex]
		&\le c_2\cdot\limsup_{n \to \infty}\Bigl(\frac{r_n}{r_{n-1}}\Bigr)^{d+\varepsilon} 
		\cdot\limsup_{n\to\infty}\frac{\kappa_H(r_n)\log r_n}{r_n^{d+\varepsilon}}
		= 0.
	\end{align*}
	This shows that $\rho_H \le d+\varepsilon$ and, by the arbitrariness of $\varepsilon$, also $\rho_H \le d$.
\end{proof}

\vspace{6pt}

\noindent
We now come to the proof of \cref{X3}.
The argument that eventually yields the upper bound relies on the monotonicity property 
from \cref{X9}\,(iii) in the form of the next lemma.

\begin{lemma}\label{X68}
	Let $c,r_0$ and $(\hat t,\hat s)$ be as in \cref{X3} and let $r>r_0$.
	Assume that we have points $s_0,s_1,s_2$ such that 
	\begin{equation}\label{X200}
		a \le s_0 < s_1 < s_2 \le b, \qquad 
		\det\Omega(s_0,s_1) \ge \frac{c}{r^2} \ge \det\Omega(s_1,s_2).
	\end{equation}
	Then, for every $\gamma>0$,
	\begin{align}
		\label{X187}
		\int_{s_1}^{s_2}\frac{h_1(t)}{\omega_{1}(\hat s(t;r),t)}\DD t
		&\le e\cdot\log\biggl(\frac{r^2\det\Omega(s_0,s_2)}c\biggr)
		\\[1ex]
		\label{X188}
		&\le \frac{1}{\gamma}\Bigl(\frac{r^2\det\Omega(s_0,s_2)}{c}\Bigr)^\gamma.
	\end{align}
\end{lemma}

\begin{proof}
	The second condition in \eqref{X200} implies that,
	for each $t\in[s_1,s_2]$, we have $\hat s(t;r)\in[s_0,s_1]$.  Let $\gamma>0$. 
	Using the definition of $\hat s$ and \cref{X9}\,(iii) we estimate 
	\begin{align*}
		& \int\limits_{s_1}^{s_2} \frac{h_1(t)}{\omega_{1}(\hat s(t;r),t)}\DD t =\int\limits_{s_1}^{s_2}
		\Bigl(\frac{r^2}c\Bigr)^\gamma\Bigl(\frac{\det\Omega(\hat s(t;r),t)}{\omega_{1}(\hat s(t;r),t)}\Bigr)^\gamma
		\frac{h_1(t)}{\omega_{1}(\hat s(t;r),t)^{1-\gamma}}\DD t
		\\[1ex]
		&\le \Bigl(\frac{r^2}{c}\Bigr)^\gamma\Bigl(\frac{\det\Omega(s_0,s_2)}{\omega_{1}(s_0,s_2)}\Bigr)^\gamma
		\cdot
		\begin{cases}
			\int\limits_{s_1}^{s_2}\frac{h_1(t)}{\omega_{1}(s_1,t)^{1-\gamma}}\DD t \CAS \gamma\in(0,1],
			\\[2ex]
			\int\limits_{s_1}^{s_2}\frac{h_1(t)}{\omega_{1}(s_0,t)^{1-\gamma}}\DD t \CAS \gamma>1
		\end{cases}
		\\[1ex]
		&= \Bigl(\frac{r^2\det\Omega(s_0,s_2)}{c}\Bigr)^\gamma\Bigl(\frac{1}{\omega_{1}(s_0,s_2)}\Bigr)^\gamma
		\cdot\frac{1}{\gamma}\cdot
		\begin{cases}
			\omega_{1}(s_1,s_2)^\gamma \CAS \gamma\in(0,1],
			\\[1ex]
			\omega_{1}(s_0,s_2)^\gamma-\omega_{1}(s_0,s_1)^\gamma \CAS \gamma>1
		\end{cases}
		\\[1ex]
		&\le \frac{1}{\gamma}\Bigl(\frac{r^2\det\Omega(s_0,s_2)}{c}\Bigr)^\gamma,
	\end{align*}
	which proves that the left-hand side in \eqref{X187} is bounded by the expression in \eqref{X188}.
	The latter is minimal when
	\[
		\gamma = \Bigl[\log\Bigl(\frac{r^2\det\Omega(s_0,s_2)}{c}\Bigr)\Bigr]^{-1}.
	\]
	Using this value for $\gamma$ we obtain the inequality in \eqref{X187}.
\end{proof}

\noindent
As a consequence we obtain the following key lemma for the upper bound, 
which is also used in the proofs of \cref{X82,X36}.

\begin{lemma}\label{X189}
	Let $\sigma_j^{(r)}$ and $\kappa_H(r)$ be as in \cref{X25} and $K_H$ and $r_0$ as in \cref{X3}.
	Then, for all $r>r_0$ and $\gamma>0$,
	\begin{align}
		\label{X191}
		\int_{\hat t(r)}^b K_H(t;r)\DD t &\le \sum_{j=2}^{\kappa_H(r)} 
		e\log\biggl(\frac{r^2\det\Omega\bigl(\sigma_{j-2}^{(r)},\sigma_j^{(r)}\bigr)}{c}\biggr)
		\\[1ex]
		\label{X192}
		&\le \sum_{j=2}^{\kappa_H(r)}\frac{1}{\gamma}\biggl(\frac{r^2\det\Omega\bigl(\sigma_{j-2}^{(r)},\sigma_j^{(r)}\bigr)}{c}\biggr)^\gamma.
	\end{align}
\end{lemma}

\begin{proof}
	First note that $\hat t(r)=\sigma_1^{(r)}$.
	For each $j\in\{2,\ldots,\kappa_H(r)\}$ we have $\det\Omega(\sigma_{j-2}^{(r)},\sigma_{j-1}^{(r)})=\frac{c}{r^2}$
	and $\det\Omega(\sigma_{j-1}^{(r)},\sigma_j^{(r)})\le\frac{c}{r^2}$.  Hence we can apply \cref{X68}
	to the triple $\sigma_{j-2}^{(r)},\sigma_{j-1}^{(r)},\sigma_j^{(r)}$ and take the sum over $j$ to obtain
	\begin{align*}
		\int_{\hat t(r)}^b K_H(t;r)\DD t
		&= \sum_{j=2}^{\kappa_H(r)}\int_{\sigma_{j-1}^{(r)}}^{\sigma_j^{(r)}}\frac{h_1(t)}{\omega_1(\hat s(t;r),t)}\DD t
		\le \sum_{j=2}^{\kappa_H(r)} 
		e\log\biggl(\frac{r^2\det\Omega\bigl(\sigma_{j-2}^{(r)},\sigma_j^{(r)}\bigr)}{c}\biggr)
		\\[1ex]
		&\le \sum_{j=2}^{\kappa_H(r)}\frac{1}{\gamma}\biggl(\frac{r^2\det\Omega\bigl(\sigma_{j-2}^{(r)},\sigma_j^{(r)}\bigr)}{c}\biggr)^\gamma.
	\end{align*}
\end{proof}

\noindent
The following lemma is the main ingredient for the lower bound in the proof of \cref{X3}.

\begin{lemma}\label{X67}
	Let $r>r_0$, where $r_0$ is as in \cref{X3},
	and assume that we have points $s_0,\ldots,s_k$ such that
	\begin{align*}	
		a \le s_0 < s_1 < \cdots < s_k \le b
	\end{align*}
	and
	\begin{equation}\label{X61}
		\det\Omega(s_{j-1},s_j) \ge \frac{c}{r^2}, \qquad j\in\{1,\ldots,k\},
	\end{equation}
	Then
	\[
		\int_{s_1}^{s_k} \frac{h_1(t)}{\omega_{1}(\hat s(t;r),t)}\DD t
		\ge k\log 2 - \Bigl[\log r+\log\frac{\tr\Omega(s_0,s_k)}{\sqrt c}\Bigr].
	\]
\end{lemma}

\begin{proof}
	Our assumption \eqref{X61} implies that 
	\[
		\hat s(t;r)\in[s_{j-2},t] \qquad\text{for all } t\in[s_{j-1},s_j],\ j=2,\ldots,k.
	\]
	Using the facts that $\omega_1$ is monotone and additive we obtain 
	\begin{align}
	\nonumber
		& \int\limits_{s_1}^{s_k} \frac{h_1(t)}{\omega_{1}(\hat s(t;r),t)}\DD t
		=\sum_{j=2}^k\,\int\limits_{s_{j-1}}^{s_j}\frac{h_1(t)}{\omega_{1}(\hat s(t;r),t)}\DD t
		\ge \sum_{j=2}^k\,\int\limits_{s_{j-1}}^{s_j}\frac{h_1(t)}{\omega_{1}(s_{j-2},t)}\DD t
		\\
	\nonumber
		&= \sum_{j=2}^k\Bigl(\log\omega_{1}(s_{j-2},s_j)-\log\omega_{1}(s_{j-2},s_{j-1})\Bigr)
		\displaybreak[0]\\
	\nonumber
		&= \sum_{j=2}^k\log\Bigl(1+\frac{\omega_{1}(s_{j-1},s_j)}{\omega_{1}(s_{j-2},s_{j-1})}\Bigr) 
		\ge \sum_{j=2}^k\log\Bigl(2\Bigl[\frac{\omega_{1}(s_{j-1},s_j)}{\omega_{1}(s_{j-2},s_{j-1})}\Bigr]^{\frac 12}\Bigr)
		\\
	\label{X70}
		&= (k-1)\log 2+\frac 12\log\frac{\omega_{1}(s_{k-1},s_k)}{\omega_{1}(s_0,s_1)}.
	\end{align}
	It follows from \eqref{X61} that
	\[
		\omega_{1}(s_{k-1},s_k) \ge \frac{\det\Omega(s_{k-1},s_k)}{\omega_{2}(s_{k-1},s_k)}
		\ge \frac c{r^2}\cdot\frac{1}{\omega_{2}(s_{k-1},s_k)},
	\]
	which yields
	\[
		\frac{\omega_{1}(s_{k-1},s_k)}{\omega_1(s_0,s_1)}
		\ge \frac{c}{r^2}\cdot\frac{1}{\omega_1(s_0,s_k)\omega_2(s_0,s_k)}
		\ge \frac{4c}{r^2(\tr\Omega(s_0,s_k))^2}.
	\]
	We can thus further estimate the expression from \eqref{X70}:
	\begin{align*}
		(k-1)\log 2+\frac 12\log\frac{\omega_{1}(s_{k-1},s_k)}{\omega_{1}(s_0,s_1)}
		&\ge k\log 2-\log 2+\frac 12\log\frac{4c}{r^2(\tr\Omega(s_0,s_k))^2}
		\\[1ex]
		&= k\log 2-\Bigl[\log r+\log\frac{\tr\Omega(s_0,s_k)}{\sqrt c}\Bigr],
	\end{align*}
	which completes the proof.
\end{proof}

\noindent
We now show that, for each fixed Hamiltonian $H$, the integral over $[a,\hat t(r)]$ of $K_H(t;r)$ 
gives only a contribution of logarithmic size.  Thus, it can be neglected if one is not interested 
in uniformity of constants w.r.t.\ $H$.  In the following proposition we return to the more general
situation of \cref{X79}, i.e.\ to general $c_-,c_+$ instead of $c,c$ and to Hamiltonians that need not 
be in the limit circle case.

\begin{proposition}\label{X62}
	Let $H\in\bb H_{a,b}$ be a definite Hamiltonian.
	Further, let $r_0\ge 0$, $c_-,c_+>0$ and assume that $(\hat t,\hat s)$ is a compatible pair 
	for $H,r_0$ with constants $c_-,c_+$.  Then
	\begin{equation}\label{X114}
		\int_a^{\hat t(r)}K_H(t;r)\DD t
		\le 2\log r + \frac{c_-}{c_+} - \log(4c_-) + 2\log\bigl[\tr\Omega\bigl(a,\hat t(r)\bigr)\bigr],
		\qquad r>r_0.
	\end{equation}
	The last term on the right-hand side of \eqref{X114} is non-increasing in $r$.
\end{proposition}

\begin{proof}
	Let $r>r_0$.  Since
	\[
		(\omega_1\omega_2)\bigl(a,\hat t(r)\bigr) \ge \det\Omega\bigl(a,\hat t(r)\bigr)
		\ge \frac{c_-}{r^2},
	\]
	there exists a unique point $\mr t(r)\in(a,\hat t(r)]$ such that 
	\[
		(\omega_1\omega_2)\bigl(a,\mr t(r)\bigr) = \frac{c_-}{r^2}.
	\]
	Let us first estimate the integral over $[a,\mr t(r)]$:
	\begin{align}
		\nonumber
		\int_a^{\mr t(r)}K_H(t;r)\DD t 
		&= \int_a^{\mr t(r)}\frac{\omega_2(a,t)h_1(t)}{\frac{c_+}{r^2}+(\omega_3(a,t))^2}\DD t
		\le \frac{r^2}{c_+}\int_a^{\mr t(r)}\omega_2(a,t)h_1(t)\DD t
		\\[1ex]
		\label{X115}
		&\le \frac{r^2}{c_+}\omega_2\bigl(a,\mr t(r)\bigr)\int_a^{\mr t(r)}h_1(t)\DD t
		= \frac{r^2}{c_+}(\omega_1\omega_2)\bigl(a,\mr t(r)\bigr) = \frac{c_-}{c_+}.
	\end{align}
	For $t\in[\mr t(r),\hat t(r)]$ we have
	\[
		\frac{c_+}{r^2} + \bigl(\omega_3(a,t)\bigr)^2 
		\ge \det\Omega\bigl(a,\hat t(r)\bigr) + \bigl(\omega_3(a,t)\bigr)^2
		\ge \det\Omega(a,t) + \bigl(\omega_3(a,t)\bigr)^2
		= (\omega_1\omega_2)(a,t)
	\]
	and hence
	\begin{align*}
		\int_{\mr t(r)}^{\hat t(r)}K_H(t;r)\DD t
		&\le \int_{\mr t(r)}^{\hat t(r)}\frac{h_1(t)}{\omega_1(a,t)}\DD t
		= \log\frac{\omega_1\bigl(a,\hat t(r)\bigr)}{\omega_1\bigl(a,\mr t(r)\bigr)}
		\le \log\biggl[\frac{\omega_1\bigl(a,\hat t(r)\bigr)}{\omega_1\bigl(a,\mr t(r)\bigr)}
		\cdot\frac{\omega_2\bigl(a,\hat t(r)\bigr)}{\omega_2\bigl(a,\mr t(r)\bigr)}\biggr]
		\\[1ex]
		&= \log\biggl[\frac{r^2}{c_-}\cdot(\omega_1\omega_2)\bigl(a,\hat t(r)\bigr)\biggr]
		= 2\log r - \log c_- + \log\bigl[(\omega_1\omega_2)\bigl(a,\hat t(r)\bigr)\bigr]
		\\[1ex]
		&\le 2\log r - \log c_- + 2\log\biggl[\frac{\tr\Omega\bigl(a,\hat t(r)\bigr)}{2}\biggr],
	\end{align*}
	which, together with \eqref{X115}, yields \eqref{X114}.
\end{proof}

\noindent
We are now ready to prove \cref{X3}.

\begin{proof}[Proof of \cref{X3}]
	Let $\sigma_j^{(r)}$ be the points defined in \cref{X25}.
	Note that $\hat t(r)=\sigma_1^{(r)}$, and thus
	\[
		K_H(t;r) = \frac{h_1(t)}{\omega_{1}(\hat s(t;r),t)}, \qquad 
		t \ge \sigma_1^{(r)}.
	\]
	In order to obtain the desired lower bound, we apply \cref{X67} with the 
	points $\sigma_1^{(r)},\ldots,\sigma_{\kappa_H(r)-1}^{(r)}$ to obtain 
	\begin{align}
		\nonumber
		\int_a^b K_H(t;r)\DD t
		&\ge \int_{\sigma_1^{(r)}}^{\sigma_{\kappa_H(r)-1}^{(r)}}\frac{h_1(t)}{\omega_{1}(\hat s(t;r),t)}\DD t
		\\
		\nonumber
		&\ge (\kappa_H(r)-1)\log 2-\Bigl[\log r+\log\frac{\tr\Omega(a,b)}{\sqrt c}\Bigr]
		\\
		\label{X99}
		&= \kappa_H(r)\log 2-\Bigl[\log r+\log\frac{2\tr\Omega(a,b)}{\sqrt c}\Bigr],
	\end{align}
	which proves the first inequality in \eqref{X157}.
	
	For an upper bound of the integral of $K_H(t;r)$ over the interval $(\sigma_1^{(r)},b)$
	we obtain from \eqref{X191} in \cref{X189} that
	\begin{align}
		\int_{\hat t(r)}^b K_H(t;r)\DD t &\le \sum_{j=2}^{\kappa_H(r)}
		e\log\biggl(\frac{r^2\det\Omega\bigl(\sigma_{j-2}^{(r)},\sigma_j^{(r)}\bigr)}{c}\biggr)
		\le e\bigl(\kappa_H(r)-1\bigr)\log\biggl(\frac{r^2\det\Omega(a,b)}{c}\biggr)
		\nonumber\\[1ex]
		&= 2e\bigl(\kappa_H(r)-1\bigr)\biggl[\log r+\log\frac{\sqrt{\det\Omega(a,b)}\,}{\sqrt{c}}\biggr]
		\nonumber\\[1ex]
		\label{X98}
		&\le 2e\bigl(\kappa_H(r)-1\bigr)\biggl[\log r+\log\frac{\tr\Omega(a,b)}{2\sqrt{c}}\biggr].
	\end{align}
	For the integral of $K_H(t;r)$ over the remaining part $(a,\hat t(r))$ we use \cref{X62},
	which yields
	\begin{equation}\label{X97}
		\int_a^{\hat t(r)}K_H(t;r)\DD t
		\le 2\log r + 1 + 2\log\frac{\tr\Omega(a,b)}{2\sqrt{c}}.
	\end{equation}
	Combining \eqref{X98} and \eqref{X97} we arrive at
	\begin{equation}\label{X158}
		\int_a^b K_H(t;r)\DD t 
		\le 2\Bigl(e\kappa_H(r)-e+1\Bigr)\Bigl[\log r+\log\frac{\tr\Omega(a,b)}{2\sqrt c}\Bigr] + 1.
	\end{equation}
	Since
	\[
		r > r_0 = \frac{\sqrt{c}}{\sqrt{\det\Omega(a,b)}\,}
		\ge \frac{2\sqrt{c}}{\tr\Omega(a,b)},
	\]
	the expression within the square brackets in \eqref{X158} is positive, and hence
	\begin{align}
		\int_a^b K_H(t;r)\DD t &\le 2e\kappa_H(r)\Bigl[\log r+\log\frac{\tr\Omega(a,b)}{2\sqrt c}\Bigr] + 1
		\nonumber\\[1ex]
		\label{X159}
		&\le 2e\kappa_H(r)\Bigl[\log r+\log\frac{\tr\Omega(a,b)}{2\sqrt c}+1\Bigr],
	\end{align}
	which proves the second inequality in \eqref{X157}.
\end{proof}

\noindent
It is a useful observation that one can get estimates for $\kappa_H(r)$ without carrying out the algorithm 
that produces the points $\sigma_j^{(r)}$ precisely. 

\begin{proposition}\label{X4}
	Assume that we have points $s_0,\ldots,s_k$ with $a=s_0<s_1<\cdots<s_{k-1}<s_k=b$. 
	\begin{Enumerate}
	\item 
		If 
		\begin{equation}\label{X166}
			\det\Omega(s_{l-1},s_l) \le \frac{c}{r^2} \qquad \text{for all} \ l\in\{1,\ldots,k\},
		\end{equation}
		then $\kappa_H(r)\le k$.
	\item 
		If 
		\begin{equation}\label{X167}
			\det\Omega(s_{l-1},s_l) \ge \frac{c}{r^2} \qquad \text{for all} \ l\in\{1,\ldots,k-1\},
		\end{equation}
		then $\kappa_H(r)\ge k$.
	\item
		$\kappa_H$ is non-decreasing.
	\item
		$\kappa_H(nr) \le n \kappa_H(r)$ for any integer $n$.
	\end{Enumerate}
\end{proposition}

\begin{proof}
	For the proof of (i) consider the pairwise disjoint intervals $(\sigma_{j-1}^{(r)},\sigma_j^{(r)}]$, 
	$j\in\{1,\ldots,\kappa_H(r)-1\}$.  By \cref{X9}\,(iv), (v) each of these intervals must contain at least one 
	of the points $s_1,\ldots,s_{k-1}$; note that $s_k=b$.
	This shows $\kappa_H(r)\le k$.
	
	The proof of (ii) is analogous: each of the pairwise disjoint intervals $(s_{l-1},s_l]$, $l\in\{1,\ldots,k-1\}$ 
	must contain at least one of the points $\sigma_1^{(r)},\ldots,\sigma_{\kappa_H(r)-1}^{(r)}$. 
	
	To see (iii), let $r\le r'$ and note that $\det\Omega\bigl(\sigma_{j-1}^{(r)},\sigma_j^{(r)}\bigr) = \frac{c}{r^2} \ge \frac{c}{(r')^2}$ 
	for $j \in \{1,\ldots,\kappa_H(r)-1\}$.  Hence $\kappa_H(r') \ge \kappa_H(r)$ by (ii).

	Finally, for item (iv), let $r>0$ and $n\in\bb N$, and set $k\DE 1+\bigl\lfloor\frac{\kappa_H(nr)-1}{n}\bigr\rfloor$. 
	Define $s_l \DE \sigma_{nl}^{(nr)}$ for $l=0,\ldots, k-1$ and $s_k \DE b$.  By \cref{X161},
	\[
		\sqrt{\det\Omega(s_{l-1},s_l)} 
		\ge \sum_{j=1}^n \sqrt{\det\Omega\bigl(\sigma_{n(l-1)+j-1}^{(nr)},\sigma_{n(l-1)+j}^{(nr)}\bigr)} 
		= n \cdot \frac{\sqrt{c}\,}{nr}=\frac{\sqrt c\,}{r}
	\]
	for $l=1,\cdots,k-1$.  Hence item (ii) implies that
	\[
		\kappa_H(r) \ge k = \Bigl\lfloor\frac{\kappa_H(nr)+n-1}{n}\Bigr\rfloor 
		\ge \frac{\kappa_H(nr)}{n}
	\]
	since $\lfloor\cdot\rfloor$ can subtract at most $\frac{n-1}{n}$ from a number in $\frac 1n \bb N$.
\end{proof}

\subsection{Cutting out pieces of the domain}
\label{X211}

It is a very intuitive fact that the growth of the monodromy matrix should not increase when one cuts out pieces of the
Hamiltonian.  A result giving meaning to this intuition is shown in \cite[Theorem~3.4 and Corollary~3.5]{pruckner.woracek:sinqA}.  
There a certain assumption is placed on the piece that is cut out, and it is shown with an operator-theoretic argument 
that the zero counting function does not increase. 
Showing monotonicity of the function $\kappa_H$ in $H$ is straightforward and does not require any additional assumption. 
As a consequence of \cref{X3,X110}, we obtain that $\int_a^b K_H(t;r) \DD t$ might increase at most by a
logarithmic factor when cutting out pieces.  In particular, the order of the monodromy matrix cannot increase.

\begin{theorem}
\label{X19}
	Let $H\in\bb H_{a,b}$ be in the limit circle case and let $\Delta\subseteq[a,b]$ be measurable.  Set
	\begin{align*}
		\lambda(t) &\DE \int_a^t\mathds{1}_{\Delta}(u)\DD u,\quad t\in[a,b], \qquad\qquad
		\tilde a\DE 0,\;\; \tilde b\DE\lambda(b),
		\\[1ex]
		\chi(s) &\DE \min\bigl\{t\in[a,b]\DS \lambda(t)\ge s\bigr\},\quad s\in[\tilde a,\tilde b],
		\\[1ex]
		\wt H(s) &\DE H(\chi(s)),\quad s\in[\tilde a,\tilde b].
	\end{align*}
	Then $\wt H\in\bb H_{\tilde a,\tilde b}$,  
	and, for $\kappa_H$ and $\kappa_{\tilde H}$ as in \cref{X25} with some constant $c>0$,
	we have 
	\[
		\kappa_H(r) \ge \kappa_{\tilde H}(r) \qquad\text{for all }r>0.
	\]
	Further, assume that $\wt H$ is definite, set $r_0\DE\bigl(\frac{c}{\det\Omega_{\tilde H}(\tilde a,\tilde b)}\bigr)^{\frac12}$
	and let $(\hat t,\hat s)$ and $(\tilde t,\tilde s)$ be compatible pairs for $H,r_0$ and $\wt H,r_0$ respectively 
	with constants $c,c$.  Then
	\begin{equation}\label{X173}
		\int_a^b K_{\wt H}(t;r)\DD t\lesssim \log r\cdot\int_a^b K_H(t;r)\DD t,
		\qquad r>r_0.
	\end{equation}
\end{theorem}

\noindent
The intuition behind the proof is simple and is expressed in the picture
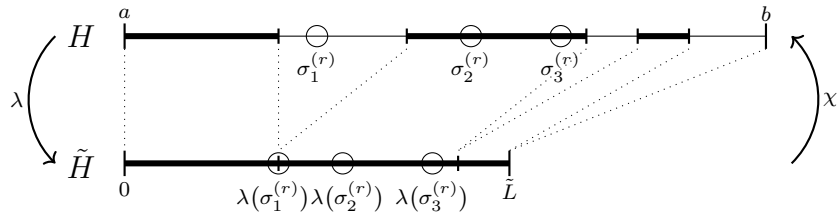
\begin{center}
	\begin{tikzpicture}[x=1.2pt,y=1.2pt,scale=0.8,font=\fontsize{8}{8}]
		\draw (10,30)--(260,30);
		\draw[line width=.8mm] (10,30)--(70,30);
		\draw[line width=.8mm] (120,30)--(190,30);
		\draw[line width=.8mm] (210,30)--(230,30);
		\draw[thick] (10,25)--(10,35);
		\draw[thick] (70,27)--(70,33);
		\draw[thick] (120,27)--(120,33);
		\draw[thick] (190,27)--(190,33);
		\draw[thick] (210,27)--(210,33);
		\draw[thick] (230,27)--(230,33);
		\draw[thick] (260,25)--(260,35);
		\draw (-7,30) node {\large $H$};
		\draw (10,39) node {${\displaystyle a}$};
		\draw (260,39) node {${\displaystyle b}$};
		
		\draw (85,30) circle (5pt);
		\draw (145,30) circle (5pt);
		\draw (180,30) circle (5pt);
		
		\draw (85,18) node {${\sigma_1^{(r)}}$};
		\draw (145,18) node {${\sigma_2^{(r)}}$};
		\draw (180,18) node {${\sigma_3^{(r)}}$};
		
		\draw[dotted] (10,25)--(10,-15);
		\draw[dotted] (70,25)--(70,-15);
		\draw[dotted] (120,25)--(70,-15);
		\draw[dotted] (190,25)--(140,-15);
		\draw[dotted] (210,25)--(140,-15);
		\draw[dotted] (230,25)--(160,-15);
		\draw[dotted] (260,25)--(160,-15);
		
		\draw[line width=.8mm] (10,-20)--(160,-20);
		\draw[thick] (10,-25)--(10,-15);
		\draw[thick] (70,-23)--(70,-17);
		\draw[thick] (140,-23)--(140,-17);
		\draw[thick] (160,-25)--(160,-15);
		\draw (-7,-20) node {\large $\tilde H$};
		\draw (10,-30) node {${\displaystyle 0}$};
		\draw (-32,5) node {$\lambda$};
		\draw (285,5) node {$\chi$};
		
		\draw (70,-20) circle (5pt);
		\draw (95,-20) circle (5pt);
		\draw (130,-20) circle (5pt);
		
		\draw (68,-33) node {${\lambda \big(\sigma_1^{(r)}\big)}$};
		\draw (97,-33) node {${\lambda \big(\sigma_2^{(r)}\big)}$};
		\draw (130,-33) node {${\lambda \big(\sigma_3^{(r)}\big)}$};		
		
		\draw (160,-30) node {${\displaystyle \tilde L}$};
		
		\draw[->, thick] (-17,30) to[out=225, in=135] (-17,-20);
		\draw[->, thick] (270,-20) to[out=45, in=315] (270,30);
	\end{tikzpicture}
\end{center}

\begin{proof}[Proof of \cref{X19}]
	The fact that $\wt H\in\bb H_{\tilde a,\tilde b}$ is shown in the proof of
	\cite[Theorem~3.4]{pruckner.woracek:sinqA}, along with the equality 
	\[
		H(t)\mathds{1}_\Delta(t) = (\wt H\circ\lambda)(t)\mathds{1}_\Delta(t),\qquad t\in[a,b]\ \text{a.e.}
	\]
	Using this we make a change of variable with the absolutely continuous function $\lambda$, and obtain 
	\[
		\int_s^t H(u)\DD u \ge \int_s^t H(u)\mathds{1}_\Delta(u)\DD u
		= \int_s^t(\wt H\circ\lambda)(u)\lambda'(u)\DD u
		= \int_{\lambda(s)}^{\lambda(t)}\wt H(v)\DD v.
	\]
	This implies that 
	\[
		\det\Omega_H(s,t) \ge \det\Omega_{\wt H}\bigl(\lambda(s),\lambda(t)\bigr),\qquad 
		a \le s \le t \le b.
	\]
	Hence, each of the pairwise disjoint intervals (quantities $\tilde\sigma_j^{(r)}$ refer to $\wt H$)
	\[
		[\tilde\sigma_{j-1}^{(r)},\tilde\sigma_j^{(r)}),\qquad 
		j\in\{1,\ldots,\kappa_{\tilde H}(r)-1\}
	\]
	contains at least one of the point $\lambda(\sigma_j^{(r)})$, and 
	it follows that $\kappa_H(r)\ge\kappa_{\wt H}(r)$.
	
	For the last assertion assume that $\wt H$ is definite.  Then clearly also $H$ is definite.
	The inequality in \eqref{X173} follows then from \cref{X3}.
\end{proof}

\subsection{Recovering finite-type results with universal constants}
\label{X149}

We revisit some upper bounds for order and type known from the literature and give an 
alternative proof by using the algorithm developed in \cref{X148}.  This leads to estimates that 
do not depend on the Hamiltonian in any other way than on the constants used in the algorithm. 
The important feature is that these bounds hold pointwise with universal constants, and not just
asymptotically as $r\to\infty$.  This feature is, in particular, exploited 
in \cref{X210}.

The logic of the arguments is appealingly simple. 
To explain it, let us analyse the proof of the upper bound in \cref{X3}, based on \cref{X189}. 
One main task was to estimate the sum in \eqref{X191}.  This can be done in different ways, 
and one of them led to the bound from \cref{X3}; namely, we used the crude 
estimate $\det\Omega(s_{j-2}^{(r)},\sigma_j^{(r)})\le\det\Omega(a,b)$ for each $j\in\{2,\ldots,\kappa_H(r)\}$.
Now we proceed differently: we use the weaker estimate \eqref{X192} instead of \eqref{X191}, but use 
a finer estimate for $\det\Omega(\cdot,\cdot)$ such that summation is possible. 

As a first illustration of this method, we discuss the classical Krein--de~Branges formula, cf.\ \cite[(3.4)]{krein:1951} and
\cite[Theorem~X]{debranges:1961}.
In the language of growth of the monodromy matrix it says that the exponential type of $w_{H,22}(b;\cdot)$ 
equals $\int_a^b\sqrt{\det H(t)}\DD t$.
The following statement is asymptotically weaker, but it is pointwise and with explicit constants.
To obtain the connection with the growth of $w_{H,22}(b;\cdot)$ one has to combine the next proposition
with \cref{X1,X62}.

\begin{proposition}\label{X82}
	Let $H\in\bb H_{a,b}$ be definite and in the limit circle case.  
	Further, let $c>0$, set $r_0 \DE \bigl(\frac{c}{\det\Omega(a,b)}\bigr)^{\frac 12}$ 
	and let $(\hat t,\hat s)$ be the unique compatible pair for $H,r_0$ with constants $c,c$.
	Then
	\begin{equation}\label{X169}
	\begin{aligned}
		\frac{\int_a^b\sqrt{\det H(t)}\DD t}{\sqrt{c}}\log 2\cdot r-\Bigl(\log r+\log\frac{2\tr\Omega(a,b)}{\sqrt{c}}\Bigr)
		&\le \int_{\hat t(r)}^b K_H(t;r)\DD t 
		\\[1ex]
		&\le \frac{4\sqrt{\det\Omega(a,b)}\,}{\sqrt{c}\,}\cdot  r
	\end{aligned}
	\end{equation}
	for $r>r_0$.
\end{proposition}

\begin{proof}
	For the upper bound we use \eqref{X192} with $\gamma=\frac12$ and Minkowski's determinant inequality, \eqref{X161}, to obtain
	\begin{align*}
		\int_{\hat t(r)}^b K_H(t;r)\DD t
		&\le \frac{2}{\sqrt c\,}\cdot r\cdot\sum_{j=2}^{\kappa_H(r)}\sqrt{\det\Omega\bigl(\sigma_{j-2}^{(r)},\sigma_j^{(r)}\bigr)}
		\\[1ex]
		&\le \frac{2}{\sqrt c\,}\cdot r\cdot\biggl[
		\det\biggl(\,\sum_{j=2}^{\kappa_H(r)}\Omega\bigl(\sigma_{j-2}^{(r)},\sigma_j^{(r)}\bigr)\biggr)\biggr]^{\frac12}
		\\[1ex]
		&\le \frac{2}{\sqrt c\,}\cdot r \cdot\sqrt{\det\bigl(2\Omega(a,b)\bigr)}
		= \frac{4}{\sqrt c\,}\cdot r\cdot\sqrt{\det\Omega(a,b)},
	\end{align*}
	where we also used that all involved matrices are non-negative.

	Let us now come to the proof of the first inequality in \eqref{X169}. 
	If $\det H=0$ a.e.\ there is nothing to prove. 
	Otherwise, choose $s_j\in[a,b]$ such that 
	\begin{equation}\label{X242}
		\int_a^{s_j}\sqrt{\det H(t)}\DD t = j\cdot\frac{\sqrt c\,}{r},
		\qquad j=0,\ldots,\bigg\lfloor\frac{r}{\sqrt c\,}\int_a^b\sqrt{\det H(t)}\DD t\bigg\rfloor\ED k.
	\end{equation}
	With several applications of the Cauchy--Schwarz inequality we obtain
	\begin{align}
		\det\Omega(s_{j-1},s_j) 
		&= \biggl[\biggl(\int_{s_{j-1}}^{s_j}h_1(t)\DD t\biggr)^{\frac12}\biggl(\,\int_{s_{j-1}}^{s_j}h_2(t)\DD t\biggr)^{\frac12}
		- \int_{s_{j-1}}^{s_j}h_3(t)\DD t\biggr]
		\nonumber\\[1ex]
		&\quad \times \biggl[\biggl(\int_{s_{j-1}}^{s_j}h_1(t)\DD t\biggr)^{\frac12}\biggl(\,\int_{s_{j-1}}^{s_j}h_2(t)\DD t\biggr)^{\frac12}
		+ \int_{s_{j-1}}^{s_j}h_3(t)\DD t\biggr]
		\nonumber\\[1ex]
		&\ge \int_{s_{j-1}}^{s_j}\Bigl(\sqrt{h_1(t)h_2(t)}-h_3(t)\Bigr)\DD t
		\cdot\int_{s_{j-1}}^{s_j}\Bigl(\sqrt{h_1(t)h_2(t)}+h_3(t)\Bigr)\DD t
		\nonumber\displaybreak[0]\\[1ex]
		&\ge \biggl[\;\int_{s_{j-1}}^{s_j}\Bigl(\sqrt{h_1(t)h_2(t)}-h_3(t)\Bigr)^{\frac12}
		\Bigl(\sqrt{h_1(t)h_2(t)}+h_3(t)\Bigr)^{\frac12}\DD t\biggr]^2
		\nonumber\\[1ex]
		&= \biggl[\int_{s_{j-1}}^{s_j}\sqrt{\det H(t)}\DD t\biggr]^2
		= \frac{c}{r^2}
		\label{X243}
	\end{align}
	for $j\in\{1,\ldots,k\}$.  In particular, $\det\Omega(a,s_1)\ge\frac{c}{r^2}$, which implies
	that $s_1\ge\hat t(r)$.
	Hence, by \cref{X67}, we obtain
	\begin{align*}
		\int_{\hat t(r)}^b K_H(t;r)\DD t &\ge \int_{s_1}^{s_k} K_H(t;r)\DD t
		\ge k\log 2-\biggl(\log r+\log\frac{\tr\Omega(a,b)}{\sqrt{c}}\biggr)
		\\[1ex]
		&\ge \biggl(\frac{r}{\sqrt c\,}\int_a^b\sqrt{\det H(t)}\DD t-1\biggr)\log 2
		-\biggl(\log r+\log\frac{\tr\Omega(a,b)}{\sqrt{c}}\biggr),
	\end{align*}
	which finishes the proof.
\end{proof}

\begin{Remark}\label{X244}
	If $\det H$ does not vanish identically, then $\kappa_H(r)\asymp r$;
	this follows from \eqref{X241}, \eqref{X242}, \eqref{X243} and \cref{X4}\,(ii).
\end{Remark}

\noindent
In the remainder of this subsection and in some later subsections we consider a 
definite Hamiltonian $H$ with zero determinant.  By reparameterisation, we can assume, 
without loss of generality, that $\tr H(t)=1$ a.e.\ and hence
\begin{equation}\label{X37}
	H(t) = \xi_{\varphi(t)}\xi_{\varphi(t)}^T
	= \begin{pmatrix} 
		\cos^2\varphi(t) & \cos\varphi(t)\cdot\sin\varphi(t) \\[1ex]
		\cos\varphi(t)\cdot\sin\varphi(t) & \sin^2\varphi(t) 
	\end{pmatrix}
\end{equation}
where $\xi_\phi$ is defined in \eqref{X104}
and $\varphi\DF[a,b]\to\bb R$ is a non-constant, measurable function with $a,b\in\bb R$, $a<b$.
In \cref{X36} we prove upper bounds for $\int_{\hat t(r)}^b K_H(t;r)\DD t$ when $\varphi$ is bounded, 
H\"older continuous, or of bounded variation, respectively.
We emphasise that our bounds hold pointwise, and one should not think of them asymptotically.
An asymptotic estimate for H\"older continuous functions, as considered in item (ii) below, 
is proved in \cite[Corollary~5.3]{pruckner.woracek:sinqA};
an estimate for the order in the situation of item (iii) is proved in \cite[Corollary~4]{romanov:2017}.
For arbitrary measurable functions $\varphi$ one obtains an upper bound also from \cref{X82};
the difference to \cref{X36}\,(i) is that, in the latter, the constant in the estimate improves 
when $\varphi$ does not vary too much.
To quantify this, let us recall the \emph{oscillation} of a bounded function $\varphi$:
\[
	\osc_\varphi(s,t) \DE \esssup(\varphi|_{[s,t]})-\essinf(\varphi|_{[s,t]}).
\]

\begin{proposition}\label{X36}
	Let $\varphi\DF[a,b]\to\bb R$ be measurable and non-constant with $a,b\in\bb R$, $a<b$, and consider the Hamiltonian 
	$H(t)\DE\xi_{\varphi(t)}\xi_{\varphi(t)}^T$.  Further, let $c>0$, set $r_0 \DE \bigl(\frac{c}{\det\Omega(a,b)}\bigr)^{\frac 12}$ 
	and let $(\hat t,\hat s)$ be the unique compatible pair for $H,r_0$ with constants $c,c$.
	Then, for $r>r_0$, the following statements hold.
	\begin{Enumerate}
	\item 
		If $\varphi\in L^\infty(a,b)$, then
		\[
			\int_{\hat t(r)}^b K_H(t;r)\DD t
			\le \frac{2}{\sqrt c\,}(b-a)\osc_\varphi(a,b)\cdot r.
		\]
	\item 
	If $\varphi$ is H\"older continuous on $[a,b]$ with exponent $\nu$ and constant $\eta$, then
		\[
			\int_{\hat t(r)}^b K_H(t;r)\DD t
			\le 4(1+\nu)(b-a)\Bigl(\frac{\eta}{2\sqrt c\,}\Bigr)^{\frac{1}{1+\nu}}\cdot r^{\frac{1}{1+\nu}}.
		\]
	\item 
	If $\varphi$ is of bounded variation on $[a,b]$, and $V_a^b(\varphi)$ denotes its total variation, then
		\[
			\int_{\hat t(r)}^b K_H(t;r)\DD t
			\le \frac{4\sqrt{2}\,}{\sqrt[4]{c}\,} \sqrt{(b-a)V_a^b(\varphi)}\cdot r^{\frac 12}.
		\]
	\end{Enumerate}
\end{proposition}

\begin{remark}\label{X78}
\rule{0ex}{1ex}
\begin{Enumerate}
\item
	Together with \cref{X62} and \cref{X7}\,(i) one immediately obtains upper bounds for the 
	eigenvalue counting function $n_H(r)$.
\item
	In \cref{X127} we provide, for any given $\nu\in(0,1)$, an example of a H\"older continuous function with exponent $\nu$ 
	where $\int_a^b K_H(t;r)\DD t\asymp r^{\frac{1}{1+\nu}}$.
\item
	\Cref{X36}\,(ii) implies, in particular, that, when $\varphi$ is a Lipschitz function, then
	$\int_{\hat t(r)}^b K_H(t;r)\DD t\lesssim r^{\frac12}$.
	If, in addition, $\varphi$ is not constant, then actually, $\int_a^b K_H(t;r)\DD t\asymp r^{\frac12}$
	as is shown in \cref{X190}.
\end{Enumerate}
\end{remark}

\medskip

\noindent
Before we prove \cref{X36}, let us consider a lemma that provides an upper bound for $\det\Omega(s,t)$ 
in terms of the oscillation of $\varphi$.

\begin{lemma}
\label{X33}
	Let $\varphi\in L^\infty(a,b)$.  Then, for all $a\le s<t\le b$, 
	\[
		\det\Omega(s,t) \le \Bigl[\frac{t-s}2\cdot\osc_\varphi(s,t)\Bigr]^2.
	\]
\end{lemma}

\begin{proof}
	Set 
	\[
		\kappa\DE\frac 12\bigl[\esssup(\varphi|_{[s,t]})+\essinf(\varphi|_{[s,t]})\bigr],
	\]
	so that 
	\[
		\|\varphi|_{[s,t]}-\kappa\|_\infty \le \frac 12 \osc_\varphi(s,t).
	\]
	Since $\xi_{\varphi(x)-\kappa}\xi_{\varphi(x)-\kappa}^T	= N_\kappa\xi_{\varphi(x)}\xi_{\varphi(x)}^TN_\kappa^T$,
	where $N_\kappa$ is as in \eqref{X203}, we obtain 
	\begin{align*}
		\det\Omega(s,t) &= \det\int_s^t \xi_{\varphi(u)}\xi_{\varphi(u)}^T\DD u
		= \det\int_s^t \xi_{\varphi(u)-\kappa}\xi_{\varphi(u)-\kappa}^T\DD u
		\\
		&\le \int_s^t \cos^2\bigl(\varphi(u)-\kappa\bigr)\DD u
		\cdot\int_s^t \sin^2\bigl(\varphi(u)-\kappa\bigr)\DD u
		\\[1ex]
		&\le (t-s)\cdot(t-s)\Bigl(\frac{\osc_\varphi(s,t)}{2}\Bigr)^2.
	\end{align*}
\end{proof}

\begin{proof}[Proof of \cref{X36}]
	Let $r>r_0$.
	The statement in (i) follows directly from \cref{X82,X33}:
	\[
		\int_{\hat t(r)}^b K_H(t;r)\DD t \le \frac{4\sqrt{\det\Omega(a,b)}\,}{\sqrt{c}}\cdot r
		\le \frac{4}{\sqrt{c}\,}\cdot\frac{b-a}{2}\osc_\varphi(a,b)\cdot r.
	\]

	For the remaining cases, (ii) and (iii), let us shorten notation and write $\sigma_j\DE\sigma_j^{(r)}$.
	Note that also here we have $\varphi\in L^\infty(a,b)$.
	Hence we can apply \cref{X33} and \eqref{X192} with arbitrary $\gamma>0$ to obtain
	\begin{align}
		\int_{\hat t(r)}^b K_H(t;r)\DD t &\le \sum_{j=2}^{\kappa_H(r)}\frac{1}{\gamma}
		\biggl(\frac{r^2\det\Omega(\sigma_{j-2},\sigma_j)}{c}\biggr)^\gamma
		\nonumber\\[1ex]
		&\le \sum_{j=2}^{\kappa_H(r)}\frac{1}{\gamma} 
		\biggl[\frac{r^2}{c}\cdot\biggl(\frac{(\sigma_j-\sigma_{j-2})\osc_\varphi(\sigma_{j-2},\sigma_j)}{2}\biggr)^2\,\biggr]^\gamma
		\nonumber\\[1ex]
		\label{X193}
		&= \sum_{j=2}^{\kappa_H(r)}\frac{1}{\gamma} 
		\biggl[\frac{r(\sigma_j-\sigma_{j-2})\osc_\varphi(\sigma_{j-2},\sigma_j)}{2\sqrt{c}}\biggr]^{2\gamma}.
	\end{align}
	We now choose $\gamma$ according to the different cases.
	
	(ii) 
	The H\"older condition implies that $\osc_\varphi(s,t)\le\eta(t-s)^\nu$.
	Hence, the estimate \eqref{X193} with $\gamma=\frac{1}{2(1+\nu)}$ yields
	\begin{align*}
		\int_{\hat t(r)}^b K_H(t;r)\DD t
		&\le 2(1+\nu)r^{\frac{1}{1+\nu}}\sum_{j=2}^{\kappa_H(r)}
		\biggl[\frac{\eta(\sigma_j-\sigma_{j-2})^{1+\nu}}{2\sqrt{c}}\biggr]^{\frac{1}{1+\nu}}
		\\[1ex]
		&= 2(1+\nu)\Bigl(\frac{\eta}{2\sqrt{c}\,}\Bigr)^{\frac{1}{1+\nu}}r^{\frac{1}{1+\nu}}
		\sum_{j=2}^{\kappa_H(r)}(\sigma_j-\sigma_{j-2})
		\\[1ex]
		&\le 4(1+\nu)(b-a)\Bigl(\frac{\eta}{2\sqrt{c}\,}\Bigr)^{\frac{1}{1+\nu}}r^{\frac{1}{1+\nu}}.
	\end{align*}
	
	(iii) 
	In the final case we use $\gamma=\frac14$ in \eqref{X193}, the estimate $\osc_\varphi(s,t)\le V_s^t(\varphi)$
	and the Cauchy--Schwarz inequality to arrive at
	\begin{align*}
		\int_{\hat t(r)}^b K_H(t;r)\DD t &\le 4\frac{1}{(2\sqrt{c})^{\frac12}}\cdot r^{\frac12}
		\sum_{j=2}^{\kappa_H(r)}\bigl[(\sigma_j-\sigma_{j-2})V_{\sigma_{j-2}}^{\sigma_j}(\varphi)\bigr]^{\frac12}
		\\[1ex]
		&\le \frac{2\sqrt{2}}{\sqrt[4]{c}}\cdot r^{\frac12}
		\Biggl(\sum_{j=2}^{\kappa_H(r)}(\sigma_j-\sigma_{j-2})\Biggr)^{\frac12}
		\Biggl(\sum_{j=2}^{\kappa_H(r)}V_{\sigma_{j-2}}^{\sigma_j}(\varphi)\Biggr)^{\frac12}
		\\[1ex]
		&\le \frac{2\sqrt{2}}{\sqrt[4]{c}}\cdot r^{\frac12}\sqrt{2(b-a)}\sqrt{2V_a^b(\varphi)},
	\end{align*}
	which implies the claimed inequality.
\end{proof}

\section{Hamiltonians that oscillate at one endpoint}
\label{X210}

In this section we consider Hamiltonians of the form $H(t)=\xi_{\varphi(t)}\xi_{\varphi(t)}^T$
(see \eqref{X37}) where $\xi_\phi=(\cos\phi,\sin\phi)^T$ and $\varphi$ is close
to a piecewise linear function, which may grow or oscillate (with regularly varying behaviour of the slopes)
towards the left endpoint $a$.
In \cref{X214} we provide bounds for $\kappa_H(r)$ and $\int_a^b K_H(t;r)\DD t$.
In many cases the actual asymptotic behaviour of the latter integral
is determined by $\kappa_H(r)$ without a logarithmic factor unlike in \cref{X3}.

In \cref{X151} we study a couple of examples.  In particular, we apply \cref{X214} to a Hamiltonian
where $\varphi$ is a chirp signal to illustrate the applicability of \cref{X214}; see \cref{X42}.
Further, in \cref{X212} we use \cref{X214} to obtain an inverse result, namely, we construct
an explicit Hamiltonian whose monodromy matrix has prescribed growth; see \cref{X213}.

\subsection{Determining the growth}
\label{X107}

In the following theorem, the main theorem in this section, we prove estimates 
for $\kappa_H(r)$ and $\int_a^b K_H(t;r)\DD t$ using an infinite partition and estimates 
as one approaches $a$.
For the formulation of the theorem recall the notion of an asymptotic inverse from \cref{X100}.

\begin{theorem}\label{X214}
	Let $a,b\in\bb R$ with $a<b$, let $\varphi:(a,b] \to \bb R$ be locally Lipschitz continuous,
	and set $H(t)\DE\xi_{\varphi(t)}\xi_{\varphi(t)}^T$ where $\xi_\phi=(\cos\phi,\sin\phi)^T$.
	Further, let $\kappa_H$ and $K_H(t;r)$ be as in \cref{X3}.

	Let $\ms f_1,\ms f_2:[1,\infty)\to(0,\infty)$ be regularly varying functions, 
	set $\rho_i\DE\Ind\ms f_i$, $i=1,2$, and assume that $\int_1^\infty \ms f_1(x)\DD x<\infty$
	\textup{(}note that hence $\rho_1\le-1$\textup{)}.
	Further, suppose that we have a strictly decreasing sequence $(x_n)_{n=1}^\infty$ in $(a,b]$ 
	with $x_1=b$ and $\lim_{n\to\infty}x_n=a$ such that
	\begin{equation}\label{X201}
		x_{n-1}-x_n \asymp \ms f_1(n), \qquad 
		|\varphi(x_{n-1})-\varphi(x_n)| \asymp \ms f_2(n), \qquad
		\big\|\varphi'|_{[x_n,x_{n-1}]}\big\|_\infty \asymp \frac{\ms f_2(n)}{\ms f_1(n)}
	\end{equation}
	for $n\in\bb N$.

	Consider the following three cases.
	\begin{Enumerate}
	\item[\textup{(A)}]
		Assume that $\ms f_2(x)\gtrsim1$ and let $\ms g_A$ be an asymptotic inverse
		of $\frac{\ms f_2}{\ms f_1}$ 
		\textup{(}note that $\rho_2\ge0$ and hence $\Ind\frac{\ms f_2}{\ms f_1}=\rho_2-\rho_1\ge1$\textup{)}.
		\begin{Itemize}
		\item
			If $\rho_1<-1$, then
			\begin{equation}\label{X260}
				\kappa_H(r) \asymp \int_a^b K_H(t;r)\DD t 
				\asymp r^{\frac12}\int_1^{\ms g_A(r)}\sqrt{(\ms f_1\ms f_2)(t)}\DD t
				\ED \ms b_{A,1}(r),
			\end{equation}
			where $\ms b_{A,1}$ is regularly varying and satisfies
			\begin{equation}\label{X275}
				\ms b_{A,1}(r) \;
				\left\{\begin{array}{ll} 
					\asymp & \text{if} \ \frac{\rho_1+\rho_2}{2}>-1  \\[1ex]
					\gg & \text{if} \ \frac{\rho_1+\rho_2}{2}\le-1
				\end{array}\right\} \;
				\ms g_A(r)\ms f_2\bigl(\ms g_A(r)\bigr)
				\asymp r\ms g_A(r)\ms f_1\bigl(\ms g_A(r)\bigr)
			\end{equation}
			and
			\[
				\Ind\ms b_{A,1} =
				\begin{cases}
					\frac{1+\rho_2}{\rho_2-\rho_1} \CAS \frac{\rho_1+\rho_2}{2}>-1,
					\\[1ex]
					\frac12 \CAS \frac{\rho_1+\rho_2}{2}\le-1.
				\end{cases}
			\]
		\item
			If $\rho_1=-1$, then
			\begin{equation}\label{X261}
				r\ms g_A(r)\ms f_1\bigl(\ms g_A(r)\bigr) \lesssim \kappa_H(r)
				\lesssim \int_a^b K_H(t;r)\DD t
				\lesssim r\int_{\ms g_A(r)}^\infty \ms f_1(t)\DD t;
			\end{equation}
			both bounds are regularly varying with index $1$.
		\end{Itemize}
	\item[\textup{(B)}]
		Assume that $\ms f_2(x)\lesssim1$ and let $\ms g_B$ be an asymptotic inverse
		of $\frac{1}{\ms f_1\ms f_2}$ \textup{(}note that $\Ind\frac{1}{\ms f_1\ms f_2}=-\rho_1-\rho_2\ge1$\textup{)}.
		Then
		\begin{equation}\label{X264}
			\ms b_{B,1}(r) \lesssim \kappa_H(r) \lesssim \int_a^b K_H(t;r)\DD t
			\lesssim \max\bigl\{\ms b_{B,1}(r),\ms b_{B,2}(r)\bigr\},
		\end{equation}
		where
		\[
			\ms b_{B,1}(r) \DE r^{\frac12}\int_1^{\ms g_B(r)}\sqrt{(\ms f_1\ms f_2)(t)}\DD t,
			\qquad
			\ms b_{B,2}(r) \DE r\int_{\ms g_B(r)}^\infty \ms f_1(t)\DD t.
		\]
		The functions $\ms b_{B,i}$ are regularly varying and satisfy
		\begin{align}
			\ms b_{B,1}(r) &
			\begin{cases}
				\;\asymp r^{\frac12} \gg \ms g_B(r) \CAS \frac{\rho_1+\rho_2}{2}<-1,
				\\[1ex]
				\;\gg \ms g_B(r) \CAS \frac{\rho_1+\rho_2}{2}=-1,
				\\[1ex]
				\;\asymp \ms g_B(r) \CAS \frac{\rho_1+\rho_2}{2}>-1,
			\end{cases}
			\label{X265}
			\\[1ex]
			\Ind\ms b_{B,1} &=
			\begin{cases}
				\frac12 \CAS \frac{\rho_1+\rho_2}{2}\le-1,
				\\[1ex]
				-\frac{1}{\rho_1+\rho_2} \CAS \frac{\rho_1+\rho_2}{2}>-1,
			\end{cases}
			\label{X266}
			\\[1ex]
			\ms b_{B,2}(r) &
			\left\{\begin{array}{ll} 
				\asymp & \text{if} \ \rho_1<-1  \\[1ex]
				\gg & \text{if} \ \rho_1=-1
			\end{array}\right\} \;
			\frac{\ms g_B(r)}{\ms f_2(\ms g_B(r))},
			\label{X267}
			\\[1ex]
			\Ind\ms b_{B,2} &= \frac{\rho_2-1}{\rho_1+\rho_2}.
			\label{X268}
		\end{align}
		\begin{Itemize}
		\item
			Assume that $\int_1^\infty\sqrt{(\ms f_1\ms f_2)(t)}\DD t<\infty$.
			Then $\ms b_{B,1}(r)\asymp r^{\frac12}$ and
			\begin{equation}\label{X270}
				\ms b_{B,2}(r) \; \left\{\begin{array}{c} \lesssim \\ \gg \end{array}\right\} \; r^{\frac12}
				\quad\Leftrightarrow\quad
				x^2\ms f_1(x) \; \left\{\begin{array}{c} \lesssim \\ \gg \end{array}\right\} \; \ms f_2(x).
			\end{equation}
		\item
			Assume that $\int_1^\infty\sqrt{(\ms f_1\ms f_2)(t)}\DD t=\infty$.
			If $\rho_2<0$ or $\rho_1=-1$, then $\ms b_{B,2}(r)\gg\ms b_{B,1}(r)$.
			If $\rho_1\in(-2,-1)$, then $\ms b_{B,2}(r)\gtrsim\ms b_{B,1}(r)$;
			if, additionally, $\ms f_2(x)\ll 1$, then $\ms b_{B,2}(r)\gg\ms b_{B,1}(r)$.
		\end{Itemize}
	\item[\textup{(B$^+$)}]
		Assume that $\ms f_2\lesssim1$ and $\varphi\in L^\infty(a,b)$,
		and let $\ms g_B$, $\ms b_{B,1}$ and $\ms b_{B,2}$ be as in case \textup{(B)}.
		Further, let $\ms f_3:[1,\infty)\to(0,\infty)$ be a regularly varying function
		such that
		\begin{equation}\label{X271}
			\osc_\varphi(a,x_n) \lesssim \ms f_3(n), \qquad n\in\bb N.
		\end{equation}
		Then
		\begin{equation}\label{X272}
			\ms b_{B,1}(r) \lesssim \kappa_H(r) \lesssim \int_a^b K_H(t;r)\DD t
			\lesssim \max\bigl\{\ms b_{B,1}(r),\ms b_{B,2}(r)\ms f_3\bigl(\ms g_B(r)\bigr)\bigr\}.
		\end{equation}
		In particular, if $\rho_1<-1$ and 
		either $\ms f_2(x)=\ms f_3(x)$ or $\ms f_3(x)=\frac{1}{x}\sqrt{\frac{\ms f_2(x)}{\ms f_1(x)}}$, then
		\begin{equation}\label{X273}
			\kappa_H(r) \asymp \int_a^b K_H(t;r)\DD t \asymp \ms b_{B,1}(r).
		\end{equation}
	\end{Enumerate}
	All asymptotic relations in $r$ are valid for $r>r_0'$ with some $r_0'\ge r_0$.
\end{theorem}

\medskip

\noindent
Before we prove \cref{X214} we provide some comments.
Examples that illustrate the application of \cref{X214} are given in \cref{X151}.

\begin{Remark}\label{X240}
\rule{0ex}{1ex}
\begin{Enumerate}
\item
	Instead of \eqref{X201} it is sufficient to assume that $\lesssim$ in the first and third relations
	and $\gtrsim$ in the second relations hold since then
	\[
		\frac{\ms f_2(n)}{\ms f_1(n)} \lesssim \frac{|\varphi(x_{n-1})-\varphi(x_n)|}{x_{n-1}-x_n}
		= \frac{1}{x_{n-1}-x_n}\bigg|\int_{x_n}^{x_{n-1}}\varphi'(t)\DD t\bigg|
		\le \big\|\varphi'|_{[x_n,x_{n-1}]}\big\|_\infty
		\lesssim \frac{\ms f_2(n)}{\ms f_1(n)}.
	\]
\item
	Assume that $\varphi'$ is continuous and $|\varphi'|$ is non-increasing on $(a,a+\varepsilon)$ with some $\varepsilon>0$.
	Then the first two relations in \eqref{X201} imply the third one.
	This can be seen as follows.  It follows from the assumptions that $\varphi'$ has constant sign on $(a,a+\varepsilon')$
	with some $\varepsilon'\le\varepsilon$.
	Let $n\in\bb N$ such that $x_{n-1}\in(a,a+\varepsilon')$ and let $t\in[x_n,x_{n-1}]$.
	Then $|\varphi'(t)|\le|\varphi'(s)|$ for all $s\in[x_{n+1},x_n]$ and hence
	\begin{align*}
		|\varphi'(t)| &= \frac{1}{x_n-x_{n+1}}\int_{x_{n+1}}^{x_n}|\varphi'(t)|\DD s
		\le \frac{1}{x_n-x_{n+1}}\int_{x_{n+1}}^{x_n}|\varphi'(s)|\DD s
		\\[1ex]
		&= \bigg|\frac{1}{x_n-x_{n+1}}\int_{x_{n+1}}^{x_n}\varphi'(s)\DD s\bigg|
		= \frac{|\varphi(x_n)-\varphi(x_{n+1})|}{x_n-x_{n+1}}
		\asymp \frac{\ms f_2(n)}{\ms f_1(n)}.
	\end{align*}
	We can now apply item (i) of this remark.
\end{Enumerate}
\end{Remark}

\medskip

\noindent
In the following let $H$ always be as in \eqref{X37}.
For the proof of \cref{X214} we start with a core lemma, which confirms the intuition 
that $\det\Omega$ is a measure for the speed of change of $\varphi(t)$.

\begin{lemma}\label{X54}
	Let $\varphi:[a,b]\to\bb R$ be measurable, and $a\le s<t\le b$.  Then
	\begin{equation}\label{X106}
		\det\Omega(s,t) = \frac12 \int_s^t \int_s^t \sin^2\bigl(\varphi(x)-\varphi(y)\bigr)\DD x \DD y.
	\end{equation}
\end{lemma}

\begin{proof}
	We have
	\begin{align*}
		\det\Omega(s,t)
		&= \int_s^t \int_s^t \bigl[\sin^2\varphi(x)\cos^2\varphi(y)
		-\sin\varphi(x)\cos\varphi(x)\sin\varphi(y)\cos\varphi(y)\big]\DD x \DD y
		\\[1ex]
		&= \int_s^t \int_s^t \sin\varphi(x)\cos\varphi(y)\sin\bigl(\varphi(x)-\varphi(y)\bigr)\DD x \DD y.
	\end{align*}
	Swapping $x$ and $y$ we obtain
	\[
		\det\Omega(s,t)
		= -\int_s^t \int_s^t \sin\varphi(y)\cos\varphi(x)\sin\bigl(\varphi(x)-\varphi(y)\bigr)\DD x \DD y.
	\]
	Adding both sides of these relations we arrive at
	\[
		2\det\Omega(s,t) = \int_s^t \int_s^t \sin^2\bigl(\varphi(x)-\varphi(y)\bigr)\DD x \DD y.
	\]
\end{proof}

\noindent
In general, it is difficult to estimate the right-hand side of \eqref{X106} from below. 
One possibility is to exploit smoothness of $\varphi$.  
In the next lemma we give an estimate of this kind, and also include a bound from above.

\begin{lemma}\label{X35}
	Let $\varphi:[a,b] \to \bb R$ be absolutely continuous and not constant. 
	The following statements hold for $a\le s<t\le b$.
	\begin{Enumerate}
	\item
		If $\varphi'\in L^\infty(a,b)$, then
		\[
			\det\Omega(s,t) \ge \frac{1}{\|\varphi'\|_\infty^2}\cdot\frac{\bigl(\varphi(t)-\varphi(s)\bigr)^2
			-\sin^2\bigl(\varphi(t)-\varphi(s)\bigr)}{4}.
		\]
	\item
		If $\varphi$ is monotone and $\frac{1}{\varphi'} \in L^\infty(a,b)$, then
		\[
			\det\Omega(s,t) \le \Big\|\frac{1}{\varphi'}\Big\|_\infty^2
			\frac{\bigl(\varphi(t)-\varphi(s)\bigr)^2
			-\sin^2\bigl(\varphi(t)-\varphi(s))}{4}.
		\]
	\end{Enumerate}
\end{lemma}

\begin{proof}
	(i)
	Let us first consider the case when $\varphi'\in L^\infty(a,b)$.  Then, by \cref{X54}
	and with a straightforward evaluation of the integrals,
	\begin{align*}
		\det\Omega(s,t) &= \frac12 \int_s^t \int_s^t
		\sin^2\bigl(\varphi(x)-\varphi(y)\bigr)\DD x \DD y
		\\[1ex]
		&\ge \frac{1}{2\|\varphi'\|_\infty^2}\int_s^t \int_s^t
		\sin^2\bigl(\varphi(x)-\varphi(y)\bigr)\varphi'(x)\varphi'(y)\DD x \DD y
		\\[1ex]
		&= \frac{1}{2\|\varphi'\|_\infty^2}\cdot\frac{\bigl(\varphi(t)-\varphi(s)\bigr)^2
		-\sin^2\bigl(\varphi(t)-\varphi(s)\bigr)}{2}.
	\end{align*}

	(ii)
	Now assume that $\varphi$ is monotone and that $\frac{1}{\varphi'}\in L^\infty(a,b)$.
	The monotonicity of $\varphi$ implies that $\varphi'(x)\varphi'(y)\ge0$
	for a.e.\ $x,y\in (a,b)$, and hence
	\begin{align*}
		\det\Omega(s,t) &= \frac12 \int_s^t \int_s^t
		\sin^2\bigl(\varphi(x)-\varphi(y)\bigr)\DD x \DD y
		\\[1ex]
		&\le \frac 12\Big\|\frac{1}{\varphi'}\Big\|_\infty^2\int_s^t \int_s^t
		\sin^2\bigl(\varphi(x)-\varphi(y)\bigr)\varphi'(x)\varphi'(y)\DD x \DD y
		\\[1ex]
		&= \frac 12\Big\|\frac{1}{\varphi'}\Big\|_\infty^2\frac{\bigl(\varphi(t)-\varphi(s)\bigr)^2
		-\sin^2\bigl(\varphi(t)-\varphi(s)\bigr)}{2}.
	\end{align*}
\end{proof}

\noindent
We present a lemma which will enable us to invoke \cref{X67,X4}.

\begin{lemma}
\label{X38}
	Let $\varphi:[a,b] \to \bb R$ be a Lipschitz function with $\varphi(a)\ne\varphi(b)$, 
	and let $r\ge\max\{\|\varphi'\|_\infty,\frac{\|\varphi'\|_\infty}{|\varphi(a)-\varphi(b)|}\}$.  Set
	\begin{equation}
	\label{X117}
		k \DE \biggl\lfloor\big|\varphi(a)-\varphi(b)\big|
		\Bigl(\frac{r}{\|\varphi'\|_\infty}\Bigr)^{\frac 12}\biggr\rfloor.
	\end{equation}
	Then there exist points $s_j \in [a,b]$, $j=0,\ldots,k$, with
	\begin{align}
		& a = s_0 < s_1 < \cdots <s_{k-1} < s_k \le b, 
		\label{X55}
		\\
		&\forall j \in \{1,\ldots,k\}\DP \;\; \det \Omega(s_{j-1},s_j) > \frac{1}{16r^2}.
		\label{X171}
	\end{align}
\end{lemma}

\begin{proof}
	Let $f(x)\DE x^2-\sin^2x$.  A (tedious) calculation shows that $\frac{f(x)}{x^4}$ is decreasing on $(0,\infty)$.  Hence
	\begin{equation}\label{X41}
		x \le 1 \quad\Longrightarrow\quad f(x) \ge \frac{f(1)}{1^4}\cdot x^4 > \frac{x^4}{4}.
	\end{equation}
	Assume for definiteness that $\varphi(a)<\varphi(b)$. Set
	\[
		h \DE \Bigl(\frac{\|\varphi'\|_\infty}{r}\Bigr)^{\frac 12}
	\]
	and choose points $s_0,\ldots,s_k$ that satisfy \eqref{X55} such that 
	$\varphi(s_j)=\varphi(a)+jh$ for $j=0,\ldots,k$. By \cref{X35} and \eqref{X41} we have
	\begin{align*}
		\det\Omega(s_{j-1},s_j) \ge \frac{f\bigl(\varphi(s_j)-\varphi(s_{j-1})\bigr)}{4\|\varphi'\|_\infty^2} 
		= \frac{f(h)}{4\|\varphi'\|_\infty^2} > \frac{h^4}{16\|\varphi'\|_\infty^2} 
		= \frac{1}{16r^2},
	\end{align*}
	which yields \eqref{X171}.
\end{proof}

\noindent
As a consequence of the previous lemma we obtain a lower bound for $\kappa_H(r)$ when $\varphi$ is 
a Lipschitz function.

\begin{proposition}\label{X190}
	Let $\varphi:[a,b] \to \bb R$ be a non-constant Lipschitz function.
	Further, let $\kappa_H$ be as in \cref{X25} with $c=\frac{1}{16}$.
	Then
	\begin{equation}\label{X227}
		\kappa_H(r) \ge \frac{\osc_\varphi(a,b)}{2\|\varphi'\|_\infty^{\frac12}}\cdot r^{\frac12},
		\qquad r \ge \max\biggl\{\frac{\|\varphi'\|_\infty}{(\osc_\varphi(a,b))^2},\|\varphi'\|_\infty\biggr\}.
	\end{equation}
\end{proposition}

\begin{proof}
	There exist $c,d\in[a,b]$ with $c<d$ such that $|\varphi(c)-\varphi(d)|=\osc_\varphi(a,b)>0$
	where the inequality follows since $\varphi$ is not constant.  Let $r$ satisfy the 
	second inequality in \eqref{X227}.  Then $r\ge\|\varphi'|_{[c,d]}\|_\infty$ and
	\[
		|\varphi(c)-\varphi(d)|\biggl(\frac{r}{\|\varphi'|_{[c,d]}\|_\infty}\biggr)^{\frac12}
		\ge \osc_\varphi(a,b)\frac{r^{\frac12}}{\|\varphi'\|_\infty^{\frac12}} 
		\ge 1.
	\]
	Hence the assumptions of \cref{X38} applied to $\varphi|_{[c,d]}$ are satisfied.
	That lemma, together with \cref{X4}\,(ii), implies
	\[
		\kappa_H(r) \ge \biggl\lfloor\osc_\varphi(a,b)\frac{r^{\frac12}}{\|\varphi'\|_\infty^{\frac12}}\biggr\rfloor
		\ge \frac12 \osc_\varphi(a,b)\frac{r^{\frac12}}{\|\varphi'\|_\infty^{\frac12}},
	\]
	where we used the inequality $\lfloor x\rfloor\ge\frac{x}{2}$ for $x\ge1$.
\end{proof}

\begin{Remark}\label{X184}
	Let $\varphi$ be a non-constant Lipschitz function on $[a,b]$.
	By putting together some previous results (\cref{X190}, \cref{X36}\,(ii), \cref{X3,X1,X62}) 
	it follows that
	\[
		\log|w_{H,22}(b;ir)| \asymp \log\Bigl(\max_{|z|=r}\|W_H(b;z)\|\Bigr)
		\asymp \kappa_H(r) \asymp \int_a^b K_H(t;r)\DD t
		\asymp r^{\frac12}.
	\]
	In this connection we would like to mention the following result, which is proved 
	in \cite[Theorem~1.6\,(c)]{remling.scarbrough:2020}: if $\varphi$ is decreasing
	and has values in $\bigl[-\frac{\pi}{2},\frac{\pi}{2}\bigr]$, then $\rho_H<\frac12$ implies $\varphi'(t)=0$ for a.e.\ $t\in[a,b]$,
	where $\rho_H$ denotes the order of the monodromy matrix $W_H(b;\cdot)$ defined at the beginning of \cref{X138}.
\end{Remark}

\noindent
The next lemma, which is also based on \cref{X38}, provides a lower estimate for $\kappa_H(r)$ and an upper estimate
for the integral of $K_H$ by using a partition of the interval.

\begin{lemma}
\label{X44}
	Let $\varphi:(a,b] \to \bb R$ be locally Lipschitz continuous and let $H(t) \DE\xi_{\varphi(t)}\xi_{\varphi(t)}^T$. 
	Suppose that we are given
	\begin{Itemize}
	\item 
		a strictly decreasing sequence $(x_n)_{n=1}^\infty$ in $(a,b]$ with $x_1=b$;
	\item 
		a function $N: (R,\infty) \to \bb N\setminus\{1\}$ with $R>0$ such that, for all $r \in (R,\infty)$ and 
		$n \in \{2,\ldots,N(r)\}$,
	\end{Itemize}
	\begin{equation}
	\label{X103}
		0 < \big\|\varphi'|_{[x_n,x_{n-1}]}\big\|_\infty
		\le r\min\bigl\{1,|\varphi(x_{n-1})-\varphi(x_n)|^2\bigr\}.
	\end{equation}
	Set $r_0\DE\bigl(\frac{1}{16\det\Omega(a,b)}\bigr)^{\frac 12}$ and let $(\hat t,\hat s)$ be 
	a compatible pair for $H,r_0$ with constants $\frac{1}{16},\frac{1}{16}$, and let 
	$\kappa_H(r)$ be defined as in \cref{X25} with $c=\frac{1}{16}$.  Then $R\ge r_0$, and for all $r>R$,
	\begin{equation}\label{X176}
		\det\Omega(a,x_{N(r)-1}) > \frac{1}{16r^2}
	\end{equation}
	and
	\begin{align}
		\label{X174}
		\kappa_H(r) &\ge \frac 12 r^{\frac 12} \sum_{n=2}^{N(r)}  
		|\varphi(x_{n-1})-\varphi(x_n)|\,\big\|\varphi'|_{[x_n,x_{n-1}]}\big\|_\infty^{-\frac 12},
		\\[1ex]
		\label{X175}
		\int_{x_{N(r)-1}}^b \hspace{-1ex} K_H(t;r)\DD t &\le 8\sqrt{2}\cdot r^{\frac 12} 
		\sum_{n=2}^{N(r)-1} (x_{n-1}-x_{n+1})\big\|\varphi'|_{[x_{n+1},x_{n-1}]}\big\|_\infty^{\frac 12}.
	\end{align}
\end{lemma}

\begin{proof}
	Let $r>R$.
	For each $n \in \{2,\ldots,N(r)\}$ set
	\[
		k(n,r) \DE \biggl\lfloor\big|\varphi(x_{n-1})-\varphi(x_n)\big|
		\Bigl(\frac{r}{\|\varphi'|_{[x_n,x_{n-1}]}\|_\infty}\Bigr)^{\frac 12}\biggr\rfloor
		\ge \frac 12 \big|\varphi(x_{n-1})-\varphi(x_n)\big|
		\Bigl(\frac{r}{\|\varphi'|_{[x_n,x_{n-1}]}\|_\infty}\Bigr)^{\frac 12}.
	\]
	Then \cref{X38} yields $k(n,r)\ge1$ disjoint subintervals $[s_{j-1}^{(n,r)},s_j^{(n,r)}]$ of $[x_n,x_{n-1}]$ with 
	\begin{equation}\label{X172}
		\det\Omega\bigl(s_{j-1}^{(n,r)},s_j^{(n,r)}\bigr) > \frac{1}{16r^2}
	\end{equation}
	for all $j\in\{1,\ldots,k(n,r)\}$. 
	In total, we obtain
	\[
		\sum_{n=2}^{N(r)} k(n,r) \ge \frac 12 r^{\frac 12}\sum_{n=2}^{N(r)} 
		|\varphi(x_{n-1})-\varphi(x_n)|\,\big\|\varphi'|_{[x_n,x_{n-1}]}\big\|_\infty^{-\frac 12}
	\]
	disjoint subintervals of $[x_{N(r)},b]$ satisfying \eqref{X172}.
	Now \cref{X4}\,(ii) yields the inequality in \eqref{X174}.
	The inequalities in \eqref{X172} also imply the relation $R\ge r_0$ and \eqref{X176}.

	Let us now come to the proof of \eqref{X175}.
	Note that, by \eqref{X176}, we have $\hat t(r)<x_{N(r)-1}$ and hence 
	\begin{align*}
		K_H(t;r) = \frac{h_1(t)}{\omega_1(\hat s(t;r),t)}
	\end{align*}
	for $t\in[x_{N(r)-1},b]$.
	For each $n\in\{2,\ldots,N(r)-1\}$ set $I_n \DE [x_{n+1},x_{n-1}]$ and $H_n \DE H|_{I_n}$. 
	Let $(\hat t_n,\hat s_n)$ be the unique compatible pair for $H_n,\bigl(\frac{1}{16\det\Omega(x_{n+1},x_{n-1})}\bigr)^{\frac12}$ 
	with constants $\frac{1}{16},\frac{1}{16}$, which exists by \cref{X39}.  Since $\det\Omega(x_{n+1},x_n) > \frac{1}{16r^2}$
	by \eqref{X172}, it follows that $\hat t_n(r) < x_n$.  Hence $\hat s_n(t;r)$ is well defined for $t \in [x_n,x_{n-1}]$ and satisfies
	\[
		\det\Omega(\hat s_n(t;r),t) = \frac{1}{16r^2} = \det\Omega(\hat s(t;r),t), \qquad t \in [x_n,x_{n-1}].
	\]
	This implies $\hat s_n(t;r)=\hat s(t;r)$ for $t \in [x_n,x_{n-1}]$. 
	Note that, for $s,t \in I_n$, we have $|\varphi(s)-\varphi(t)| \le \|\varphi'|_{I_n}\|_\infty |s-t|$, 
	i.e.\ $\varphi|_{I_n}$ is Lipschitz continuous with constant $\|\varphi'|_{I_n}\|_\infty$. 
	Applying \cref{X36}\,(ii) to $H_n$ we obtain
	\begin{align*}
		&\int_{x_n}^{x_{n-1}} \mkern-2mu \frac{h_1(t)}{\omega_1(\hat s(t;r),t)} \DD t 
		= \int_{x_n}^{x_{n-1}} \mkern-2mu\frac{h_1(t)}{\omega_1(\hat s_n(t;r),t)} \DD t 
		\le \int_{\hat t_n(r)}^{x_{n-1}} \frac{h_1(t)}{\omega_1(\hat s_n(t;r),t)} \DD t 
		\\[1ex]
		&\le 8\sqrt{2}\,(x_{n-1}-x_{n+1})\big\|\varphi'|_{I_n}\big\|_\infty^{\frac 12}\cdot r^{\frac 12}.
	\end{align*}
	We now take the sum over $n\in\{2,\ldots,N(r)-1\}$ to complete the proof of \eqref{X175}.
\end{proof}

\bigskip

\noindent
We are now ready to prove \cref{X214}.

\begin{proof}[Proof of \cref{X214}]
	First note that we can choose, without loss of generality, $c=\frac{1}{16}$ in the definitions
	of $\kappa_H$ and $K_H$ by \cref{X219,X1}.

	Let $\ms f_{A/B}:[1,\infty)\to(0,\infty)$ be a continuous, strictly increasing, regularly varying function
	such that
	\[
		\ms f_{A/B}(x) \sim 
		\begin{cases}
			\frac{\ms f_2(x)}{\ms f_1(x)} & \text{in case (A)},
			\\[1ex]
			\frac{1}{\ms f_1(x)\ms f_2(x)} & \text{in cases (B) and (B$^+$)},
		\end{cases}
	\]
	as $x\to\infty$.
	Further, set $\ms g_{A/B}\DE\ms g_A$ in case (A) and $\ms g_{A/B}\DE\ms g_B$ in cases (B) and (B$^+$).
	Then $\ms g_{A/B}$ is an asymptotic inverse of $\ms f_{A/B}$.
	It follows from \eqref{X201} and the assumptions in cases (A), (B) and (B$^+$) that
	\[
		\frac{\big\|\varphi'|_{[x_n,x_{n-1}]}\big\|_\infty}{\min\bigl\{1,|\varphi(x_{n-1})-\varphi(x_n)|^2\bigr\}}
		\asymp \frac{\ms f_2(n)}{\ms f_1(n)}\cdot\frac{1}{\min\{1,\ms f_2(n)^2\}}
		\asymp \ms f_{A/B}(n).
	\]
	Hence there exists $\eta>0$ such that the left-hand side is bounded from above by $\frac{1}{\eta}\ms f_{A/B}(n)$
	for all $n\in\{2,3,\ldots\}$.
	Set $N(r)\DE\lfloor \ms f_{A/B}^{-1}(\eta r)\rfloor$ for $r>R$ where $R>0$ is chosen such that $\ms f_{A/B}^{-1}(\eta R)$
	is well defined and $N(R)\ge2$.
	Let $r>R$.  Then, for $n\in\{2,\ldots,N(r)\}$, we have $\ms f_{A/B}(n)\le\eta r$ and hence 
	the inequalities in \eqref{X103} are satisfied.
	Therefore we can apply \cref{X44}, which yields that \eqref{X174} and \eqref{X175} hold.
	We obtain from \eqref{X201} that
	\begin{align*}
		|\varphi(x_{n-1})-\varphi(x_n)|\,\big\|\varphi'|_{[x_n,x_{n-1}]}\big\|_\infty^{-\frac 12} &\asymp \sqrt{(\ms f_1\ms f_2)(n)},
		\\[1ex]
		(x_{n-1}-x_{n+1})\big\|\varphi'|_{[x_{n+1},x_{n-1}]}\big\|_\infty^{\frac 12} &\asymp \sqrt{(\ms f_1\ms f_2)(n)}.
	\end{align*}
	This implies that the two bounds in \eqref{X174} and \eqref{X175} coincide, and consequently, we arrive at
	\begin{equation}\label{X250}
		\int_{x_{N(r)-1}}^b K_H(t;r)\DD t \lesssim r^{\frac12}\int_1^{N(r)}\sqrt{(\ms f_1\ms f_2)(x)}\DD x
		\lesssim \kappa_H(r).
	\end{equation}
	It follows from \cref{X3} that
	\begin{align}
		\kappa_H(r) &\lesssim \log r + \int_a^b K_H(t;r)\DD t
		\\[1ex]
		&= \log r + \int_a^{\hat t(r)}K_H(t;r)\DD t + \int_{\hat t(r)}^{x_{N(r)-1}}K_H(t;r)\DD t + \int_{x_{N(r)-1}}^b K_H(t;r)\DD t.
		\label{X269}
	\end{align}
	Since the first integral on the right-hand side grows at most logarithmically by \cref{X62}
	and $\kappa_H(r)\gtrsim r^{\frac12}$ by \eqref{X250}, we obtain from \eqref{X269} and 
	the first inequality in \eqref{X250} that
	\[
		\kappa_H(r) \lesssim \int_a^b K_H(t;r)\DD t 
		\lesssim \int_{\hat t(r)}^{x_{N(r)-1}}K_H(t;r)\DD t + r^{\frac12}\int_1^{N(r)}\sqrt{(\ms f_1\ms f_2)(x)}\DD x.
	\]
	The fact that $\ms f_{A/B}^{-1}$ is regularly varying implies that 
	$N(r)\asymp \ms f_{A/B}^{-1}(\eta r)\asymp \ms f_{A/B}^{-1}(r)\asymp \ms g_{A/B}(r)$.
	With
	\[
		\ms b_{A/B,1}(r) \DE r^{\frac12}\int_1^{\ms g_{A/B}(r)}\sqrt{(\ms f_1\ms f_2)(x)}\DD x
	\]
	we therefore obtain
	\begin{align}
		\kappa_H(r) &\lesssim \int_a^b K_H(t;r)\DD t 
		\lesssim \int_{\hat t(r)}^{x_{N(r)-1}}K_H(t;r)\DD t + \ms b_{A/B,1}(r),
		\label{X251}
		\\[1ex]
		\ms b_{A/B,1}(r) &\lesssim \kappa_H(r),
		\label{X262}
	\end{align}
	where the second inequality follows from \eqref{X250}.
	Let us now estimate the first integral on the right-hand side of \eqref{X251}.
	\Cref{X82} implies that
	\begin{align}
		\int_{\hat t(r)}^{x_{N(r)-1}}K_H(t;r)\DD t 
		&\lesssim r(x_{N(r)-1}-a)
		= r\sum_{k=N(r)}^\infty(x_{k-1}-x_k) \asymp r\sum_{k=N(r)}^\infty \ms f_1(k)
		\nonumber\\[1ex]
		&\asymp r\int_{\ms g_{A/B}(r)}^\infty \ms f_1(x)\DD x
		\ED \ms b_{A/B,2}(r).
		\label{X263}
	\end{align}
	Let us now distinguish the different cases.
	
	\medskip
	
	\noindent
	Case (A).
	We write $\ms b_{A,i}$ for $\ms b_{A/B,i}$, $i=1,2$
	and use \cref{X95}, the relation $r\asymp\ms f_A(\ms g_A(r))$
	and \eqref{X262} to obtain
	\begin{align}
		\ms b_{A,2}(r) \;
		& \left\{\begin{array}{ll} 
			\asymp & \text{if} \ \rho_1<-1  \\[1ex]
			\gg & \text{if} \ \rho_1=-1
		\end{array}\right\} \;
		r\ms g_A(r)\ms f_1\bigl(\ms g_A(r)\bigr)
		\nonumber\\[1ex]
		&\asymp \ms g_A(r)\ms f_2\bigl(\ms g_A(r)\bigr)
		\asymp r^{\frac12}\ms g_A(r)\sqrt{(\ms f_1\ms f_2)(\ms g_A(r))}
		\nonumber\\[1ex]
		& \left\{\begin{array}{ll}
			\asymp & \text{if} \ \frac{\rho_1+\rho_2}{2}>-1 \\[1ex]
			\ll & \text{if} \ \frac{\rho_1+\rho_2}{2}\le-1
		\end{array}\right\} \;
		\ms b_{A,1}(r)
		\lesssim \kappa_H(r).
		\label{X256}
	\end{align}
	If $\rho_1<-1$, then $\ms b_{A,2}(r)\lesssim\ms b_{A,1}(r)\lesssim\kappa_H(r)$,
	and therefore \eqref{X260} and \eqref{X275} follow from \eqref{X251}, \eqref{X263} and \eqref{X256}.
	
	Now assume that $\rho_1=-1$.  Then $\frac{\rho_1+\rho_2}{2}\ge-\frac12>-1$ and hence
	$\ms b_{A,2}(r)\gg\ms b_{A,1}(r)$.
	This, together with \eqref{X251} and \eqref{X263} and another application of \eqref{X256}
	yields \eqref{X261}.
	The statements about the indices of the bounds are easy to check.

	\medskip

	\noindent
	Cases (B) and (B$^+$).
	With $\ms b_{B,i}\equiv\ms b_{A/B,i}$, the inequalities in \eqref{X264} follow directly 
	from \eqref{X262}, \eqref{X251} and \eqref{X263}.
	Using \cref{X95,X206} one can easily show the relations in \eqref{X265}--\eqref{X268}.

	Let us first consider the case when $\int_1^\infty\sqrt{(\ms f_1\ms f_2)(x)}\DD x<\infty$.
	First note that $\ms b_{B,1}(r)\asymp r^{\frac12}$.
	Since $r\asymp\ms f_B(\ms g_B(r))$, we obtain from \cref{X95} that
	\begin{equation}\label{X274}
		\frac{\ms b_{B,2}(r)}{r^{\frac12}} \;
		\left\{\begin{array}{ll} 
			\asymp & \text{if} \ \rho_1<-1  \\[1ex]
			\gg & \text{if} \ \rho_1=-1
		\end{array}\right\} \;
		r^{\frac12}\ms g_B(r)\ms f_1\bigl(\ms g_B(r)\bigr)
		\asymp \ms g_B(r)\sqrt{\frac{\ms f_1(\ms g_B(r))}{\ms f_2(\ms g_B(r))}}.
	\end{equation}
	If $\rho_1<-1$, then
	\[
		\ms b_{B,2}(r) \; \left\{\begin{array}{c} \lesssim \\ \gg \end{array}\right\} \; r^{\frac12}
		\quad\Leftrightarrow\quad
		x\sqrt{\frac{\ms f_1(x)}{\ms f_2(x)}} \; \left\{\begin{array}{c} \lesssim \\ \gg \end{array}\right\} \; 1,
	\]
	which implies \eqref{X270}.
	If $\rho_1=-1$, then $x^2\ms f_1(x)\gg\ms f_2(x)$ and hence $\ms b_{B,2}(r)\gg r^{\frac12}$ again
	by \eqref{X274}.

	Next let us consider the case when $\int_1^\infty\sqrt{(\ms f_1\ms f_2)(x)}\DD x=\infty$.
	First, assume that $\rho_2<0$.  Then $\Ind\ms b_{B,2}>\Ind\ms b_{B,1}$ and 
	hence $\ms b_{B,2}(r)\gg\ms b_{B,1}(r)$.
	Assume now that $\rho_2=0$ and $\rho_1>-2$.  Then, by \cref{X95},
	\begin{align*}
		\ms b_{B,2}(r) &= r\int_{\ms g_B(r)}^\infty \ms f_1(t)\DD t \;
		\left\{\begin{array}{ll} 
			\asymp & \text{if} \ \rho_1<-1  \\[1ex]
			\gg & \text{if} \ \rho_1=-1
		\end{array}\right\} \;
		r\ms g_B(r)\ms f_1\bigl(\ms g_B(r)\bigr)
		\\[1ex]
		&\asymp \frac{\ms g_B(r)}{\ms f_2(\ms g_B(r))} \;
		\left\{\begin{array}{ll} 
			\gg & \text{if} \ \ms f_2(x)\ll 1  \\[1ex]
			\gtrsim & \text{otherwise}
		\end{array}\right\} \;
		\ms g_B(r)
		\\[1ex]
		&\asymp \ms g_B(r)r^{\frac12}\sqrt{(\ms f_1\ms f_2)(\ms g_B(r))} \asymp \ms b_{B,1}(r),
	\end{align*}
	which proves the last claim in case (B).

	\medskip

	\noindent
	Case (B$^+$).
	Let us now assume, in addition, that $\varphi\in L^\infty(a,b)$ and that \eqref{X271} holds.
	In this situation we estimate the integral in \eqref{X251} with the help of \cref{X36}\,(i)
	instead of \cref{X82}.  Using part of the estimate in \eqref{X263} we obtain
	\begin{align*}
		\int_{\hat t(r)}^{x_{N(r)-1}}K_H(t;r)\DD t 
		&\lesssim r(x_{N(r)-1}-a)\cdot\osc_\varphi(a,x_{N(r)-1})
		\lesssim \ms b_{B,2}(r)\ms f_3\bigl(N(r)-1\bigr)
		\\[1ex]
		&\asymp \ms b_{B,2}(r)\ms f_3\bigl(\ms g_B(r)\bigr),
	\end{align*}
	which, together with \eqref{X251}, yields the upper estimate in \eqref{X272}.
	The lower estimate is already proved in the general case (B).
	
	Let us prove the last statement in (B$^+$).  Assume that $\rho_1<-1$.
	First, we consider the case when $\ms f_3=\ms f_2$.
	It follows from \eqref{X267} and \eqref{X265} that
	\[
		\ms b_{B,2}(r)\ms f_3\bigl(\ms g_B(r)\bigr) \asymp \frac{\ms g_B(r)}{\ms f_2(\ms g_B(r))}\cdot\ms f_2\bigl(\ms g_B(r)\bigr)
		= \ms g_B(r) \lesssim \ms b_{B,1}(r),
	\]
	which, together with \eqref{X272}, proves \eqref{X273}.
	Finally, assume that $\ms f_3(x)=\frac{1}{x}\sqrt{\frac{\ms f_2(x)}{\ms f_1(x)}}$.
	It follows from \eqref{X267} that
	\[
		\ms b_{B,2}(r)\ms f_3\bigl(\ms g_B(r)\bigr) 
		\asymp \frac{\ms g_B(r)}{\ms f_2(\ms g_B(r))}\cdot\frac{1}{\ms g_B(r)}\sqrt{\frac{\ms f_2(\ms g_B(r))}{\ms f_1(\ms g_B(r))}}
		= \frac{1}{\sqrt{(\ms f_1\ms f_2)(\ms g_B(r))}\,}
		\asymp r^{\frac12} \lesssim \ms b_{B,1}(r),
	\]
	which again proves \eqref{X273}.
\end{proof}

\subsection{Examples}
\label{X151}

In this subsection we consider some examples to illustrate \cref{X214}.
The first example, presented in \cref{X42}, occurs already in \cite[Example~5.6]{pruckner.woracek:sinqA}, 
where an upper bound for the order of the monodromy matrix is proved.
We determine the exact growth of the monodromy matrix (up to multiplicative constants)
demonstrating the applicability of \cref{X214}.
Further, we study two examples in connection with some boundary cases in \cref{X214}.

\begin{theorem}\label{X42}
	Let $(\gamma,\beta)\in\bb R^2 \setminus\{(0,0)\}$ and let $\varphi\DF(0,1]\to\bb R$ be the chirp signal 
	\[
		\varphi(t) \DE t^\gamma\sin\Bigl(\frac{1}{t^\beta}\Bigr).
	\]
	Consider the Hamiltonian $H(t)\DE\xi_{\varphi(t)}\xi_{\varphi(t)}^T$ where $\xi_\phi=(\cos\phi,\sin\phi)^T$, 
	let $W_H(1;z)$ be its monodromy matrix, and let $\kappa_H(r)$ be as in \cref{X25}.
	Then
	\begin{equation}\label{X197}
		\log\Bigl(\max_{|z|=r}\|W_H(1;z)\|\Bigr) \asymp \log|w_{H,22}(1;ir)| \asymp \kappa_H(r) \asymp
		\begin{cases}
			r^{\frac12} \CAS \beta<\gamma+1,
			\\[0.5ex]
			r^{\frac12}\log r \CAS \beta=\gamma+1,
			\\[0.5ex]
			r^{\hat\rho} \CAS \beta>\gamma+1,
		\end{cases}
	\end{equation}
	as $r\to\infty$, where
	\begin{equation}\label{X195}
		\hat\rho =
		\begin{cases}
			\frac{\beta}{\beta+\gamma+1} \CAS \gamma\ge0,
			\\[1.5ex]
			\frac{\beta-\gamma}{\beta-\gamma+1} \CAS \gamma<0.
		\end{cases}
	\end{equation}
\end{theorem}

\begin{figure}[ht]
\begin{center}
	\begin{tikzpicture}[x=1.2pt,y=1.2pt,scale=0.8,font=\fontsize{8}{8}]
		\draw[dashed,->] (16,0)--(232,0);
		\draw[dashed] (130,-94)--(130,20);
		\draw[->] (130,20)--(130,117);
		\draw     (20,-90)--(130,20);
		\draw (130,20) -- (220,110) node[midway, sloped, below]{ ${\displaystyle r^{\frac 12}\log r}$};
		
		\draw (122,114) node {${\displaystyle \beta}$};
		\draw (232,-8) node {${\displaystyle \gamma}$};
		
		\draw (160,-30) node {\large $r^{\frac 12}$};
		\fill[pattern=dots, opacity=0.1] (20,-90)--(220,110)--(220,-90)--(20,-90);
		
		\draw (80,40) node {\Large $r^{\frac{\beta-\gamma}{\beta-\gamma+1}}$};
		\fill[pattern=crosshatch dots, opacity=0.25] (20,-90)--(130,20)--(130,110)--(20,110)--(20,-90);
		
		\draw (155,80) node {\Large $r^{\frac{\beta}{\beta+\gamma+1}}$};
		\fill[pattern=crosshatch, opacity=0.15] (130,20)--(130,110)--(220,110)--(130,20);
	\end{tikzpicture}
	\caption{The asymptotic behaviour of $\log|w_{H,22}(t;ir)|$ and $\kappa_H(r)$ depending
	on the parameters $\beta$ and $\gamma$.}\label{X196}
\end{center}
\end{figure}
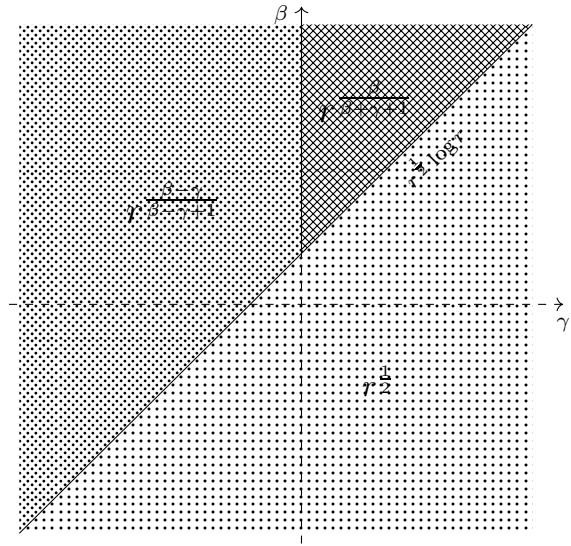

\begin{remark}\label{X194}
	The exponent $\hat\rho$ in \eqref{X195} is increasing in $\beta$ and decreasing in $\gamma$.
	Moreover, it converges to $\frac12$ when $(\gamma,\beta)$ approaches the line $\beta=\gamma+1$
	and it converges to $1$ when $\beta\to\infty$ or $\gamma\to-\infty$.
	See also \cref{X196}.
\end{remark}

\begin{proof}[Proof of \cref{X42}]
	We distinguish several cases.
	\begin{Elist}
	\item
		Assume that $\beta\le0$ and $\gamma\ge\beta$. \\
		In this case the function $\varphi$ is of bounded variation.  Hence \cref{X36}\,(iii) and \cref{X62},
		together with \cref{X1}, imply that
		\[
			\log|w_{H,22}(1;ir)| \asymp \int_0^1 K_H(t;r)\DD t \lesssim r^{\frac12}.
		\]
		It now follows from \cref{X45} that
		\[
			\log\Bigl(\max_{|z|=r}\|W_H(1;z)\|\Bigr) \lesssim r^{\frac12}.
		\]
		Since $\varphi$ is a Lipschitz function on $\bigl[\frac12,1\bigr]$, we obtain from \cref{X3,X19,X190} that 
		\[
			\int_0^1 K_H(t;r)\DD t \gtrsim \kappa_H(r) \ge \kappa_{H|_{[\frac12,1]}}(r) \gtrsim r^{\frac12},
		\]
		which proves the claim in this case.
		
		In the remaining cases we use \cref{X214}.
		
	\item
		Assume that $\gamma<\beta\le0$. \\
		In this situation we have $\varphi(t)\sim t^\alpha$ as $t\to0$ with $\alpha\DE\gamma-\beta<0$.
		Set $x_n\DE n^{\frac{1}{\alpha}}$, $n\in\bb N$.  A simple computation shows that
		\begin{align*}
			& x_{n-1}-x_n \sim \frac{1}{|\alpha|}n^{\frac{1}{\alpha}-1},
			\\[1ex]
			& |\varphi(x_{n-1})-\varphi(x_n)| \sim 1,
			\\[1ex]
			& \big\|\varphi'|_{[x_n,x_{n-1}]}\big\|_\infty \sim |\varphi'(x_n)| \sim |\alpha|n^{\frac{\alpha-1}{\alpha}}
		\end{align*}
		as $n\to\infty$.  
		Hence \eqref{X201} is satisfied with
		\[
			\ms f_1(x) = x^{\frac{1}{\alpha}-1}, \qquad \ms f_2(x) = 1, \qquad 
			\rho_1=\frac{1}{\alpha}-1, \qquad \rho_2=0.
		\]
		and we can choose $\ms g_A(r)=r^{\frac{\alpha}{\alpha-1}}$.
		It follows from \cref{X214}, case (A), and \cref{X45} that
		\begin{align*}
			\log\Bigl(\max_{|z|=r}\|W_H(1;z)\|\Bigr) \asymp \log|w_{H,22}(1;ir)| \asymp \kappa_H(r) \asymp
			\begin{cases}
				r^{\frac12} \CAS \alpha>-1,
				\\[1ex]
				\displaystyle r^{\frac12}\int_1^{r^{\frac12}}\frac{1}{x}\DD x \CAS \alpha=-1,
				\\[3ex]
				r^{\frac{\alpha}{\alpha-1}} \CAS \alpha<-1,
			\end{cases}
		\end{align*}
		which proves \eqref{X197} in this case.
		
	\item
		Assume that $\beta>0$. \\
		In this case $\varphi$ oscillates.  Let us choose
		\[
			x_1 \DE 1; \qquad x_n \DE \bigl(\pi(n+\tfrac12)\bigr)^{-\frac{1}{\beta}}, \quad n\ge2.
		\]
		One can easily show that
		\begin{align*}
			& x_{n-1}-x_n \sim \frac{1}{\beta}\pi^{-\frac{1}{\beta}}n^{-\frac{\beta+1}{\beta}},
			\\[1ex]
			& |\varphi(x_{n-1})-\varphi(x_n)| 
			= \bigl(\pi\bigl(n-\tfrac12\bigr)\bigr)^{-\frac{\gamma}{\beta}} + \bigl(\pi\bigl(n+\tfrac12\bigr)\bigr)^{-\frac{\gamma}{\beta}}
			\sim 2\pi^{-\frac{\gamma}{\beta}} n^{-\frac{\gamma}{\beta}}
		\end{align*}
		as $n\to\infty$.
		Since
		\[
			|\varphi'(t)| = \big|t^{\gamma-\beta-1}\bigl(-\beta\cos(t^{-\beta})+\gamma t^\beta\sin(t^{-\beta})\bigr)\big|
			\le (\beta+|\gamma|)t^{\gamma-\beta-1}
		\]
		for $t\in(0,1]$, we have
		\[
			\big\|\varphi'|_{[x_n,x_{n-1}]}\big\|_\infty \le \max_{t\in[x_n,x_{n-1}]}(\beta+|\gamma|)t^{\gamma-\beta-1}
			\sim (\beta+|\gamma|)(\pi n)^{-\frac{\gamma-\beta-1}{\beta}}.
		\]
		Hence \eqref{X201} is satisfied with
		\[
			\ms f_1(x) = x^{-\frac{\beta+1}{\beta}}, \qquad \ms f_2(x) = x^{-\frac{\gamma}{\beta}}, \qquad 
			\rho_1=-\frac{\beta+1}{\beta}, \qquad \rho_2=-\frac{\gamma}{\beta}.
		\]
		Let us first consider the case when $\gamma\le0$.  Then we have $\ms f_2(x)\gtrsim1$, 
		and we can choose $\ms g_A(r)=r^{\frac{\beta}{\beta-\gamma+1}}$.  \Cref{X214}, case (A), and \cref{X45} yield
		\begin{align*}
			\log\Bigl(\max_{|z|=r}\|W_H(1;z)\|\Bigr) \asymp \log|w_{H,22}(1;ir)| \asymp \kappa_H(r) \asymp
			\begin{cases}
				r^{\frac12} \CAS \beta<\gamma+1,
				\\[1ex]
				\displaystyle r^{\frac12}\int_1^{r^{\frac{\beta}{2}}}\frac{1}{x}\DD x \CAS \beta=\gamma+1,
				\\[3ex]
				\ms g_A(r)\ms f_2\bigl(\ms g_A(r)\bigr) \CAS \beta>\gamma+1;
			\end{cases}
		\end{align*}
		in the third case we can rewrite the expression on the right-hand side:
		\[
			\ms g_A(r)\ms f_2\bigl(\ms g_A(r)\bigr) = \bigl(r^{\frac{\beta}{\beta-\gamma+1}}\bigr)^{1-\frac{\gamma}{\beta}}
			= r^{\frac{\beta-\gamma}{\beta-\gamma+1}}.
		\]
		Let us finally assume that $\gamma>0$.  Then $\ms f_2(x)\lesssim1$, $\varphi\in L^\infty(0,1)$,
		\[
			\osc_\varphi(0,x_n) \le 2x_n^\gamma \sim 2(\pi n)^{-\frac{\gamma}{\beta}} \asymp \ms f_2(n),
		\]
		and therefore \eqref{X271} with $\ms f_3=\ms f_2$ is also satisfied.
		With $\ms g_B(r)=r^{\frac{\beta}{\beta+\gamma+1}}$ we obtain from \cref{X214}, case (B$^+$), and \cref{X45} that
		\begin{align*}
			\log\Bigl(\max_{|z|=r}\|W_H(1;z)\|\Bigr) \asymp \log|w_{H,22}(1;ir)| \asymp \kappa_H(r) \asymp
			\begin{cases}
				r^{\frac12} \CAS \beta<\gamma+1,
				\\[1ex]
				\displaystyle r^{\frac12}\int_1^{r^{\frac12}}\frac{1}{x}\DD x \CAS \beta=\gamma+1,
				\\[3ex]
				\ms g_B(r) \CAS \beta>\gamma+1,
			\end{cases}
		\end{align*}
		which proves \eqref{X197} also in this case.
	\end{Elist}
\end{proof}

\noindent
The next example shows that in case~(B) with $\int_1^\infty\sqrt{(\ms f_1\ms f_2)(t)}\DD t=\infty$,
any of the two functions $\ms b_{B,1}$, $\ms b_{B,2}$ can dominate the other one.

\begin{Example}\label{X177}
	Assume that \eqref{X201} holds with $\ms f_1(x)=\frac{1}{x^2}$ and $\ms f_2(x)=(\log x)^{-\alpha}$
	with $\alpha\in[0,2)$.
	We are then in case (B) in \cref{X214} with $\int_1^\infty\sqrt{(\ms f_1\ms f_2)(t)}\DD t=\infty$.
	It is easy to see that one can choose $\ms g_B(r)=2^{\frac{\alpha}{2}}r^{\frac12}(\log r)^{-\frac{\alpha}{2}}$
	and that
	\[
		\ms b_{B,1}(r) \asymp r^{\frac12}(\log r)^{1-\frac{\alpha}{2}}, \qquad
		\ms b_{B,2}(r) \asymp r^{\frac12}(\log r)^{\frac{\alpha}{2}}.
	\]
	Hence
	\[
		\ms b_{B,2}(r) \;
		\left\{\begin{array}{ll} 
			\ll & \text{if} \ \alpha\in[0,1)  \\[1ex]
			\asymp & \text{if} \ \alpha=1 \\[1ex]
			\gg & \text{if} \ \alpha\in(1,2)
		\end{array}\right\} \;
		\ms b_{B,1}(r).
	\]
	This shows that in the case when $\rho_1=-2$, $\rho_2=0$, everything can happen.
\end{Example}

\medskip

\noindent
Let us consider an example where $H$ oscillates at the left endpoint with exponential speed.
The order of $W_H(1;\cdot)$ is $1$ in this situation, and we can determine the growth of $\log|w_{H,22}(1;ir)|$ 
up to a logarithmic factor.

\begin{Example}\label{X180}
	Let $H(t)=\xi_{\varphi(t)}\xi_{\varphi(t)}^T$ where
	\[
		\varphi(t) = e^{t^{-\gamma}}, \qquad t\in(0,1],
	\]
	with $\gamma>0$.  Choose a strictly decreasing sequence $(x_n)_{n=1}^\infty$ 
	such that $x_1=1$ and $x_n=(\log n)^{-\frac{1}{\gamma}}$ for large enough $n$.
	Then \eqref{X201} is satisfied with
	\[
		\ms f_1(x) = \frac{1}{x(\log x)^{\frac{1}{\gamma}+1}}, \qquad \ms f_2(x) = 1, \qquad
		\rho_1=-1, \qquad \rho_2=0.
	\]
	We can choose $\ms g_A$ such that $\ms g_A(r)=\frac{r}{(\log r)^{\frac{1}{\gamma}+1}}$ for large enough $r$.
	It follows from \cref{X214}, case~(A), that
	\[
		\frac{r}{(\log r)^{\frac{1}{\gamma}+1}} \lesssim \kappa_H(r)
		\lesssim \int_a^b K_H(t;r)\DD t
		\lesssim \frac{r}{(\log r)^{\frac{1}{\gamma}}}.
	\]
\end{Example}

\subsection{An inverse construction}
\label{X212}

In \cref{X213} we construct a Hamiltonian whose monodromy matrix has a 
prescribed growth for any given regularly varying function with index
in $(\frac12,1)$ and some functions with index $\frac12$.  
This substantially improves \cite[Theorem~6.1]{pruckner.woracek:sinqA}
where the growth is prescribed only up to some slowly varying function.
Recall from \cref{X74} that it is no restriction in the following
theorem that $\ms f$ or $v$ are smoothly varying functions instead of just regularly
varying functions.

\begin{theorem}\label{X213}
	Let either $\ms f:[1,\infty)\to(0,\infty)$ be a strictly increasing, smoothly varying function
	with $\Ind\ms f\in(\frac12,1)$ or let 
	\begin{equation}\label{X215}
		\ms f(r) = r^{\frac12}\int_1^r v(s)\DD s, \qquad r\in[1,\infty),
	\end{equation}
	where $v$ is a smoothly varying function with index $-1$ such that $\int_1^\infty v(s)\DD s=\infty$
	\textup{(}in which case $\Ind\ms f=\frac12$\textup{)}.
	In the first case set
	\[
		\ms g(x) \DE \frac{x}{\ms f^{-1}(x)}, \qquad x\in[\ms f(1),\infty);
	\]
	in the second case set
	\[
		\ms h(r) \DE \int_1^r s^{\frac12}v(s)\DD s \qquad\text{and}\qquad 
		\ms g(x) \DE \frac{x}{\ms h^{-1}(x)}, \qquad x\in[1,\infty).
	\]
	In both cases there exists $\hat x>0$ such that $\ms g$ is strictly decreasing on $[\hat x,\infty)$.
	Define the function
	\[
		\varphi(t) \DE \ms g^{-1}(t), \qquad t\in(0,b],
	\]
	with $b\DE\ms g(\hat x)$ and set $H(t)\DE\xi_{\varphi(t)}\xi_{\varphi(t)}^T$ where $\xi_\phi$ 
	is defined in \eqref{X104}.  Then
	\[
		\log\Bigl(\max_{|z|=r}\|W_H(b;z)\|\Bigr) \asymp \log|w_{H,22}(b;ir)| \asymp \ms f(r).
	\]
\end{theorem}

\begin{Remark}\label{X218}
\rule{0ex}{1ex}
\begin{Enumerate}
\item
	The function $\varphi$ is smoothly varying at $0$ with index $\frac{\rho}{\rho-1}$ where $\rho\DE\Ind\ms f$.
\item
	\Cref{X213} is bound to growth functions $\ms f$ that satisfy $\ms f(r)\gg r^{\frac12}$,
	which is an intrinsic restriction.
	On the rough scale of order it is complemented by \cite[Theorem~1.6]{remling.scarbrough:2020}:
	there, for every $\rho\in[0,\frac12]$, a Hamiltonian is constructed whose monodromy matrix has order $\rho$.
	In that paper no assertions about more detailed growth are made.
\end{Enumerate}
\end{Remark}

\medskip

\begin{proof}[Proof of \cref{X213}]
	Set $\rho\DE\Ind\ms f$.
	The function $\ms g$ is smoothly varying with $\Ind \ms g=1-\frac{1}{\rho}\in[-1,0)$, 
	and hence $\ms g$ is strictly decreasing on some interval $[\hat x,\infty)$ by \eqref{X222}.

	Choose $x_n$, $n\in\bb N$, strictly decreasing such that $x_1=b$ and $x_n=\ms g(n)$ for large enough $n$.
	We have $x_n\to0$ as $n\to\infty$ since $\Ind\ms g<0$.
	Moreover, $\varphi(x_n)=n$ for large enough $n$, which yields
	\[
		|\varphi(x_{n-1})-\varphi(x_n)| = 1.
	\]
	The smooth variation of $\ms g$ and \eqref{X222} imply that
	\[
		x_{n-1}-x_n = \ms g(n-1)-\ms g(n) \sim -\ms g'(n) \sim \Bigl(\frac{1}{\rho}-1\Bigr)\cdot\frac{\ms g(n)}{n}.
	\]
	Set $\ms f_1(x)\DE\frac{\ms g(x)}{x}$ and $\ms f_2(x)\DE1$.  Then the first two relations in \eqref{X201} are satisfied.
	Since $\varphi$ is smoothly varying at $0$ with $\Ind\varphi=\frac{\rho}{\rho-1}<0$, 
	we have $\varphi''(t)>0$ for small enough $t$ by \eqref{X222}.
	Hence \cref{X240}\,(ii) implies that also the third relation in \eqref{X201} holds.
	With the notation from \cref{X214} we have $\rho_1=\Ind\ms g-1=-\frac{1}{\rho}$, $\rho_2=0$ and
	\[
		\frac{\ms f_2(x)}{\ms f_1(x)} = \frac{x}{\ms g(x)} = 
		\begin{cases}
			\ms f^{-1}(x) \CAS \rho\in(\frac12,1),
			\\[1ex]
			\ms h^{-1}(x) \CAS \rho=\frac12.
		\end{cases}
	\]
	We can choose $\ms g_A=\ms f$ when $\rho\in(\frac12,1)$ and $\ms g_A=\ms h$ when $\rho=\frac12$.
	Hence \cref{X214}, case (A), and \cref{X45} imply that
	\[
		\log\Bigl(\max_{|z|=r}\|W_H(b;z)\|\Bigr) \asymp \log|w_{H,22}(b;ir)|
		\asymp \ms b_{A,1}(r)
	\]
	with the notation from \cref{X214}.
	If $\rho\in(\frac12,1)$, then $\ms b_{A,1}(r)\asymp \ms g_A(r) = \ms f(r)$ by \eqref{X275}.
	Let us now assume that $\rho=\frac12$, i.e.\ $\ms f$ is given by \eqref{X215}.
	Then the substitution $x=\ms h(s)$ yields
	\begin{align*}
		\ms b_{A,1}(r)
		&\asymp r^{\frac12}\int_1^{\ms g_A(r)}\sqrt{\ms f_1(x)}\DD x
		= r^{\frac12}\int_1^{\ms h(r)}\frac{1}{\sqrt{\ms h^{-1}(x)}\,}\DD x
		\\[1ex]
		&= r^{\frac12}\int_{\ms h^{-1}(1)}^r s^{-\frac12}\ms h'(s)\DD s
		= r^{\frac12}\int_{\ms h^{-1}(1)}^r v(s)\DD s
		\sim \ms f(r),
	\end{align*}
	which finishes the proof also in this case.
\end{proof}

\section{Additions and examples}
\label{X135}

In this section we derive a comparison results (\cref{X150}) and study 
two examples (\cref{X152,X153}); the first example shows the sharpness of the estimate in \cref{X36}\,(ii),
the second example proves the sharpness of the upper bound in \cref{X3}.

\subsection{A comparison result}
\label{X150}

Comparison theorems for Sturm--Liouville and Schr\"odinger equations have a long history
going back to C.F.~Sturm.
For canonical systems, for instance, spectral properties of a Hamiltonian $H$ have been compared 
with those of the diagonal Hamiltonian that is obtained by setting the off-diagonal entries equal to zero;
see \cite[Theorems~3.4 and 3.6]{romanov.woracek:ideal} and \cite[Theorem~1.1]{remling.scarbrough:2020a}.
In \cref{X8} we show that the function $K_H(t;r)$, or at least the integral over $K_H(t;r)$---and with it the 
quantities on the respective left-hand sides of \eqref{X90} and \eqref{X65}--- depends monotonically on the values of $H$.

\begin{theorem}
\label{X8}
	Let $H,\wt H\in\bb H_{a,b}$ be definite, assume that $\int_a^b h_1(t)\DD t<\infty$ and that either
	both $H,\wt H$ are in the limit circle case or both are in the limit point case.
	Further, let $r_0\ge 0$, assume that we have compatible pairs for $H,r_0$ and $\wt H,r_0$ with constants $c_-,c_+$, 
	and denote by $K_H(t;r)$ and $K_{\wt H}(t;r)$ the corresponding functions in \eqref{X15}.
	\begin{Enumerate}
	\item
		Let $\gamma_-,\gamma_+>0$ and assume that 
		\begin{alignat}{2}
			\label{X223}
			& \det\Omega_{\wt H}(s,t) \le \gamma_+^2\det\Omega_H(s,t),\qquad && s,t\in\dom H \text{ with } s<t,
			\\[1ex]
			\label{X224}
			& \omega_{\wt H,2}(a,t)\le \gamma_+\omega_{H,2}(a,t),\qquad && t\in(a,b)\text{ a.e.}, 
			\\[1ex]
			\label{X225}
			& \gamma_-h_1(t) \le \tilde h_1(t) \le \gamma_+h_1(t),\qquad && t\in(a,b)\text{ a.e.}
		\end{alignat}
		Then 
		\begin{equation}
		\label{X164}
			K_{\wt H}(t;r) \lesssim K_H(t;r), \qquad r>r_0,\;t\in\dom H.
		\end{equation}
		The constant implicit in the relation `\/$\lesssim$' depends on $\gamma_-,\gamma_+,c_-,c_+$ but not on $H,\wt H$.
	\item
		Assume that $H$ and $\wt H$ are in the limit circle case,
		let $\gamma_+>0$ and assume that \eqref{X223} holds.  Then
		\[
			\int_a^b K_{\wt H}(t;r)\DD t \lesssim \log r\cdot \int_a^b K_H(t;r)\DD t, \qquad r>r_0',
		\]
		with $r_0'\DE\max\{r_0,\gamma_+^2r_0,e\}$ and where the implicit constant in $\lesssim$ depends
		only on $\gamma_+,c_-,c_+$,$\int_a^b\tr H(t)\DD t$ and $\int_a^b\tr \wt H(t)\DD t$.
	\end{Enumerate}
\end{theorem}

\begin{Remark}\label{X226}
	Note that if $\wt H(t)\le H(t)$ for a.e.\ $t\in\dom H$, then \eqref{X223}, \eqref{X224}
	and the second inequality in \eqref{X225} are satisfied.
\end{Remark}

\noindent
In order to prove \cref{X8}\,(i), a monotonicity property for the Weyl coefficients of Hamiltonians 
in the limit point case is needed. 
The lemma below is a generalisation of \cite[Corollary~2.4]{reiffenstein:imq}. 

\begin{lemma}\label{X50}
	Let $H,\wt H\in\bb H_{a,b}$ be in the limit point case, and assume that $(a,b)$ is 
	not indivisible of type $0$ for $H$ or $\wt H$. 
	Further, let $\gamma_-,\gamma_+>0$ and assume that 
	\begin{alignat}{2}
		\label{X162}
		& \det\Omega_{\wt H}(a,t) \le \gamma_+^2\det\Omega_H(a,t),\qquad && t\in[a,b),
		\\[1ex]
		\label{X163}
		& \omega_{\wt H,2}(a,t) \ge \gamma_-\omega_{H,2}(a,t),\qquad && t\in[a,b).
	\end{alignat}
	Then 
	\begin{equation}
	\label{X59}
		\Im q_{\wt H}(ir) \lesssim \Im q_H(ir),\qquad r\in(0,\infty).
	\end{equation}
	The constant implicit in the relation `\/$\lesssim$' depends on $\gamma_-,\gamma_+$, but not on $H,\wt H$.
\end{lemma}

\begin{proof}
	If $\wt H$ is not definite, then $q_{\wt H}$ is a real constant, and \eqref{X59} holds with any constant in
	`\/$\lesssim$'.  If $H$ is not definite, then 
	\[
		\det\Omega_{\wt H}(a,b) \le \gamma_+^2\cdot\det\Omega_H(a,b) = 0,
	\]
	and hence $\wt H$ is also not definite.  Again \eqref{X59} holds with any constant in `\/$\lesssim$'.  

	Assume now that $H,\wt H$ are both definite. 
	Let $\hat t_H(r)$ and $\hat t_{\wt H}(r)$ be compatible functions for $H,0$ and $\wt H,0$, respectively, 
	with constants $1,1$, which exist and are unique by \cref{X39}.  It follows from assumption \eqref{X162} that 
	\[
		\det\Omega_H\bigl(a,\hat t_{\wt H}(r)\bigr)
		\ge \frac{1}{\gamma_+^2}\det\Omega_{\wt H}\bigl(a,\hat t_{\wt H}(r)\bigr)
		= \frac{1}{(\gamma_+r)^2} = \det\Omega_H\bigl(a,\hat t_H(\gamma_+r)\bigr),
	\]
	and hence $\hat t_{\wt H}(r) \ge \hat t_H(\gamma_+r)$.  From this, \eqref{X163} and \cref{X49} we obtain 
	\begin{align*}
		\Im q_{\wt H}(ir)
		&\asymp \frac{1}{r\omega_{\wt H,2}\bigl(a,\hat t_{\wt H}(r)\bigr)}
		\le \frac{1}{\gamma_-r\omega_{H,2}\bigl(a,\hat t_{\wt H}(r)\bigr)}
		\\[1ex]
		&\le \frac{1}{\gamma_-r\omega_{H,2}\bigl(a,\hat t_H(\gamma_+r)\bigr)}
		\asymp \Im q_H(i\gamma_+r) \asymp \Im q_H(ir),
	\end{align*}
	where the last asymptotic relation follows from the inequalities
	\[
		\min\biggl\{\gamma_+,\frac{1}{\gamma_+}\biggr\} \le \frac{\Im q_H(i\gamma_+r)}{\Im q_H(ir)}
		\le \max\biggl\{\gamma_+,\frac{1}{\gamma_+}\biggr\},
	\]
	which can easily be shown with the help of the Herglotz integral representation of $q_H$. 
\end{proof}

\begin{proof}[Proof of \cref{X8}]
	(i)
	For $t=a$ both sides of \eqref{X164} vanish.
	For the rest of the proof let $t\in\dom H\setminus\{a\}$.
	Define the Hamiltonians $H_{(t)}$ and $\wt H_{(t)}$ as in \eqref{X27} corresponding to $H$ and $\wt H$.
	Combined with \eqref{X112} and \eqref{X182}, our assumption implies that, for all $s \in (a,\infty)$,
	\[
		\det\Omega_{\wt H_{(t)}}(a,s) \le \max \{\gamma_+,\gamma_+^2 \}\cdot\det\Omega_{H_{(t)}}(a,s), \qquad 
		\omega_{\wt H_{(t)},2}(a,s) \ge \min\{\gamma_-,1\}\omega_{H_{(t)},2}(a,s).
	\]
	An application of \cref{X50} with $H_{(t)}$ and $\wt H_{(t)}$ yields 
	\[
		\Im q_{\wt H_{(t)}}(ir) \lesssim \Im q_{H_{(t)}}(ir),
		\qquad r>0,
	\]
	where the constant implicit in `\/$\lesssim$' depends on $\gamma_+,\gamma_-$, but not on $H,\wt H,t$. 
	It follows from this, \cref{X52}, \eqref{X63} and the assumption $\tilde h_1 \le \gamma_+h_1$ that,
	\[
		K_{\wt H}(t;r) \asymp r\bigl(\Im q_{\tilde H_{(t)}}(ir)\bigr)\tilde h_1(t)
		\lesssim r\bigl(\Im q_{H_{(t)}}(ir)\bigr)h_1(t)
		\asymp K_H(t;r)
	\]
	for $r>r_0$.
	
	(ii)
	Let $\hat t_H$ be the given compatible function for $H,r_0$ with constants $c_-,c_+$.
	It follows from the definition of a compatible function and the monotonicity of $\det\Omega_H$ that, 
	for all $r>r_0$,
	\[
		\det\Omega_H(a,b) \ge \det\Omega_H(a,\hat t_H(r)) \ge \frac{c_-}{r^2},
	\]
	and hence
	\[
		r_0 \ge \Bigl(\frac{c_-}{\det\Omega_H(a,b)}\Bigr)^{\frac12},
	\]
	and similarly for $\wt H$ instead of $H$.

	In the following we indicate the dependence on constants in the notation.
	For $c>0$, we denote by $K_{H,c}$ the kernel that corresponds to the unique compatible pair for $H,r_0$
	with constants $c,c$ from \cref{X39}.  Moreover, let $\kappa_{H,c}$ be the function from \cref{X25}
	corresponding to $H$ and the constant $c$.
	
	Now let $\sigma_j^{(r)}$, $j=0,\ldots,\kappa_{H,c_-}(r)$, be the sequence as in \cref{X25} corresponding
	to $H$ and the constant $c_-$.  It follows from the definition of $\sigma_j^{(r)}$ and assumption \eqref{X223} that
	\[
		\det\Omega_{\wt H}\bigl(\sigma_{j-1}^{(r)},\sigma_j^{(r)}\bigr)
		\le \gamma_+^2\det\Omega_H\bigl(\sigma_{j-1}^{(r)},\sigma_j^{(r)}\bigr)
		\le \frac{\gamma_+^2c_-}{r^2},
	\]
	for $j\in\{1,\ldots,\kappa_{H,c_-}(r)\}$; note that the second inequality is an equality for $j<\kappa_{H,c_-}(r)$.
	Hence \cref{X4}\,(i) implies that 
	\[
		\kappa_{\wt H,\gamma_+^2c_-}(r) \le \kappa_{H,c_-}(r), \qquad r>0.
	\]
	We now obtain from \cref{X1,X3} that, for $r>r_0'$,
	\begin{align*}
		\int_a^b K_{\wt H}(t;r)\DD t 
		&\asymp \log|w_{\wt H,22}(t;ir)| 
		\asymp \int_a^b K_{\wt H,\gamma_+^2c_-}(t;r)\DD t
		\lesssim \kappa_{\wt H,\gamma_+^2c_-}(r)\cdot\log r
		\\[1ex]
		&\le \kappa_{H,c_-}(r)\cdot\log r
		\lesssim \log r\cdot\int_a^b K_{H,c_-}(t;r)\DD t
		\asymp \log r\cdot\int_a^b K_H(t;r)\DD t,
	\end{align*}
	where the implicit constants depend only on $c_-,c_+,\gamma_+$, $\int_a^b\tr H(t)\DD t$ and $\int_a^b\tr\wt H(t)\DD t$.
\end{proof}

\subsection{Rotation given by the Weierstra{\ss} function}
\label{X152}

In \cref{X36}\,(ii) we considered Hamiltonians of the form $H(t)\DE\xi_{\varphi(t)}\xi_{\varphi(t)}^T$, 
where the function $\varphi$ is assumed to be H\"older continuous with exponent $\nu$. 
This result, together with \cref{X45}, implies that
\begin{equation}\label{X125}
	\log\Bigl(\max_{|z|=r}\|W_H(b;z)\|\Bigr) \lesssim r^{\frac{1}{1+\nu}},
\end{equation}
which is precisely the statement of \cite[Corollary~5.3]{pruckner.woracek:sinqA}. 
It is known that this upper bound is sharp up to a multiplicative factor that is slowly varying
(see \cite[Theorem~6.1 and Example~6.2]{pruckner.woracek:sinqA}), 
but it was unclear whether there is a function $\varphi$ for which `$\asymp$' holds in \eqref{X125}. 
We show in \cref{X127} that there exist such functions $\varphi$ for every $\nu\in(0,1)$. 
In this example the lower bound in \cref{X3} is attained (up to a multiplicative constant). 
It is given by the Weierstra{\ss} function, which is defined as
\begin{equation}\label{X160}
	\varphi(t) \DE \sum_{n=0}^\infty \alpha^n\cos(\beta^n\pi t), \qquad t\in\bb R,
\end{equation}
where $0<\alpha<1$, $\alpha\beta \ge 1$.  The Weierstra{\ss} function is an early example (but not the earliest) of a continuous 
but nowhere differentiable function.  For $\alpha\beta>1$, it is H\"older continuous with exponent
\begin{equation}\label{X126}
	\nu \DE -\frac{\log\alpha}{\log\beta};
\end{equation}
see, e.g.\ \cite[Theorem~1.33]{hardy:1916}. 

\begin{theorem}
\label{X127}
	Let $\alpha\in(0,1)$, and let $\beta$ be an even integer such that $\alpha\beta > 1+\frac{\pi}{2}$.
	Consider the Hamiltonian $H(t)\DE\xi_{\varphi(t)}\xi_{\varphi(t)}^T$, $t \in (0,1)$, 
	where $\varphi$ is the Weierstra{\ss} function in \eqref{X160} with parameters $\alpha,\beta$
	and where $\xi_\phi$ is defined in \eqref{X104}.  Then
	\begin{equation}\label{X168}
		\log\Bigl(\max_{|z|=r}\|W_H(1;z)\|\Bigr) \asymp r^{\frac{1}{1+\nu}} 
		\asymp\kappa_H(r),
		\qquad r>1,
	\end{equation}
	and 
	\[
		0 < \limsup_{r\to\infty}\frac{n_H(r)}{r^{\frac{1}{1+\nu}}} < \infty,
	\]
	where $\nu$, given by \eqref{X126}, is the exponent of H\"older continuity of $\varphi$.
\end{theorem}

\noindent
Let us point out that we can start with any $\nu\in(0,1)$ and construct $\alpha,\beta$ satisfying the 
assumption of \cref{X127} and such that \eqref{X126} holds.  Hence, we have, for any $\nu\in(0,1)$, 
an example of a H\"older continuous function with exponent $\nu$ for which `$\asymp$' holds in \eqref{X125} 
and $\kappa_H(r)$ gives the exact growth.
We expect \cref{X127} to hold for arbitrary $\alpha,\beta$ with $0<\alpha<1$ and $\alpha\beta>1$, but cannot give a proof due to 
some technical complications.

\begin{proof}
	Due to \eqref{X125}, it suffices to estimate $\log |w_{H,22}(1;ir)|$ from below. 
	We claim that, for $m\in\bb N$ and $k\in\{0,\ldots,\beta^m-1\}$,
	\begin{equation}\label{X128}
		\det\Omega\bigl(k\beta^{-m},(k+1)\beta^{-m}\bigr) \ge c\cdot\Bigl(\frac{\alpha}{\beta}\Bigr)^{2m}
	\end{equation}
	with some $c>0$ independent of $m$ and $k$, which we prove below.  

	Let $\kappa_H(r)$ be as in \cref{X25} for $c$ as in \eqref{X128};
	note that, by \cref{X219}, the asymptotic properties of $\kappa_H$ are independent of $c$.
	Set $r_m \DE \bigl(\frac{\beta}{\alpha}\bigr)^m$, $m\in\bb N_0$.
	For $m\in\bb N$, \eqref{X128} implies that \eqref{X167} is satisfied with $s_j=j\beta^{-m}$, $j\in\{0,\ldots,\beta^m\}$
	and $r=r_m$.  Hence \cref{X4}\,(ii) yields that $\kappa_H(r_m)\ge\beta^m$.
	The latter inequality is also true for $m=0$.
	Define $m(r)\DE\bigl\lfloor\frac{\log r}{\log\beta-\log\alpha}\bigr\rfloor$ for $r>1$.
	Then, for arbitrary $r>1$, we have $r\in[r_{m(r)},r_{m(r)+1})$ and thus
	\[
		\kappa_H(r) \ge \kappa_H(r_{m(r)}) \ge \beta^{m(r)} > \beta^{\frac{\log r}{\log\beta-\log\alpha}-1}
		= \frac{1}{\beta}r^{\frac{\log\beta}{\log\beta-\log\alpha}} = \frac{1}{\beta}r^{\frac{1}{1+\nu}}.
	\]
	Now it follows from \cref{X1,X3} that
	\[
		\log|w_{H,22}(t;ir)| \asymp \int_0^1 K_H(t;r)\DD t
		\ge \log 2\cdot\kappa_H(r)-\BigO(\log r)
		\gtrsim r^{\frac{1}{1+\nu}},
	\]
	which proves \eqref{X168}.
	
	We now come to the proof of the claim \eqref{X128}.  Since $\beta$ is an even integer, we have,
	for $k\in\{0,\ldots,\beta^m\}$,
	\begin{align*}
		\varphi(k\beta^{-m}) &= \sum_{n=0}^{m-1}\alpha^n\cos\bigl(\pi k\beta^{n-m}\bigr) + \alpha^m\cos(\pi k)
		+ \sum_{n=m+1}^\infty \alpha^n\cos\bigl(\pi k\beta^{n-m}\bigr)
		\\[1ex]
		&= \sum_{n=0}^{m-1}\alpha^n\cos\bigl(\pi k\beta^{n-m}\bigr) + (-1)^k\alpha^m + \frac{\alpha^{m+1}}{1-\alpha},
	\end{align*}
	and hence, for $k\in\{0,\ldots,\beta^m-1\}$,
	\begin{align}
		\nonumber
		& \big|\varphi\bigl((k+1)\beta^{-m}\bigr)-\varphi\bigl(k\beta^{-m}\bigr)+2(-1)^k\alpha^m\big|
		= \Bigg|\sum_{n=0}^{m-1} \alpha^n\Bigl[\cos\bigl(\pi(k+1)\beta^{n-m}\bigr)-\cos\bigl(\pi k\beta^{n-m}\bigr)\Bigr]\Bigg|
		\\[1ex]
		\nonumber
		&= 2\,\Bigg|\sum_{n=0}^{m-1} \alpha^n\sin\Bigl(\frac{\pi}{2}(2k+1)\beta^{n-m}\Bigr)\cdot\sin\Bigl(\frac{\pi}{2}\beta^{n-m}\Bigr)\Bigg|
		\le 2\sum_{n=0}^{m-1} \alpha^n\sin\Bigl(\frac{\pi}{2}\beta^{n-m}\Bigr)
		\\[1ex]
		\label{X165}
		&\le \pi\sum_{n=0}^{m-1}\alpha^n\beta^{n-m} 
		= \pi\beta^{-m}\frac{(\alpha\beta)^m-1}{\alpha\beta-1}
		\le \frac{\pi}{\alpha\beta-1}\alpha^m
		= (2-\varepsilon)\alpha^m
	\end{align}
	with $\varepsilon\in(0,2)$, where the latter follows from the assumption $\alpha\beta>1+\frac{\pi}{2}$.
	The inequality in \eqref{X165} implies that
	\begin{equation}\label{X129}
		\big|\varphi\bigl((k+1)\beta^{-m}\bigr)-\varphi\bigl(k\beta^{-m}\bigr)\big| \ge \varepsilon \alpha^m.
	\end{equation}
	By \cite[Theorem~1.33]{hardy:1916} the function $\varphi$ is H\"older continuous 
	with exponent $\nu$, i.e.\ $|\varphi(x)-\varphi(y)|\le\eta|x-y|^\nu$ with some $\eta>0$.
	Let $\delta\DE\bigl(\frac{\varepsilon}{3\eta}\bigr)^{\frac{1}{\nu}}$.
	For every $x\in[k\beta^{-m},(k+\delta)\beta^{-m}]$ we have
	\[
		\big|\varphi(x)-\varphi(k\beta^{-m})\big| \le \eta(\delta\beta^{-m})^\nu
		= \frac{\varepsilon}{3}\beta^{-m\nu} = \frac{\varepsilon}{3}\alpha^m
	\]
	and similarly $|\varphi(y)-\varphi((k+1)\beta^{-m})|\le\frac{\varepsilon}{3}\alpha^m$
	for every $y\in[(k+1-\delta)\beta^{-m},(k+1)\beta^{-m}]$,
	which, together with \eqref{X129}, implies that
	\[
		|\varphi(x)-\varphi(y)| \ge \frac{\varepsilon}{3}\alpha^m.
	\]
	Now we obtain from \cref{X54} that
	\begin{align*}
		& \det\Omega\bigl(k\beta^{-m},(k+1)\beta^{-m}\bigr) 
		= \frac 12\int\limits_{k\beta^{-m}}^{(k+1)\beta^{-m}}\mkern8mu\int\limits_{k\beta^{-m}}^{(k+1)\beta^{-m}}\mkern-12mu 
		\sin^2\bigl(\varphi(x)-\varphi(y)\bigr)\DD x\DD y 
		\\[1ex]
		&\ge \frac 12\int\limits_{(k+1-\delta)\beta^{-m}}^{(k+1)\beta^{-m}}\int\limits_{k\beta^{-m}}^{(k+\delta)\beta^{-m}}\mkern-12mu 
		\sin^2\Bigl(\frac{\varepsilon}{3}\alpha^m\Bigr)\DD x\DD y 
		= \frac12 \delta^2 \beta^{-2m}\sin^2\Bigl(\frac{\varepsilon}{3}\alpha^m\Bigr)
		\asymp \Bigl(\frac{\alpha}{\beta}\Bigr)^{2m},
	\end{align*}
	which finishes the proof of the claim \eqref{X128}.

	The assertion for the counting function $n_H(r)$ follows from \cref{X1,X7}.
\end{proof}

\subsection{Sharpness of the upper bound}
\label{X153}

We already saw several instances where the lower bound in \cref{X3} is attained (up to a multiplicative constant). 
Now we give an example where 	the upper bound is attained.

\begin{proposition}
\label{X5}
	Let $H$ be the Hamiltonian defined on the interval $[0,1]$ a.e.\ by 
	\[
		H(t)
		\DE
		\begin{cases}
			\smmatrix 1000 \CAS e^{-(2n+1)} < t \le e^{-2n},
			\\[1.5ex]
			\smmatrix 0001 \CAS e^{-(2n+2)} < t \le e^{-(2n+1)}
		\end{cases}
		\qquad\text{for }n\in\bb N_0.
	\]
	Let $\kappa_H(r)$ and $K_H(t;r)$ be as in \cref{X3} with $c=1$.  Then 
	\[
		\kappa_H(r) \asymp \log r, \qquad \int_0^1 K_H(t;r)\DD t \asymp (\log r)^2,
	\]
	for $r>r_0=\frac{1}{\sqrt{\det\Omega(0,1)}\,}$.
\end{proposition}

\begin{proof}
	As a first step we compute $\det\Omega$ on initial sections of the interval $[0,1]$.
	Let $n\in\bb N$; then 
	\begin{align*}
		\det\Omega(0,e^{-2n})
		&= \sum_{l=n}^\infty\bigl[e^{-2l}-e^{-(2l+1)}\bigr]\cdot\sum_{l=n}^\infty\bigl[e^{-(2l+1)}-e^{-(2l+2)}\bigr]
		\\
		&= \frac{1}{e}\Bigl(1-\frac 1e\Bigr)^2\biggl(\sum_{l=n}^\infty e^{-2l}\biggr)^2
		= \alpha e^{-4n},
	\end{align*}
	where $\alpha\DE\frac{e}{(e+1)^2}$.  This formula yields that 
	\[
		\det\Omega(0,e^{-2n}) \le \frac{1}{r^2}
		\quad\Longleftrightarrow\quad
		n \ge \frac 12\log r+\frac 14\log\alpha.
	\]
	For sufficiently large $r$, the interval 
	\[
		\Bigl[\frac 12\log r+\frac 14\log\alpha,\,\log r\Bigr]
	\]
	contains at least one integer.  Let $n(r)$ be such an integer and set 
	\[
		s_0 \DE 0; \quad 
		s_l \DE e^{-2n(r)+l-1} \;\;\text{for} \ l=1,\ldots,2n(r)+1.
	\]
	Then
	\[
		\det\Omega(s_0,s_1) \le \frac{1}{r^2}; \qquad 
		\det\Omega(s_{l-1},s_l) = 0 \quad\text{for} \ l\in\{2,\ldots,2n(r)+1\}.
	\]
	Hence we can apply \cref{X4}\,(i) to obtain $\kappa_H(r) \le 2n(r)+1 \le 2\log r+1$, 
	which, in turn, implies $\int_0^1 K_H(t;r)\DD t\lesssim(\log r)^2$ by \cref{X3}.

	The reverse inequalities are obtained by explicit computation of the integral of $K_H(t;r)$ over certain subintervals. 
	Let $t\in[e^{-(2n+1)},e^{-2n}]$.  Then we clearly have $\hat s(t;r) \le e^{-(2n+1)}$.  Furthermore, 
	\[
		\det\Omega\bigl(e^{-(2n+2)},t\bigr) = (e-1)e^{-(2n+2)}\bigl(t-e^{-(2n+1)}\bigr),
	\]
	and hence
	\[
		\hat s(t;r) \ge e^{-(2n+2)}
		\quad\Longleftrightarrow\quad
		\det\Omega\bigl(e^{-(2n+2)},t\bigr) \ge \frac{1}{r^2}
		\quad\Longleftrightarrow\quad
		t \ge e^{-(2n+1)}+\frac{e^{2n+2}}{r^2(e-1)}.
	\]
	The interval 
	\[
		\Bigl[e^{-(2n+1)}+\frac{e^{2n+2}}{r^2(e-1)},\,e^{-2n}\Bigr]
	\]
	is non-empty if and only if 
	\[
		n \le \frac 12\log r-\frac 14\log\beta,
	\]
	where $\beta\DE\frac{e^3}{(e-1)^2}$.  For such $n$ we can estimate
	\begin{align*}
		\int\limits_{e^{-(2n+1)}}^{e^{-2n}} K_H(t;r)\DD t
		&\ge \int\limits_{e^{-(2n+1)}+\frac{e^{2n+2}}{r^2(e-1)}}^{e^{-2n}}\frac{h_1(t)}{\omega_1(\hat s(t;r),t)}\DD t
		\\[1ex]
		&= \int\limits_{e^{-(2n+1)}+\frac{e^{2n+2}}{r^2(e-1)}}^{e^{-2n}}\frac{1}{t-e^{-(2n+1)}}\DD t
		= 2\log r-4n-\log\beta.
	\end{align*}
	Set $N\DE\bigl\lfloor\frac 14\log r\bigr\rfloor$.  For large enough $r$ we have $N\le\frac12\log r-\frac14\log\beta$
	and hence
	\begin{align*}
		\int_0^1 K_H(t;r)\DD t
		&\ge \sum_{n=0}^N\int_{e^{-(2n+1)}}^{e^{-2n}}K_H(t;r)\DD t
		\ge \sum_{n=0}^N\bigl(2\log r-4n-\log\beta\bigr)
		\\[1ex]
		&= 2(N+1)\Bigl[\log r-N-\frac12\log\beta\Bigr]
		\\[1ex]
		&\ge \frac12\log r\cdot\Bigl[\frac34\log r-\frac12\log\beta\Bigr]
		\gtrsim (\log r)^2.
	\end{align*}
	Together with \cref{X3}, this implies that $\kappa_H(r)\gtrsim\log r$. 
\end{proof}

\begin{remark}\label{X156}
	We observe that for the Hamiltonians in \cref{X5}, where the upper bound in \cref{X3} is attained, the monodromy matrix
	has order zero. Contrasting this, for the Hamiltonians in \cref{X214,X127}, where the lower bound is attained, the
	monodromy matrix has positive order. This phenomenon occurs in all examples we know.
\end{remark}


{\footnotesize
\begin{flushleft}
	M.~Langer \\
	Department of Mathematics and Statistics \\
	University of Strathclyde \\
	26 Richmond Street \\
	Glasgow G1 1XH \\
	UNITED KINGDOM \\
	email: \texttt{m.langer@strath.ac.uk} \\[5mm]
\end{flushleft}
\begin{flushleft}
	J.~Reiffenstein \\
	Department of Mathematics \\
	Stockholms universitet \\
	106 91 Stockholm \\
	SWEDEN \\
	email: \texttt{jakob.reiffenstein@math.su.se} \\[5mm]
\end{flushleft}
\begin{flushleft}
	H.~Woracek\\
	Institute for Analysis and Scientific Computing\\
	Vienna University of Technology\\
	Wiedner Hauptstra{\ss}e\ 8--10/101\\
	1040 Wien\\
	AUSTRIA\\
	email: \texttt{harald.woracek@tuwien.ac.at}\\[5mm]
\end{flushleft}
}

\end{document}